%% file: inverse_swe_siam_revised.tex
\documentclass[onefignum,onetabnum]{siamart171218}

\input{inverse_swe_shared}

\ifpdf
\hypersetup{
  pdftitle={An Example Article},
  pdfauthor={D. Doe, P. T. Frank, and J. E. Smith}
}
\fi

\usepackage{soul}

\newcommand{\hly}[1]{{#1}}
\newcommand{\hlc}[1]{{#1}}

\newcommand{\hmy}[1]{{#1}}

\begin{document}

\maketitle

\input{abstract}

\begin{keywords}
  Inverse Problems, Adjoint Methods, Discontinuous Galerkin Methods, 
  Hyberbolic Balance Laws, Shallow Water Equations, 
  $L^1$ and $H^1$ Regularization
\end{keywords}

\begin{AMS}
  35L60, 35L65, 35Q35, 65M30, 65M32, 65M60, 65K10 
\end{AMS}

\input{inverse_swe_body_revised}


\bibliographystyle{siamplain}
\bibliography{references_inverseSWE}
\end{document}

%% file: inverse_swe_shared.tex
\usepackage[ruled,vlined,algo2e]{algorithm2e}
\usepackage{bm}
\usepackage{booktabs}
\usepackage{graphicx}
\usepackage{subcaption}
\usepackage[utf8]{inputenc}
\usepackage[english]{babel}
\usepackage{xcolor}
\usepackage{rotating}

\usepackage{lipsum}
\usepackage{amsfonts}
\usepackage{epstopdf}
\ifpdf
  \DeclareGraphicsExtensions{.eps,.pdf,.png,.jpg}
\else
  \DeclareGraphicsExtensions{.eps}
\fi

\newsiamremark{remark}{Remark}
\newsiamremark{hypothesis}{Hypothesis}
\crefname{hypothesis}{Hypothesis}{Hypotheses}
\newsiamthm{claim}{Claim}

\newcommand{\U}{\bm{U}}
\newcommand{\tU}{\bm{\tilde{U}}}

\newcommand{\sig}{\bm{\sigma}}
\newcommand{\F}{\bm{F}}
\newcommand{\FU}{\bm{F}(\bm{U})}
\newcommand{\G}{\bm{G}}
\newcommand{\bS}{\bm{S}}

\newcommand{\B}{\bm{B}}
\newcommand{\tB}{{\tilde{B}}}

\newcommand{\tp}{\tilde{p}}

\newcommand{\A}{\bm{A}}
\newcommand{\C}{\bm{C}}
\newcommand{\J}{\bm{J}}

\newcommand{\dJ}{\bm{\delta J}}
\newcommand{\R}{\bm{R}}

\newcommand{\bv}{\bm{v}}

\newcounter{subsecval}[section]
\renewcommand{\thesubsecval}{\thesubsection}
\newcounter{subsecvalnext}[subsecval]
\renewcommand{\thesubsecvalnext}{\thesubsecval\alph{subsecvalnext}}

\headers{Recovery of a Bottom Topography Function}{J. Britton, Y. T. Chow, W. Chen, and Y. Xing}

\title{Recovery of a Time-Dependent Bottom Topography Function from the Shallow Water Equations via an Adjoint Approach\thanks{Submitted to the editors \today{}.
\funding{The work of Y. Xing is partially supported by the NSF grant DMS-1753581. W. Chen is partially supported by the NSF grant DMS-1853701 and DMS-1762063 through the joint NSF DMS/NIH NIGMS Initiative to Support Research at the Interface of the Biological and Mathematical Sciences.}}}

\author{Jolene Britton\thanks{Department of Mathematics, University of California Riverside, Riverside, CA 92521, USA 
  (\email{jhout001@ucr.edu}, \email{yattinc@ucr.edu}, \email{weitaoc@ucr.edu}).}
\and Yat Tin Chow\footnotemark[2]
\and Weitao Chen\footnotemark[2]
\and Yulong Xing\thanks{Department of Mathematics, The Ohio State University, Columbus, OH 43210, USA
(\email{xing.205@osu.edu})}}

\usepackage{amsopn}

\makeatletter
\newcommand*{\addFileDependency}[1]{
  \typeout{(#1)}
  \@addtofilelist{#1}
  \IfFileExists{#1}{}{\typeout{No file #1.}}
}
\makeatother


%% file: abstract.tex
\begin{abstract}
    We develop an adjoint approach for recovering the topographical function included in the source term of one-dimensional hyperbolic balance laws. We focus on a specific system, namely the shallow water equations, in an effort to recover the riverbed topography. The novelty of this work is the ability to robustly recover the bottom topography using only noisy boundary data from one measurement event and the inclusion of two regularization terms in the iterative update scheme. The adjoint scheme is determined from a linearization of the forward system and is used to compute the gradient of a cost function. The bottom topography function is recovered through an iterative process given by a three-operator splitting method which allows the feasibility to include two regularization terms. Numerous numerical tests demonstrate the robustness of the method regardless of the choice of initial guess and in the presence of discontinuities in the solution of the forward problem. 
\end{abstract}

%% file: inverse_swe_body_revised.tex
\section{Introduction} \label{sec_intro}

Various phenomena arising frequently in natural, engineering and socio-economical applications can be modeled by hyperbolic conservation and balance laws. Examples of conservation laws include models for traffic flow \cite{Garavello2006TrafficModel}, fluid dynamics \cite{Eiseman1980ConservationSurvey} and supply chains \cite{Brunnermeier2002InteroperabilityChain}. Moreover, conservation laws with source terms, also known as balance laws, are used in different models, e.g. the gas pipeline flow \cite{Gugat2011GasStabilization}, shallow water flow \cite{Xing2017NumericalEquations, Kurganov2018Finite-volumeEquations}, gas dynamics under gravitational field \cite{Xing2013HighFields}, and blood flow through arteries \cite{Young1808XIII.Blood}. A variety of theoretical studies have been conducted to understand the underlying mathematical structure and a wide range of numerical approaches dedicated to solving for the state variables have been developed in the past few decades. On the other hand, optimization, control, and recovery of the system parameters is a problem of great interest due to its high practical value. In this work, we derive an inverse problem algorithm for a specific system of hyperbolic balance laws, in which a time-dependent topographical profile is recovered based on the measurement from the boundary.

Various optimization and control methods have been developed for conservation and balance laws such as backstepping \cite{Vazquez2011LocalBackstepping,Glass2007OnEquation}, Lyapunov-based \cite{Vazquez2011LocalBackstepping}, derivative-free  \cite{Larson2019Derivative-freeMethods}, and optimal control methods \cite{Jacquet2006OptimalLaws}. In this work we employ an adjoint approach, which is often praised for its efficiency. Its computational cost (of each iteration) is comparable to that of solving a  partial differential equation (PDE) once, instead of depending on the number of control variables or design parameters as in other optimization methods.

Adjoint methods were introduced for optimal control problems in 1971 by Lions \cite{Lions1971OptimalEquations} in the context of shape optimization of aerodynamic bodies. A few years later the methods were extended to fluid dynamics by Pironneau \cite{Pironneau1974OnMechanics}. Jameson popularized the techniques for potential flow and the Euler equations \cite{Jameson1988AerodynamicTheory,Jameson1995OptimumTheory}. The methods have also been applied to biological systems in the search of parameter identification \cite{Raffard2008AnSignaling}. Estimation of open water states \cite{Castaings2006AutomaticModeling} and traffic states on the freeways \cite{Jacquet2005TrafficWaves} have also been accomplished via the adjoint optimization method.

\hly{The adjoint method for computing the gradient of a cost function can be accomplished in a variety of manners. The discrete adjoint approach} \cite{Gugat2005OptimalNetworks,Kotsialos2004NonlinearMetering,Giles2000AnDesign} \hly{begins with discretizing the forward system of PDEs while the continuous adjoint approach} \cite{Gugat2005OptimalNetworks,Jacquet2005OPTIMALMETHODS,Moin1994FeedbackTurbulence,Reuther1996AerodynamicFormulation} \hly{begins with the continuous forward system. Both cases result in discrete adjoint equations. The continuous approach allows one to obtain an explicit gradient formulation} \cite{Reilly2015Adjoint-BasedMetering} \hly{while automatic differentiation techniques can be applied to ease discrete approach solvers} \cite{Muller2005OnDifferentiation,Giering1998RecipesConstruction}.
With regards to scalar conservation laws, Holden et al. \cite{Holden2014OnLaws} developed a reconstruction procedure for the coefficient inverse problem in which a spatially dependent coefficient of the flux term is recovered. B{\"{u}}rger et al. \cite{Burger2009NumericalSedimentation} solved the inverse scalar conservation law modelling sedimentation numerically by assuming a variational form of the problem.  The flux function of a scalar conservation law was reconstructed using the information from the shock that forms in the work by Kang and Tanuma \cite{Kang2005InverseLaws}.
In a more general setting for balance laws, Montecinos et al. \cite{Montecinos2005AProblems} derived a unified scheme for solving the forward and adjoint problems simultaneously. Methodology for the scalar \hly{Burgers'} equation was presented by Lellouche et al. \cite{Lellouche1994BoundaryApproach} in which the authors aimed to find the best approximation for the measured data by means of boundary control and an adjoint approach. Ferlauto \cite{Ferlauto2015AEquations} obtained optimal geometric shapes for aerodynamic bodies by solving an inverse problem for the three-dimensional incompressible Euler equations.

Numerical computation of the optimization problems for conservation laws have been studied extensively due to the theoretical and numerical challenges that arise. As the exact solution of conservation laws often contains discontinuities, one challenge in the related optimization problems is that non-negligible numerical errors may occur in capturing the discontinuities. Some of these difficulties are mitigated via introduction of the Lax-Friedrichs schemes \cite{Giles2010ConvergenceExtensions} or relaxation methods \cite{Banda2012AdjointLaws}, for instance.
Convergence analyses have been provided for optimization problems in the aforementioned works.
The fact that many conservation laws are nonlinear presents another challenge because this can lead to non-convex formulations of the optimization problem. One method to tackle this difficulty is to use linear programming methods once the discretization scheme is `relaxed' \cite{Gomes2006OptimalModel,Ziliaskopoulos2000LinearProblem}. This allows for a global optimum to be found and reduces computational cost, but the linearization may not be a good physical representation of the original system \cite{Reilly2015Adjoint-BasedMetering}. To maintain the nonlinearity of the system, a more expensive method, such as gradient descent, can be used but may not ensure a global optimum is achieved. 

In this paper, we focus on the inverse problem that arises in first order nonlinear hyperbolic balance laws. Many difficulties arise in this research field as a result of uncertainties in data, measurements, and the use of complex data. It is very important to develop highly accurate, easy to implement, and cost-efficient methods with high resolution to study fluvial environments numerically. Although the method we employ can be quite general, for the sake of simplicity and better illustrative purposes, we focus on a prototype example of hyperbolic balance laws in this work, namely the nonlinear shallow water equations (SWEs), also referred to as the Saint-Vernant equations. This system models fluvial environments such as flood plane dynamics, coastal and tidal flows, and flow and sediment transport. It has wide applications in ocean, environmental, hydraulic engineering and atmospheric modeling. The model parameter we aim to recover in this context is the riverbed topography, however other terms such as friction may also be of interest to reconstruct.

Ground surveys of riverbeds or direct topographical data collection are not always effective for determining underwater topography because the operations are costly and time consuming. 
Numerical techniques of constructing riverbed topography can offer many benefits over aerial and ground techniques in terms of cost reduction, efficiency, as well as flexibility. Heining and Aksel \cite{Heining2009BottomStability} used a direct approach to reconstruct the bottom topography of steady-state thin-film flow. Castaings et al. \cite{Castaings2006AutomaticModeling} presented an automatic differentiation technique and free surface information to reconstruct river bed topology. Honnorat et al. \cite{Honnorat2007DassflowFlows} derived a method for recovering channel topology from a steady-state solution of the forward problem using an optimization technique called variational data assimilation. A direct approach from the one-dimensional SWEs was used by Gessese et al. \cite{Gessese2011ReconstructionFlow} to reconstruct the river bed from free surface data. Lastly, a stable finite volume scheme in the presence of wetting-drying fronts and inverse computational algorithms (based on variational approach) were presented in \cite{Monnier2016InverseDynamics}. All methods developed in these studies dealt with time-independent bottom topography functions. 

The main objective of this paper is to develop a robust algorithm that requires less data to reconstruct a dynamic bottom function, thereby allowing the construction of a more accurate and inexpensive model. We assume the measurements, possibly with noise, are taken only on two boundaries of the spatial domain in a given time period. Traditionally, the bottom topography in the SWEs is a function of space alone within the framework of the inverse problem construction. Here, we consider the time-dependent bottom topography function, which
allows for the recovery of the bottom topography with less data considering the fact that movement in the forward problem solution coming from the change of the bottom topography allows more information propagation to the boundary measurements.
Usually inverse problems are more difficult when we need to recover both the temporal and spatial profile. In this work, as a first step, we assume a special form of the topographical profile which represents two known spatial profiles and an unknown temporal interaction. This time-dependent bottom topography can practically describe a physical phenomenon when two platonic plates with known topography are moving against one another, e.g. an earthquake, underwater volcanic buildup, or a moving sand bottom. We aim to recover only the temporal profile representing, e.g. the pulse of the earthquake. After constructing the adjoint formulation, we present the cost function with two regularization terms added to suppress noise and to handle the ill-posedness of this problem. An iterative update scheme based on a three-operator splitting scheme is employed to update the targeting function. This splitting scheme requires each operator to be computed only once per iteration and is straightforward to implement.

This paper is organized as follows. In Section \ref{sec_forward} we introduce the primal equations used throughout the paper. The discontinuous Galerkin (DG) numerical scheme is also presented as the method used for solving the forward problem. Section \ref{sec_inverse} includes a discussion on the formulation of the inverse problem. The adjoint equations are derived from a linearized system and used to derive the gradient formulation of the cost function. This section also contains the description of the iterative updating procedure for determining the desired source term, as well as a discussion on the choice of regularization terms. Numerical examples are presented in Section \ref{sec_examples_swe}, and conclusions are discussed in Section \ref{sec_conclusion}. 

\section{Forward Problem} \label{sec_forward}

In this section, we present the hyperbolic PDE system used to define the forward (or primal) problem. The forward system will be used in Section \ref{sec_inverse} to determine the adjoint formulation, which in turn is employed to derive the gradient of a cost function and iteratively update the time-dependent bottom topography function. We will discuss the PDE system as well as the numerical scheme used to solve the forward problem.

\subsection{Forward Problem Formulation}

Hyperbolic balance laws are conservation laws with a source term. A few examples of such systems include the nonlinear SWEs, the arterial blood flow model, the Euler equations under gravity, and the telegrapher's equations. A source term might arise as a result of many factors, such as a friction term or a topographical term. 

In this work, we only consider the one-dimensional systems of $m$ hyperbolic balance laws which take the form
\begin{equation} \label{swe_balancelaw}
    \begin{cases}
    \partial_t \U + \partial_x \FU = \hat{\bS}(\U,\B), & (x,t) \in (x_0,x_L) \times (0,T], \\
    \U(x,0) = \U_0(x), & x \in (x_0,x_L), \\
    \end{cases}
\end{equation}
where $\U$ are the state variables, $\FU$ are the fluxes, and $\U_0(x)$ are the initial conditions. The vector $\B$ represents the model parameters 
we aim to recover in the inverse problem and is only present in the source term, $\hat{\bS}(\U,\B)$. 
The source term can be rewritten in the form of $\hat{\bS}(\U,\B) = \bS(\U,\B)\U$, in which $\bS(\U,\B) \in \mathbb{R}^m \times \mathbb{R}^m$ is a matrix. 

Alternatively, the system can be written in quasi-linear form using the Jacobian matrix $\A(\U)= \frac{\partial \F}{\partial \U}$,
\begin{equation}\label{swe_jacobiansys}
    \begin{cases}
    \partial_t \U + \A(\U)\partial_x \U = \bS(\U,\B)\U,& (x,t) \in (x_0,x_L) \times (0,T], \\
    \U(x,0) = \U_0(x), & x \in (x_0,x_L). \\
    \end{cases}
\end{equation}
To determine the solution of a forward problem, one seeks to determine the state variables $\U$ with the model parameters $\B$ given. In this work, we will only consider the case when we have a single topographical function, denoted by $B$. While the information in $B$ is traditionally a function of space alone, here we consider $B$ as a function depending on both space and time as follows:
\begin{equation}\label{bottom}
    B(x,t) = B_0(x) + p(t)B_1(x),
\end{equation}
where we assume $B_0(x)$ and $B_1(x)$ to be known and $p(t)$ to be the component we wish to recover.

In this paper we will focus on the SWEs with a non-flat bottom topography, one of the most well known systems of hyperbolic balance laws. In particular, we aim to recover the riverbed topography function, denoted by $b$. This term occurs only in the source term of the momentum equation in the form of its derivative, $\partial_x b$, so we define
\begin{equation}\label{B_SWE}
    B=\partial_x b(x,t) = \partial_x b_0(x) + p(t)\partial_x b_1(x).
\end{equation} 
The state variables, flux terms, and source term for the SWEs are given by
\begin{equation}\label{eq:SWE}
    \U=\begin{bmatrix}h \\ hu\end{bmatrix},
    \quad
    \FU = \begin{bmatrix}hu \\ hu^2+\frac{1}{2}gh^2\end{bmatrix},
    \quad
    \hat{\bS}(\U,B) = \begin{bmatrix}0 \\ -ghB\end{bmatrix}
    =\begin{bmatrix}0 \\ -gh\partial_x b\end{bmatrix},
\end{equation}
by following the form \eqref{swe_balancelaw}. Here, $h(x,t) \geq 0$ is the water height, $(hu)(x,t)$ is the water discharge with $u(x,t)$ being the depth averaged velocity, and $g = 9.812$ is the gravitational constant. On the other hand, under the quasi-linear formulation, we write
\begin{equation}
    \A(\U) = \begin{bmatrix}0 & 1 \\ gh-u^2 & 2u\end{bmatrix},
    \quad\quad\quad
    \bS(\U,B) = \begin{bmatrix}0 & 0\\ -g\partial_x b & 0\end{bmatrix}.
    \label{forward_prob}
\end{equation}

Time-dependent bottom topography functions have been considered in the literature. In a more complicated model, e.g. \cite{Sleath1995SedimentCurrents, Heinrich1992NonlinearLandslides, Hu2012NumericalChina, Li2011FullyRivers, Kurganov2018Finite-volumeEquations}, the change of the bottom function in the SWEs may depend on other state variables. For instance, the bottom function may be determined by erosion, sediment transport, dam breaks, or landslides due to floods. In such cases, additional equations to model the evolution of bottom topography may arise in the system in order to better describe these dependence. However, these models are different from the one we consider in this work.



\subsection{Discontinuous Galerkin Method for the Forward Problem}\label{sec_forward_numerical_dg}

The DG method will be used to solve the forward problem (\ref{forward_prob}). It is a high order accurate scheme that has gained significant attention in the last decade. The method is advantageous for hyperbolic conservation laws because it is both stable, similar to the finite volume method, and flexible, like the finite element method. The arbitrary order feature of the DG method can provide accurate results on a coarse mesh.   
In particular, in an inverse problem algorithm, a forward solver is usually employed during each iteration, hence the use of a coarser mesh is ideal in an effort to reduce computation cost in the iterative process. Lastly, the DG scheme is able to capture the discontinuous solutions well and help us locate the interfaces accurately.

The computational domain, $I=[x_0,x_L]$, is first discretized into $N$ cells. The $j^{\text{th}}$ cell is denoted by $I_j = [x_{j-\frac{1}{2}}, x_{j+\frac{1}{2}}]$ with $x_j$ being the center of the cell. The size of the $j^{\text{th}}$ cell is denoted by $\Delta x_j = x_{j+\frac{1}{2}}-x_{j-\frac{1}{2}}$ and we let $h = \max_j {\Delta x_j}$. We seek an approximation $\U_{h}$ of the solution $\U$, in which $\U_{h}^{(i)}$ for $i=1,...,m$ belongs to the finite dimensional piecewise polynomial space
\begin{equation} \label{fespace}
\mathbb{V}^k_{h} = \{v:v|_{I_j} \in P^k(I_j), j=1,...,J\},
\end{equation}
where $P^k(I_j)$ denotes the space of polynomials of degree up to $k$ on $I_j$. The topographical source term variable, $B$, is also projected into $\mathbb{V}_{h}^k$, and is written as $B_h$. The scheme does not require $\U_{h}$ to be continuous at the cell interface $x_{j+\frac{1}{2}}$, so we introduce the notation $\U^+_{h,j+\frac{1}{2}}$ as the limit of the solution $\U_h$ at $x_{j+\frac12}$ from the right cell $I_{j+1}$, and $\U^-_{h,j+\frac{1}{2}}$ as the limit from the left cell $I_j$. 

The DG scheme in each cell $I_j$ is based on a modification of the weak formulation of the PDE,
\begin{equation} \label{eq:dgform}
\int_{I_j}\partial_t \U_h \bv \ dx 
- \int_{I_j} \F(\U_h) \partial_x \hmy{\bv} \ dx 
+ \hat{\F}_{j+\frac{1}{2}}\bv^-_{j+\frac{1}{2}}
-\hat{\F}_{j-\frac{1}{2}}\bv^+_{j-\frac{1}{2}}
=\int_{I_j} \hat{\bS}(\U_h,B_h) \bv \ dx, 
\end{equation}
where $\bv(x)$ is a vector of test functions from the test space $\mathbb{V}_{h}^k$, and the numerical flux 
$\hat{\F}_{j+\frac{1}{2}} = f\left(\U^-_{h,j+\frac{1}{2}}, \U^+_{h,j+\frac{1}{2}}\right)$ is a function that takes information from both the left and right side of the cell interface. We implement the simple Lax-Friedrichs flux
\begin{equation} \label{lfflux}
f(a,b)=\frac{1}{2}\Big(\F(a)+\F(b)-\alpha(b-a)\Big),
 \end{equation}
where $\alpha$ is chosen to be the maximum of the largest eigenvalue of $\A(\U)$ over the entire computational domain or locally in each cell. In the case of the SWEs, let $\alpha = \max{(|u| + \sqrt{gh})}$.

The semi-discrete method \eqref{eq:dgform} can be rewritten in the ODE form as 
\begin{equation*}
    \partial_t \U_h = \mathcal{F}(\U_h),
\end{equation*}
after representing $\U_h$ as a linear combination of the basis functions of $\mathbb{V}^k_{h}$.
In order to advance the scheme in time, we use the high order Strong-Stability Preserving Ruge-Kutta (SSP-RK) temporal discretization \cite{Gottlieb2001StrongMethods}. Throughout this paper, the third order SSP-RK method of the form
\begin{equation} \label{rkscheme}
\begin{split}
\U_h^{(1)} & = \U_h^n + \Delta t \mathcal{F}(\U_h^n), \\
\U_h^{(2)} & = \frac{3}{4}\U_h^n + \frac{1}{4}\Big(\U_h^{(1)} + \Delta t \mathcal{F}\big(\U_h^{(1)}\big)\Big), \\
\U_h^{n+1} & = \frac{1}{3}\U_h^n + \frac{2}{3}\Big(\U_h^{(2)} + \Delta t \mathcal{F}\big(\U_h^{(2)}\big)\Big).
\end{split}
\end{equation}
is used. \hlc{Furthermore, the equation for the space and time dependent bottom topography function} \eqref{bottom} \hlc{must also involve a temporal discretization. As $\mathcal{F}(\U_h^n)$ is evaluated at time $t^n$, $\mathcal{F}\big(\U_h^{(1)}\big)$ is evaluated at time $t^n + \Delta t$, and $\mathcal{F}\big(\U_h^{(2)}\big)$ is evaluated at time $t^n + \frac{1}{2}\Delta t$, we evaluate $B(x,t)$ at the same time values for each Runge-Kutta step. When the function values at $t^n+\frac{1}{2}\Delta t$ are unavailable, they are determined via quadratic interpolation.}

\section{Inverse Problem} \label{sec_inverse}

In this section, we describe the inverse problem of our focus, which is the reconstruction of the topographical source function $B$ from boundary data of the hyperbolic conservation law \eqref{swe_balancelaw} from a single measurement event.
We reduce the inverse problem to an optimization problem of a residual functional coming from boundary measurements, with an addition of two regularization terms, which will be described more concretely later in the section. 

In our work, we adopt the adjoint method to numerically obtain a gradient of our functional. We will describe the cost function we wish to minimize, the derivation of the adjoint formulation for the gradient calculation, and the iterative update scheme for determining the time component, $p(t)$, of the source function $B$. The numerical solution to the adjoint problem will be calculated using the DG method.

We assume that noisy measurements of $\U$ are taken only on both boundaries of the spatial domain, given a period of time $[0,T]$ in one single measurement event.
During the numerical reconstruction process, we assume that only these noisy solutions at the boundary of the spatial domain is known to us. For notational sake, we denote these noisy measurements as $\hat{\Lambda}^{\text{noisy}} = \hat{\Lambda}\mu$ where the multiplicative noise is uniformly distributed, $\mu \sim \mathcal{U}\left[1-\frac{1}{2}\eta_{meas},1+\frac{1}{2}\eta_{meas}\right]$, with a given noise level $\eta_{meas}$. The goal of the inverse scheme is to find the function, $p(t)$, that provides the best approximation $\Lambda(B(p)) \approx \hat{\Lambda}^{\text{noisy}}$. The map $\Lambda(B(p))= \U|_{\{x_0,x_L\} \times [0,T]} $ represents the forward map with the input, $B(p)$, as the topographical function and the output as the solutions, $\U$, at the boundary points, $\{x_0,x_L\}$, over the time interval $[0,T]$. \hlc{We thus declare the control to be $p(t)$ and the number of control variables to be equivalent to the number of time steps in the numerical scheme. On the other hand, the observation values are $\Lambda(B(p))= \U|_{\{x_0,x_L\} \times [0,T]} $. The number of observation values is determined by taking the product of the number of state variables ($m$), the number of boundary points (two in 1D), and the number of time steps in the numerical scheme.}  

Finding the best approximation reduces to minimizing the error or residue of the predicted and measured data for all time at the boundary of the computational domain. This corresponds to minimizing the functional 
\begin{equation} \label{cost_funct}
    \J(p) := \J_0(p) + \R(p) := \int_0^T \frac{1}{2}\left|\left[\mathcal{E}(p)\right](x_0,t)\right|^2 + \frac{1}{2}\left|\left[\mathcal{E}(p)\right](x_L,t)\right|^2 dt + \R(p),
\end{equation}
where the error function, $\mathcal{E}$, for a given $p$ is defined as
\begin{equation} \label{error_funct}
    \left[\mathcal{E}(p)\right] (x,t) = \Lambda(B(p)) (x,t)-\hat{\Lambda}^{\text{noisy}}(x,t).
\end{equation}
The term $\R(p)$ is a regularization term that will be discussed in Section \ref{sec_regularization}.

The optimization problem becomes
\begin{equation} \label{min_prob}
    \begin{split}
        & \text{minimize } \J(p) := \J_0(p) + \R(p) \quad \text{subject to } \eqref{swe_jacobiansys}.
    \end{split}
\end{equation}
We employ a descent method to minimize the above functional, in which the (formal) gradient  $\nabla \J$ will be obtained via the adjoint method following a linearization process of the equation \eqref{swe_jacobiansys}. 

\subsection{Gradient Derivation}
In an effort to determine $\nabla \J_0$, we begin by calculating the variational derivative (in the sense of the Gateaux differential) of $\J_0$ and dualize it using $L^2$-pivoting.
In what follows, we would like to denote, for a functional $\mathcal{F}$, the variational derivative of $\mathcal{F}$ at $p$ along $\tp$ as 
\begin{equation}\label{definition}
    \delta \mathcal{F}(p;\tp) := \lim_{\epsilon \rightarrow 0} \frac{\mathcal{F}(p+\epsilon \tp) - \mathcal{F}(p)}{\epsilon} 
\end{equation}
whenever it exists.  Furthermore, whenever $\delta \mathcal{F}(p;\tp)$ is linear with respect to $\tp$, we (formally) dualize the variational derivative $\delta \mathcal{F}(p;\tp)$ using $L^2$-pivoting and define the gradient, $\nabla \mathcal{F}(p)$, such that it satisfies the relation
\begin{equation}\label{grad_def}
    \delta \mathcal{F} (p;\tp) := \int_{0}^{T} \left[\nabla \mathcal{F}(p)\right](t) \, \tp (t)\ dt.
\end{equation}

With these notions at hand, we readily compute that
\begin{equation}\label{var_derv}
\begin{split}
    \dJ_0(p;\tp) 
    & = \lim_{\epsilon \rightarrow 0} \frac{\J_0(p+\epsilon \tp) - \J_0(p)}{\epsilon} \\
   & = \int_0^T  \bigg(\left[\delta\mathcal{E}^T (p;\tp) \, \mathcal{E} (p)   \right](x_0,t)
    + \left[\delta\mathcal{E}^T (p;\tp) \,  \mathcal{E} (p) \,  \right](x_L,t)\bigg) \ dt \,,
\end{split}
\end{equation}
where the superscript $^T$ now represents the transpose of a matrix (and not the adjoint operator).
From the definition of $\mathcal{E}$, we quickly realize that $\delta\mathcal{E}(p;\tp)(x,t) = \delta\Lambda (B(p);B(\tp))$. Hence, \eqref{var_derv} reduces to
\begin{equation}\label{var_derv3}
\begin{split}
      \dJ_0(p;\tp) & = \int_0^T   \left(\left[\delta\Lambda^T (B(p);B(\tp)) \, \mathcal{E} (p) \, \right](x_0,t)
    + \left[\delta\Lambda^T (B(p);B(\tp)) \, \mathcal{E} (p) \right](x_L,t)\right) \ dt \,.
\end{split}
\end{equation}
We now see the necessity of evaluating the term $\delta\Lambda (B(p);B(\tp)) $ explicitly. Albeit seemingly complicated, the difficulty of the evaluation will be mitigated via solving a related adjoint equation, which will be described in the next subsection.

\subsection{Linearization \& Adjoint Formulation} \label{subsec_adjoint}
The adjoint formulation can be understood from multiple perspectives. One way is the Lagrange framework in which the adjoint variables are Lagrange multipliers. This method is commonly used in the aeronautical community, popularized by Jameson \cite{Jameson1988AerodynamicTheory}, because it provided a solid connection to theories of constrained optimal control and optimization. Another type of approach, the duality framework, requires one to linearize the system in order to derive the adjoint equations. We will use the duality framework in this paper, however the Lagrange framework provides the exact same adjoint formulation. 

\subsubsection{Linearization of the forward system}
In this subsection, we aim to linearize the forward system \eqref{swe_jacobiansys} as follows.
We consider an $\epsilon$-perturbation of $B$, $B^{\varepsilon}:= B + \epsilon\tB$ along the direction $\tB$, and see how the resulting
$\U$ that satisfies \eqref{swe_jacobiansys} is perturbed. We denote $\U^{\epsilon}$ as the solution to \eqref{swe_jacobiansys} given $B^{\varepsilon}$ and define
\begin{eqnarray}\label{u_tilde}
 \tU := \lim_{\epsilon \rightarrow 0} \frac{\U^{\epsilon} - \U}{\epsilon} 
 \end{eqnarray}
whenever it exists.
Now we quickly realize that
\begin{equation} \label{pert_ics}
\tU(x,0) = \lim_{\epsilon \rightarrow 0} \frac{\U^{\epsilon}(x,0) - \U(x,0)}{\epsilon}  = 0
\end{equation}
as the initial conditions of $\U^{\epsilon}$ and $ \U$ shall coincide.
Moreover, taking the differences of the respective equations coming from \eqref{swe_jacobiansys} for $\U^{\epsilon}$ and $ \U$ directly gives
\begin{equation} \label{linear_eq1}
    \begin{split}
        0 & = \lim_{\epsilon \rightarrow 0} \frac{1}{\epsilon}\left(\partial_t [ \U^{\epsilon} - \U ] + \left[ \A(\U^{\epsilon})\partial_x\U^{\epsilon} - \A(\U)\partial_x\U\right] - \left[\bS(\U^{\epsilon},B^\epsilon) \U^{\epsilon}- \bS(\U, B)\U\right]\right), 
    \end{split}
\end{equation}
and each term in the bracket can be simplified whenever they exist.  For instance, we directly have
\begin{eqnarray} 
     \lim_{\epsilon \rightarrow 0} \frac{1}{\epsilon}\partial_t \left( \U^{\epsilon} - \U \right) & = & \partial_t \tU. \label{linear_eq11}
\end{eqnarray}
Meanwhile, we may simplify the flux term as
\begin{align}
&\lim_{\epsilon \rightarrow 0} \frac{1}{\epsilon} \left[ \A(\U^\epsilon)\partial_x\U^\epsilon - \A(\U)\partial_x\U\right]
\label{linear_flux} \\
&\quad = \lim_{\epsilon \rightarrow 0} \frac{1}{\epsilon} \left[\A(\U)\left(\partial_x\left(\U + \epsilon \tU\right)-\partial_x(\U)\right) + \left(\A(\U + \epsilon \tU)-\A(\U)\right)\partial_x(\U + \epsilon \tU)\right] \notag\\
&\quad =  \A(\U)\partial_x \tU + \left(\sum_{k =1}^m \partial_{\U_k}A(\U)\,\tU_k \right)\partial_{x} \U \notag \\
       &\quad = \A(\U)\partial_x \tU + \left(\sum_{j=1}^m \partial_{\U_j}A(\U)\partial_{x} \U \right)\tU \notag \\
        &\quad = \A(\U)\left[\partial_x \tU\right] + \left[\partial_x \A(\U)\right] \tU \notag \\
        &\quad =  \partial_x \left[\A(\U) \tU\right]. \notag
\end{align}
\hly{Here, the first equality is derived from rewriting the flux term after adding and subtracting $\A(\U)\partial_x\U^\epsilon$ and applying the definition} \eqref{u_tilde}. \hly{The second equality follows from applying the limit in which the second term is the result of partial derivatives involving the chain rule. The third equality is determined from the symmetry relation $\partial_{U_k} \A_{ij}(U)  = \partial_{U_k} \partial_{U_j} \F_i (U) = \partial_{U_j} \A_{ik}(U)$. The last equality follows from the product rule.}

Likewise, we can simplify the source term and obtain
\begin{align} 
\label{linear_source2}
& \lim_{\epsilon \rightarrow 0} \frac{1}{\epsilon}\left(\bS(\U^\epsilon, B^\epsilon)\U^{\epsilon} - \bS(\U,B)\U\right)  \\
& \qquad\qquad = \lim_{\epsilon \rightarrow 0} \frac{1}{\epsilon} \left(\bS(\U,B)\left(\U^{\epsilon} - \U\right) + \left(\bS(\U^\epsilon,B^\epsilon)-\bS(\U,B)\right)\U^\epsilon\right) \notag\\
&\qquad\qquad =  \bS(\U,B)\tU 
        + \lim_{\epsilon \rightarrow 0} \frac{1}{\epsilon}\left(\bS(\U+ \epsilon\tU,B+ \epsilon\tB)
       -\bS(\U,B)\right)\left(\U+\epsilon\tU\right) \notag \\
        &\qquad\qquad= \bS(\U,B)\tU 
        + \left(\sum_{i=1}^m \partial_{\U_i}\bS(\U,B)\tU_i 
        + \partial_{B}\bS(\U,B)\tB \right) \U \notag \\
        &\qquad\qquad :=\left(\bS(\U,B) 
        + \bm{C}(\U,B)\right)\tU 
        + \partial_{B}\bS(\U,B)\U\tB, \notag
\end{align} 
where $\bm{C}$ denotes the matrix $\bm{C}_{ij} = \sum_{k=1}^m \frac{\partial \bS_{ik}}{\partial \U_j}\U_k$. \hly{The first equality is derived from rewriting the original equality after adding and subtracting the term $\bS(\U,B)\U^\epsilon$. The second equality follows from applying the limit to the first term and from applying the definitions in} \eqref{u_tilde} \hly{and of $B^\epsilon$ to the second term. The third equality follows from taking partial derivatives which involves the chain rule and applying the limit. Lastly, the fourth equality is the result of a symmetry relation.}

Substituting \eqref{linear_eq11}, \eqref{linear_flux} and \eqref{linear_source2} into
\eqref{linear_eq1}, and combining that with the initial condition \eqref{pert_ics}, we therefore obtain the following linear system for $\tU$
\begin{equation} \label{linear_eq2}
    \begin{cases}
    \left(\partial_t - \bS -\bm{C} \right) \tU +\partial_x\left(\A\tU\right) = \left(\partial_{B}\bS(\U,B)\U\right)\tB,  & (x,t) \in (x_0,x_L) \times (0,T], \\
    \tU(x,0) = 0, & x \in [x_0,x_L] \,, \\
    \end{cases}
\end{equation}
which serves as the linearization of the forward system \eqref{swe_jacobiansys}.

\subsubsection{The adjoint system}

With the linearization process given in the previous subsection, we may proceed to obtain $\delta\Lambda (B(p);B(\tp))$ at the boundary points, and thereby evaluate $\dJ_0 (p;\tp)$ appropriately.

We start by considering $\sig$ which satisfies the following adjoint system with final time condition and boundary conditions
\begin{equation} \label{adjoint_eq}
    \begin{cases}
    (\partial_t+ \A^T\partial_x + \bS^T +\bm{C}^T)\sig = 0, & x\in (x_0,x_L) \times (0,T], \\
    \sig(x,T) = 0, & x\in (x_0,x_L), \\
    \sig(x_0,t) = -(\A^T)^{-1}(x_0,t)\left[\mathcal{E}(p)\right](x_0,t), & t \in (0,T], \\
    \sig(x_L,t) = (\A^T)^{-1}(x_L,t)\left[\mathcal{E}(p)\right](x_L,t) & t \in (0,T]\,. \\ 
    \end{cases}
\end{equation}
In the particular case of the SWEs, the matrices appearing in \eqref{adjoint_eq} are given by
\begin{equation}
    \A^T = \begin{bmatrix}0 & gh-u^2 \\ 1 & 2u\end{bmatrix},
    \quad\quad\quad 
    \bS^T = \begin{bmatrix}0 & -g\partial_x b \\ 0 & 0\end{bmatrix},
    \quad\quad\quad 
    \C^T = 0.
\end{equation}

Taking the inner product of the solution $\sig$ of \eqref{adjoint_eq}  and the weak formulation of the linearized system in \eqref{linear_eq2}, we get
\begin{equation} 
\begin{split}
\label{ibp_eq}
 &\int_0^T \int_{x_0}^{x_L} \sig^T \left(\partial_{B}\bS(\U,B)\U\right)\tB \ dxdt \\
    &\qquad =\int_0^T \int_{x_0}^{x_L} \sig^T \left(\partial_t+\partial_x \A - \bS -\bm{C} \right)\tU \ dxdt  \\
    &\qquad =
    -\int_0^T \int_{x_0}^{x_L} \tU^T \left(\partial_t+ \A^T\partial_x +\bS^T + \bm{C}^T\right)\sig \ dxdt  \\
    &\qquad\qquad + \int_{x_0}^{x_L} \tU^T \sig\bigg|_{t=0}^{t=T} \ dx
    + \int_0^T \tU^T \A^T\sig\bigg|_{x=x_0}^{x=x_L} \ dt ,
\end{split}
\end{equation}
where we simplify further, with \eqref{adjoint_eq}, to obtain
\begin{equation} \label{new_eq}
\begin{split}
    & \int_0^T \int_{x_0}^{x_L} \sig^T(x,t) \left(\partial_{B}\bS(\U,B)\U\right) (x,t) \tB(x,t) \ dxdt \\
    & = \int_0^T \left[\delta{\Lambda}^T(B(p);\tB(p)) \, \mathcal{E}(p)\right](x_L,t) \ dt \\
    & \qquad + \int_0^T\left[\delta{\Lambda}^T(B(p);\tB(p)) \, \mathcal{E}(p)\right](x_0,t)\ dt. \\
\end{split}
\end{equation}
Here the last equality follows from the choice of boundary conditions described in \eqref{adjoint_eq} and the fact that $\U = \Lambda(B(p)) $ implies $\tU = \delta{\Lambda}( B(p) ;\tB(p) )$. 
We may now readily substitute \eqref{new_eq} into the expression \eqref{var_derv3} to obtain
\begin{equation} \label{var_derv4}
\begin{split}
      \dJ_0(p;\tp) & = \int_0^T \int_{x_0}^{x_L} \sig^T(x,t) \left(\partial_{B}\bS(\U,B)\U\right) (x,t) \tB(x,t) \ dxdt. 
\end{split}
\end{equation}
By utilizing the fact that
\begin{equation} \label{tB}
    \left[\tB(p)\right] (x,t)= \delta B ( p ; \tilde p ) = [ \delta ( B_0 + p B_1 )  ] ( p ; \tilde p )  =  B_1(x) \, \tp(t) ,
\end{equation}
we further simplify \eqref{var_derv4} to
\begin{equation}
\begin{split}
\label{new_eq1}
      \dJ_0(p;\tp) & =  \int_0^T \left(\int_{x_0}^{x_L}  \sig^T(x,t) \left(\partial_{B}\bS(\U,B)\U\right) (x,t) B_1(x) \ dx\right) \tp(t) \ dt.
\end{split}
\end{equation}
Therefore from definition \eqref{grad_def}, we obtain the following (formal) gradient from \eqref{new_eq1}
\begin{equation} \label{gradient2}
    \nabla \J_0(p)(t) 
    = \int_{x_0}^{x_L} \sig^T(x,t) \left(\partial_{B}\bS(\U,B)\U\right) (x,t) B_1(x) \ dx \,.
\end{equation}

We again remark that, in the case of the SWEs, we have
\begin{equation}
    \partial_{B}\bS(\U,B)\U = \begin{bmatrix}0 & 0 \\ -g & 0\end{bmatrix}\begin{bmatrix}h \\ hu\end{bmatrix} = \begin{bmatrix}0 \\-gh \end{bmatrix}, 
    \qquad B_1 = \partial_x b_1,
\end{equation}
and therefore the gradient is simplified to the form
\begin{equation}
\begin{split}
    \nabla \J_0 (p) 
    & = \int_{x_0}^{x_L} \left(\begin{bmatrix} \sigma_1 \\ \sigma_2 \end{bmatrix}^T \begin{bmatrix}0 \\ -gh\end{bmatrix}\right)(x,t)\partial_x b_1(x)\ dx \\
    & = \int_{x_0}^{x_L} -g \sigma_2(x,t) h(x,t) \partial_x b_1(x) \ dx,
\end{split}
\end{equation}
where $\sig = \begin{bmatrix} \sigma_1 \\ \sigma_2 \end{bmatrix}$ is the solution of the adjoint equation \eqref{adjoint_eq}.

%

\subsection{Numerical Scheme for the Inverse Problem}
In this subsection we will discuss the numerical algorithms for the inverse problem. The DG scheme will be employed to solve the adjoint problem \eqref{adjoint_eq} and an iterative method will be presented to update the function $p$ with the suitably chosen regularization terms.

\subsubsection{Discontinuous Galerkin Method for the Adjoint Problem}\label{sec_adjoint_numerical_dg}
Noting that the spatial derivative in the adjoint problem \eqref{adjoint_eq} is not in the conservative form, 
we start by reformulating the adjoint problem as a balance law of the form 
\begin{equation}
    \partial_t \sig + \partial_x\left(\A^T\sig\right)  = (\partial_x\A^T-\bS^T - \C^T)\sig,
\end{equation}
where $\A=\A(\U)$ does not depend on the unknown $\sig$.
Following the same discretization strategy as presented in Section \ref{sec_forward_numerical_dg}, we seek an approximate solution $\sig_h$ in which $\sig_h^{(i)}$ for $i=1,..., m$ belong to $\mathbb{V}_h^k$. The DG method in cell $I_j$ becomes
\begin{equation} \label{eq:dgform_adjoint}
\begin{split}
\int_{I_j}\partial_t \sig_h \bv \ dx 
& - \int_{I_j} \A(\U_h)^T\sig_h \partial_x \hmy{\bv} \ dx 
+ \hat{\G}_{j+\frac{1}{2}}\bv^-_{j+\frac{1}{2}}
-\hat{\G}_{j-\frac{1}{2}}\bv^+_{j-\frac{1}{2}} \\
& =\int_{I_j} \left(\partial_x \A^T(\U_h) - \bS^T(B_h) - \C^T\right)\sig_h \bv \ dx, 
\end{split}
\end{equation}
where $\bv \in \mathbb{V}^k_h$ is a vector of test functions and the Lax-Friedrichs numerical flux takes the form
\begin{equation}
    \hat{\G}_{j+\frac{1}{2}} 
    = \frac{1}{2}\left(
    \A(\U_{h,j+\frac{1}{2}}^-)^T\sig_{h,j+\frac{1}{2}}^- 
    + \A(\U_{h,j+\frac{1}{2}}^+)^T\sig_{h,j+\frac{1}{2}}^+ 
    - \alpha\left(\sig_{h,j+\frac{1}{2}}^+-\sig_{h,j+\frac{1}{2}}^-\right)\right),
\end{equation}
with the value of $\alpha$ being the same as in the forward DG scheme, described in Section \ref{sec_forward_numerical_dg}.

\subsubsection{Regularization \& Update Scheme} \label{sec_regularization}

In this subsection, we describe the numerical method designed to recover the function $p(t)$ via an iterative scheme.   Usually, either a descent type \cite{CAUCHY1847MethodeSimultanees}, Newton type \cite{Hintermuller2010SemismoothApplications}, or a trust region algorithm \cite{Sun2006OptimizationProgramming} is employed.   A Newton type algorithm usually provides a certain acceleration to the convergence, but as a trade off, it is usually more computationally expensive.  In this work, we employ a descent type algorithm, in light of the fact that our functional is highly non-linear and highly non-convex, applying a higher order method may result in getting stuck at a local optimum even more easily. 

We employ an operator splitting algorithm to update the function $p(t)$.   The scheme is initialized with a random initial guess for $p$, denoted by $p^{0,\text{noisy}}$. We use multiplicative noise following a uniform distribution, i.e., $p^{0,\text{noisy}} = p^{0}\nu$, where $\nu \sim \mathcal{U}\left[1-\frac{1}{2}\eta_{p},1+\frac{1}{2}\eta_{p}\right]$, to define the random initial guess. Not only is the update of $p$ dependent on $\nabla \J_0$, but it also relies on a regularization term. The regularization term ensures the optimization problem is locally convex and makes it possible to solve an ill-posed problem efficiently by incorporating a-priori knowledge of the profile to be reconstructed. Various regularization terms have been constructed for different purposes. For instance, $L^1$ regularization \cite{Yin2008BregmanSensing,Goodfellow2016DeepLearning} results in a simpler sparse solution. On the other hand, TV regularization \cite{Rudin1992NonlinearAlgorithms} favors piecewise constant functions of the coefficients to be recovered, whereas Sobolev regularization \cite{Fischer2017SobolevAlgorithm} favors smoothness of the coefficients to be reconstructed.

Before we focus on our choice of regularization, we first discuss the update algorithm.
To better motivate our choice of algorithm, we start by simplifying our discussion and considering the situation when there is only one regularization term.  In this case, the proximal gradient descent method (or the forward-backward splitting) \cite{Beck2017First-OrderOptimization,Parikh2014ProximalAlgorithms} is a common choice. The explicit term is usually assigned as the term coming from the gradient of a more complicated functional. The implicit term is typically chosen so that the proximal map is easy to evaluate and the stability of the algorithm is increased. Consequently, the $k^{th}$ iteration is given by
\begin{equation} \label{time_update}
    \begin{split}
        & p^{k+1} 
        = p^k - \ell_k \nabla \J_0\left(p^k\right) - \ell_k \partial \R\left(p^{k+1}\right), \\
    \end{split}
\end{equation}
where $\ell_k$ is the step size or learning rate and $\partial$ represents the subgradient when the proximal map of $\R$ can be computed. In our work, we will choose a constant step size, i.e., $\ell_k = \ell$ for all $k$.  The scheme can be rewritten so that the update for iteration $k+1$ only depends on the information from iteration $k$,
\begin{equation} \label{time_update1}
        p^{k+1} = \left(\mathcal{I} + \ell \partial \R\right)^{-1}\left(p^{k} - \ell \nabla \J_0\left(p^{k}\right)\right), \\
\end{equation}
where $\mathcal{I}$ is the identity matrix and
\begin{equation} \label{time_update_2}
    \left(\mathcal{I} + \ell \partial \R\right)^{-1}(w) 
    = \text{argmin}_y \left\{\R(y) + \frac{1}{2\ell} \|w-y\|_2^2\right\} 
    = \text{prox}_{R,\ell}(w).
\end{equation}
This leads to the formulation
\begin{equation}
     p^{k+1} 
     = \text{argmin}_y \left\{\R(y) + \frac{1}{2\ell} \|p^{k}-\ell \nabla \J_0(p^{k}) - y\|_2^2\right\}. 
\end{equation}
A common choice for the regularization is $L^1$ regularization, where $\R(y) = \gamma \|y - p_0 \|_1$ with $p_0$ being a chosen coefficient of homogeneous background and $\gamma$ a scalar parameter, aiming to impose sparsity of the difference between the resulting optimum and $p_0$.  
The proximal gradient method coming from this choice of regularizer is
\begin{equation}
\begin{split}
     p^{k+1} 
     & = \text{argmin}_y \left\{\gamma \|y - p_0 \|_1 + \frac{1}{2\ell} \|p^{k}-\ell \nabla \J_0\left(p^{k}\right) - y\|_2^2\right\} \\
     & = \mathcal{S}_{\gamma \ell}\left(p^{k}-\ell \nabla \J_0\left(p^{k}\right) - p_0 \right) + p_0,
\end{split}
\end{equation}
where the shrinkage operator $\mathcal{S}_{\gamma \ell}$ \cite{Darbon2016AlgorithmsElsewhere,Darbon2015OnEquations,Yin2008BregmanSensing} is given as follows
\begin{equation}\label{shrinkage}
\mathcal{S}_{\gamma \ell }\left( p \right) =  \text{sign}\left( p \right)\max \left\{\left|p\right| - \ell \gamma,0\right\}\,.
\end{equation}

After briefly describing the simple motivating example which carries only one regularization term, we now describe the combination of regularization terms that we use in our work, and how we perform the operator splitting in our algorithm.
In this paper, the regularization term is taken as a sum of two regularizers
\begin{equation}
    \R(p) = \R_{L^1}(p - p_0) + \R_{H^1}(p),
\end{equation}
where $\R_{L^1}(p) = \|p\|_1$ represents $L^1$ regularization and $\R_{H^1}(p) = \|\nabla p\|_2^2$ represents $H^1$ regularization. The $L^1$ regularization term will aid in removing the noise by sparsifying it, while the $H^1$ regularization term will be beneficial for the purpose of smoothing out the noisy data, an advantage over total variation ($\|\nabla p\|_1$) regularization. $H^1$ regularization has been shown to be good for flow control problems \cite{Collis2001NumericalEquations, Heinkenschloss1998FormulationFlow} as well as image reconstruction and deblurring \cite{Ng2000CosineReconstruction, Kindermann2005DeblurringFunctionals}. 

Now, we wish to minimize $\J_0(p) + \R_{L^1}(p) + \R_{H^1}(p)$, which reduces to finding $p$ so that
\begin{equation}
    0 \in \nabla \J_0 (p) + \partial \R_{L^1}(p) + \nabla \R_{H^1}(p).
\end{equation}
We must be careful in our approach and employ a more complicated splitting scheme than the proximal gradient descent since we now have an additional operator. In our work, we adapt the three-operator splitting algorithm \cite{Davis2017AApplications}, which we will describe in detail in Algorithm \ref{alg_oper_splitting}. For simplicity, we introduce the following notation corresponding to each regularization term
\begin{equation}
\begin{split}
    \mathcal{J}_{\ell \gamma_L \R_{L^1}}(\omega)  
    & = \left(\mathcal{I} + \ell \gamma_L \R_{L^1}\right)^{-1}(\omega) \\ 
    & = \text{sign}(\omega - p_0)\max\left\{\left|\omega - p_0\right| - \ell\gamma_L, 0\right\} + p_0,
\end{split}
\end{equation}
and
\begin{equation}
    \mathcal{J}_{ \ell \gamma_H \R_{H^1}}(\omega) 
   = (\mathcal{I} + \ell \gamma_H \R_{H^1})^{-1}(\omega) = \left(\mathcal{I} - \ell\gamma_H \Delta \right)^{-1}(\omega),
\end{equation}
where $\gamma_L$ is the $L^1$ regularization parameter and $\gamma_H$ is the $H^1$ regularization parameter. Furthermore, we would like to note that the gradient of the cost function is time-dependent, i.e. $\nabla \J_0(p^k,t)$, but we denote it as $\nabla \J_0(p^k)$ for the sake of simplifying notation. The update becomes
\begin{equation}
    p^{k+1} = \mathcal{J}_{\ell \gamma_L \R_{L^1}} \circ \left[z^k + \lambda_k\left(\mathcal{J}_{ \ell \gamma_H \R_{H^1}} \circ \left[2p^k - z^k - \ell\nabla \J_0(p^k)\right]-p^k\right)\right],
\end{equation}
where $z^0$ is originally initialized to be $p^0$ and $\lambda_k$ is the relaxation parameter which can be used to help speed up the rate of convergence of the iterative solutions.
We are now ready to introduce our algorithm.

\begin{algorithm}[H]
\SetAlgoLined
 initialize $p^0$ to be the random initial guess\; \\
 initialize $z^0 = p^0$ \; \\
 set regularization parameters $\gamma_L, \gamma_H$\; \\
 set relaxation parameter $(\lambda_k)_{k\geq 0}$\; \\
 set learning rate $\ell$\; \\
 \For{$k=0,1,...$}{
    compute $\Lambda(B(p^k))$ from $B(p^k)$ by solving the forward problem \eqref{swe_jacobiansys} \; \\
    compute $\sig$ from $(B(p^k), \Lambda(B(p^k)) )$ by solving the adjoint problem \eqref{adjoint_eq}  \; \\
    evaluate $\nabla \J_0(p^k) = \int_{x_0}^{x_L}
     \sig^T(x,t) \left(\partial_{B}\bS(\U,B)\U\right)(x,t) B_1 (x) \ dx$\; \\
    define $\kappa^k = 2p^k - z^k - \ell\nabla \J_0(p^k)$ \; \\
    evaluate $\omega^k = \mathcal{J}_{\ell \gamma_H \R_{H^1}}(\kappa^k) = (\mathcal{I}- \ell \gamma_H \Delta )^{-1} ( \kappa^k)$\; \\
    update $z^{k+1} = z^k + \lambda_k(\omega^k - p^k)$\; \\
    update $p^{k+1} = \mathcal{J}_{ \ell \gamma \R_{L^1}}(z^{k+1}) = \text{sign}(z^{k+1}-p_0)\max\{|z^{k+1}-p_0| - \ell\gamma_L,0\} +p_0$\;
    }
 \caption{Three-Operator Splitting Algorithm}
 \label{alg_oper_splitting}
\end{algorithm}

\section{Numerical Examples for the Shallow Water Equations} \label{sec_examples_swe}
In this section, we will be considering the one-dimensional nonlinear SWEs \eqref{eq:SWE}. We aim to recover the temporal component $p(t)$ in the bottom topography function $b(x,t)$, see \eqref{bottom}. 

In all the numerical tests, we use a relaxation parameter of $\lambda_k = 1$, a noise parameter for the measured data of $\eta_{meas}=0.1$ or $5\%$ noise, a noise parameter for the initial guess $p^0$ of $\eta_p=0.25$ or $12.5\%$ noise. A coefficient of homogeneous background is assumed to be known and taken as $p_0=1$. Each test is ran for 1000 iterations. The iteration with the smallest residue, $\J_0(p^k)$, is selected as the best recovered representation for the true temporal component of the bottom function $p(t)$. The measured data is computed using a high order accurate DG method with a uniform mesh of 400 cells and $P^3$ piecewise polynomials, with noise added to the DG solutions to represent noisy measurement.
A mesh of $50$ uniform cells with $P^{2}$ piecewise polynomials are used to solve the forward problem and a uniform mesh of 25 cells with $P^{1}$ piecewise polynomials are used to solve the adjoint problem, unless stated otherwise. The measured data, forward, and adjoint solvers are designed with different meshes and polynomial degree approximations in an effort to avoid committing ``inverse crime''\cite{Wirgin2004TheCrime}.  

\subsection{Tests for Recovering Different Time Profiles of \texorpdfstring{$p(t)$}{p(t)}} \label{sec_test_swe_1}
In this subsection, we will perform numerical experiments aiming to recover several unknown time profiles, $p_{true}(t)$, from noisy boundary measurements. 

We solve the forward problem \eqref{swe_balancelaw} with the DG method described in Section \ref{sec_forward_numerical_dg} where our computational domain is chosen to be $[x_0,x_L]=[0,1]$, the initial conditions are given by
\begin{equation} \label{ics_swe_smooth0}
    h(x,0) = 7 + \exp(\sin(2\pi x)), 
    \quad\quad 
    hu(x,0) = \cos(2 \pi x),
\end{equation}
and the spatial components of the bottom topography function are defined as
\begin{equation}\label{bvals_swe_smooth0}
    b_0(x) = \cos(\sin(2 \pi x), 
    \quad\quad 
    b_1(x) = \sin^2\left(\pi x\right).
\end{equation}
The final time is set as $T=0.05$ and periodic boundary conditions are used. 

We examine several choices for the true value of $p(t)$ and the corresponding initial guesses, which are outlined in Table \ref{tab:swe_p_1}. A constant learning rate of $\ell =0.6$ is used in each test. The regularization parameters are fixed with $\gamma_L = 1 \times 10^{-6}$ in all examples and $\gamma_H = 5 \times 10^{-8}$ in cases \eqref{case_swe_2a},\eqref{case_swe_2b}, \eqref{case_swe_2d}, \eqref{case_swe_2e}, $\gamma_H = 1 \times 10^{-8}$ in cases \eqref{case_swe_2c}, \eqref{case_swe_2d} and $\gamma_H = 5 \times 10^{-9}$ in case \eqref{case_swe_2f}.

\begin{table}[h]
    \centering
    \begin{tabular}{c c c}
    \toprule
        Case & $p_{true}(t)$ & $p^0(t)$ \\
        \midrule
        \stepcounter{subsecval}\refstepcounter{subsecvalnext}(\thesubsecvalnext)\label{case_swe_2a} 
        & $e^{\beta\left(t-\frac{1}{3}T\right)^2}+1$ 
        & $e^{\beta\left(t-\frac{2}{3}T\right)^2}+1$ \\
        \refstepcounter{subsecvalnext}(\thesubsecvalnext)\label{case_swe_2b}  
        & $e^{\beta\left(t-\frac{2}{3}T\right)^2}+1$ 
        & $e^{\beta\left(t-\frac{1}{3}T\right)^2}+1$\\
        \refstepcounter{subsecvalnext}(\thesubsecvalnext)\label{case_swe_2c}  
        & $e^{2\beta\left(t-\frac{1}{4}T\right)^2}
        +e^{2\beta\left(t-\frac{3}{4}T\right)^2}+1$ 
        & $\frac{3}{2}e^{\beta\left(t-\frac{1}{2}T\right)^2}+1$\\
        \refstepcounter{subsecvalnext}(\thesubsecvalnext)\label{case_swe_2d}  
        & $e^{\beta\left(t-0.3T\right)^2}
        +\frac{3}{2}e^{2\beta\left(t-0.7T\right)^2}+1$ 
        & $3\cos^2\left(\frac{10\pi}{T}t\right) + \frac{3}{4}$\\
        \refstepcounter{subsecvalnext}(\thesubsecvalnext)\label{case_swe_2e}  
        & $\frac{3}{2}e^{\beta\left(t-0.3T\right)^2}
        +e^{2\beta\left(t-0.7T\right)^2}+1$ 
        & $3\cos^2\left(\frac{10\pi}{T}t\right) + \frac{3}{4}$\\
        \refstepcounter{subsecvalnext}(\thesubsecvalnext)\label{case_swe_2f}  
        & $e^{4\beta\left(t-\frac{1}{4}T\right)^2}
        +\frac{3}{2}e^{4\beta\left(t-\frac{1}{2}T\right)^2}
        -\frac{1}{2}e^{4\beta\left(t-\frac{3}{4}T\right)^2}+1$
        & $3\cos^2\left(\frac{10\pi}{T}t\right) + \frac{3}{4}$\\
    \bottomrule
    \end{tabular}
    \caption{The true function for $p(t)$ denoted as $p_{true}$ and the corresponding initial guess used, $p^0$ with $\beta = -10,000$. Multiplicative noise is applied to $p^0$ in the simulations.}
    \label{tab:swe_p_1}
\end{table}

\begin{figure}[h]
    \begin{subfigure}[t]{0.32\textwidth}
        \centering
        \includegraphics[width=\linewidth]{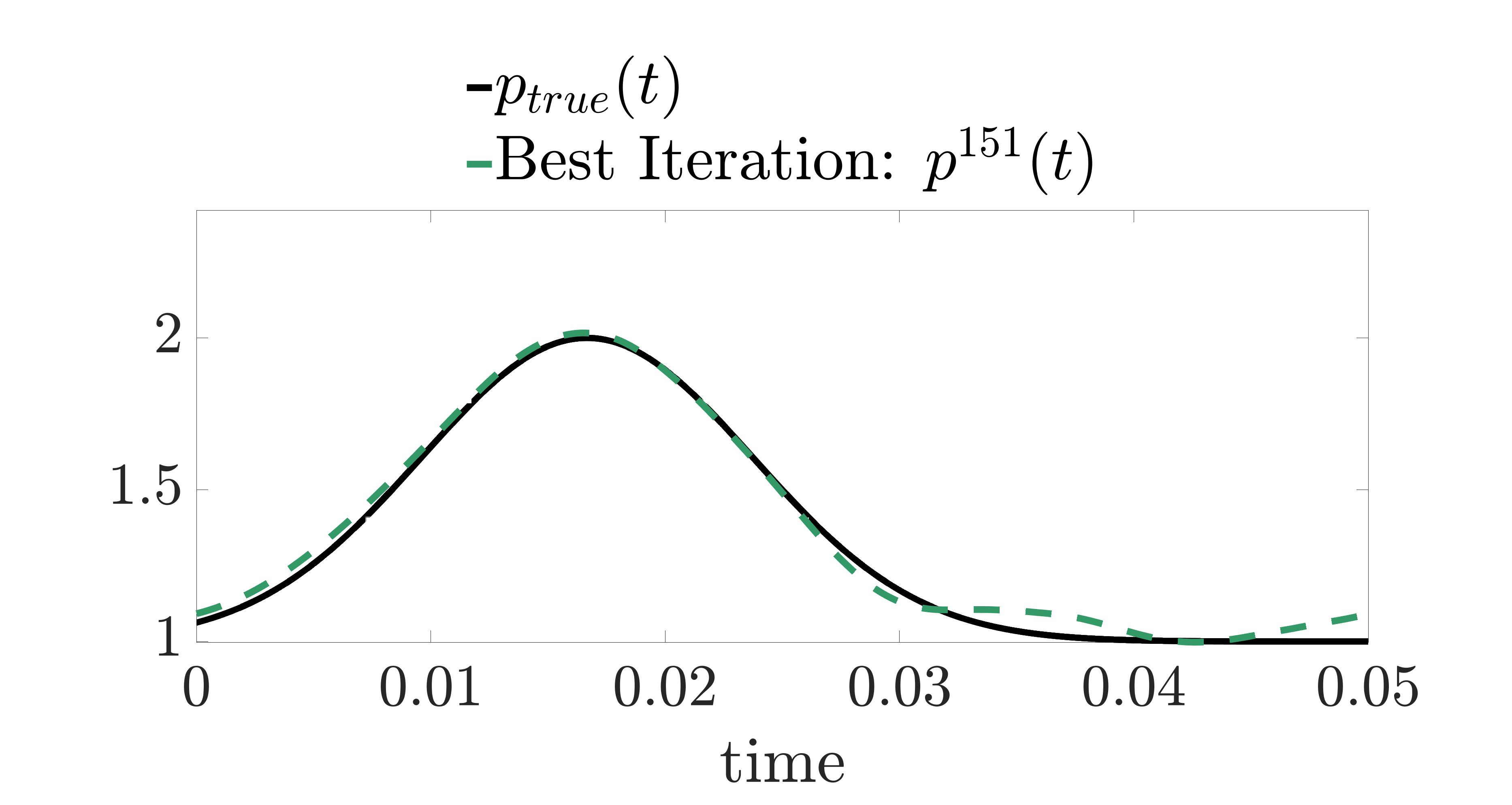}
    \end{subfigure}
    \begin{subfigure}[t]{0.32\textwidth}
        \centering
        \includegraphics[width=\linewidth]{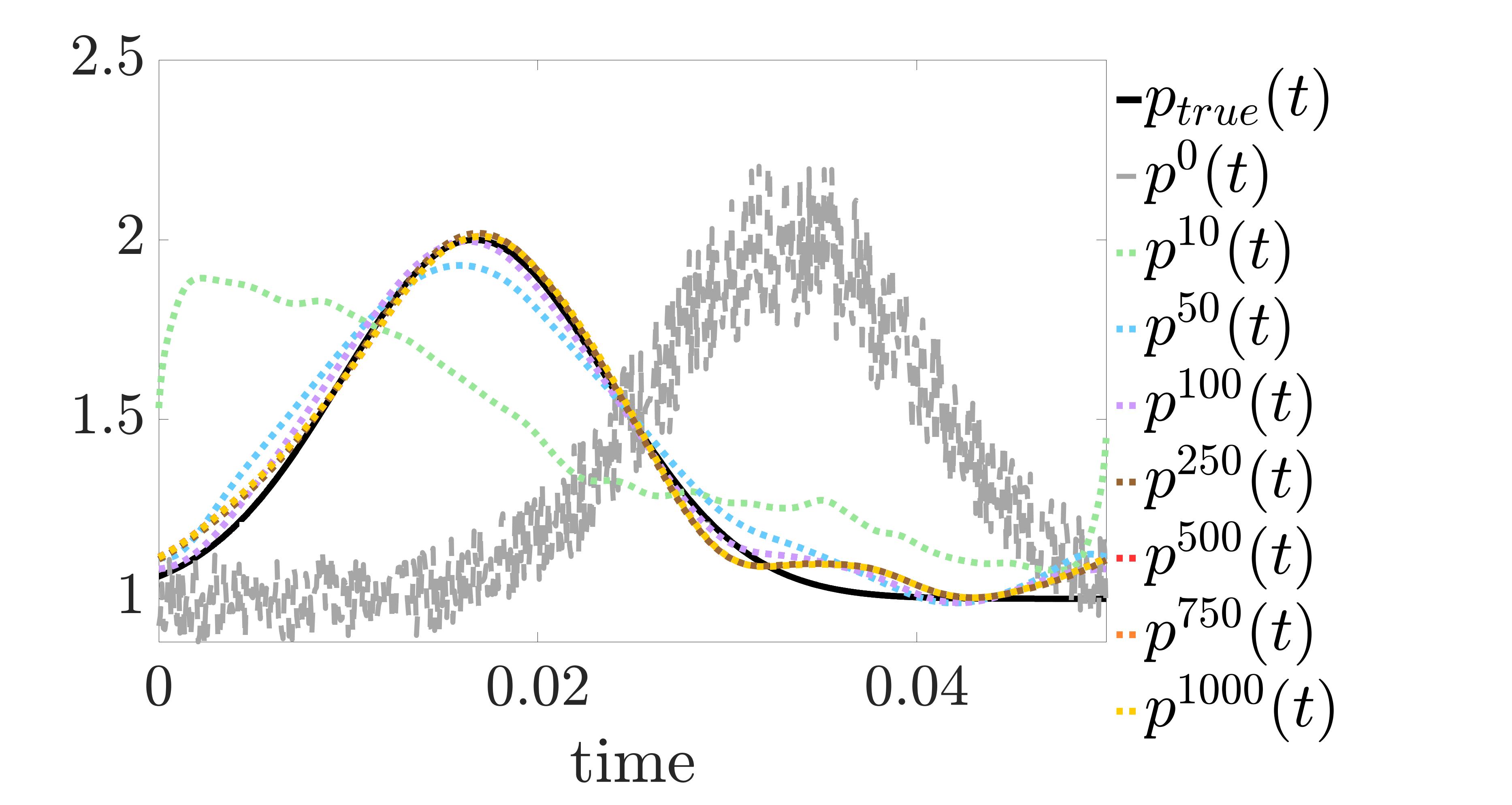}
    \end{subfigure}
    \begin{subfigure}[t]{0.32\textwidth}
        \includegraphics[width=\linewidth]{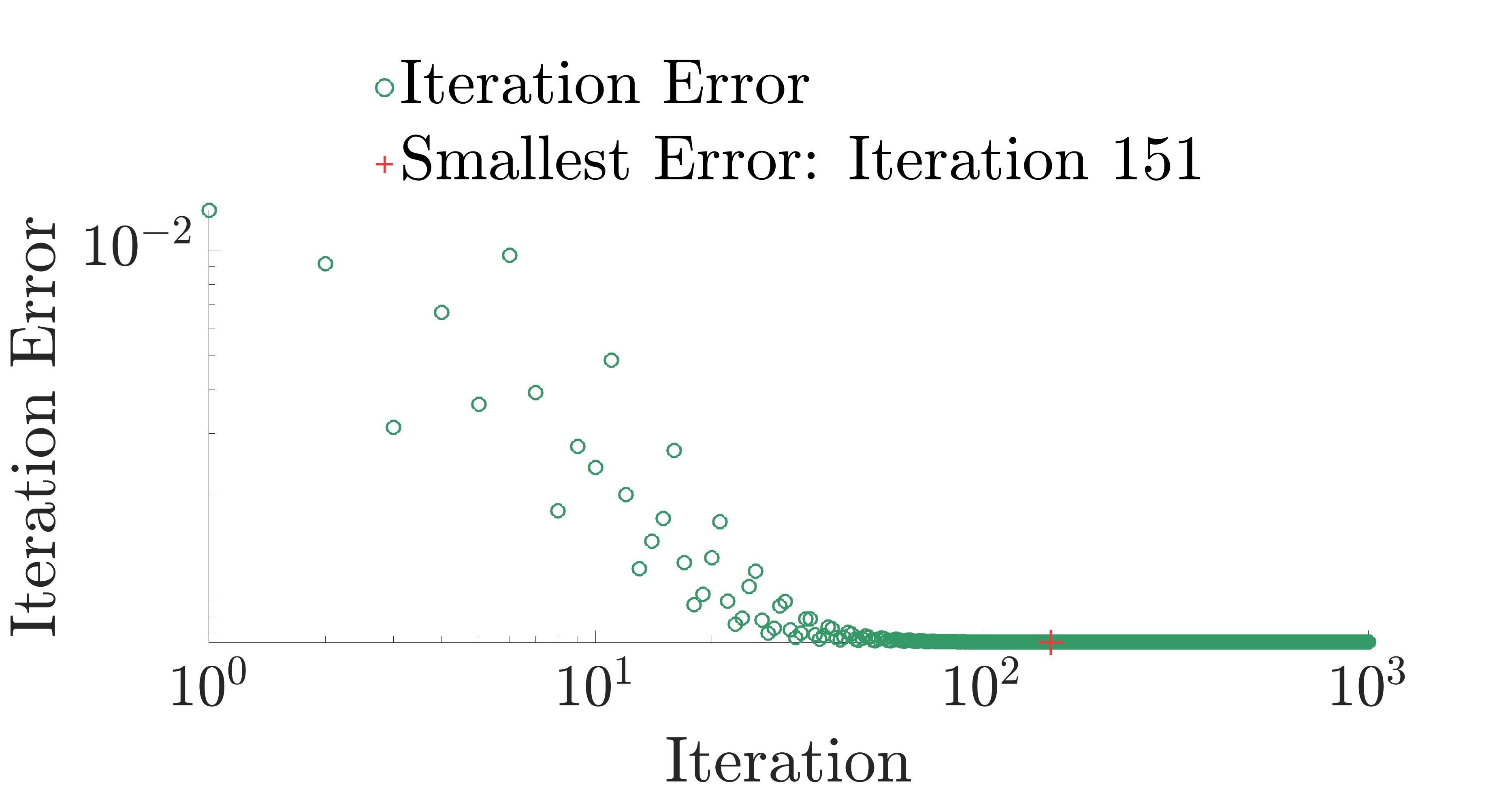}
    \end{subfigure}
    \caption{Results for Case \eqref{case_swe_2a}. Left: plots of the true $p$ and the numerically recovered $p$ at iteration 151, corresponding with the smallest residue; Middle: plots of the true $p$, the noisy initial guess, and various iteration values for $p$; Right: iteration errors on a log-log scale.}
    \label{fig:swe_2A}
\end{figure}
\begin{figure}[h]
    \begin{subfigure}[t]{0.32\textwidth}
        \centering
        \includegraphics[width=\linewidth]{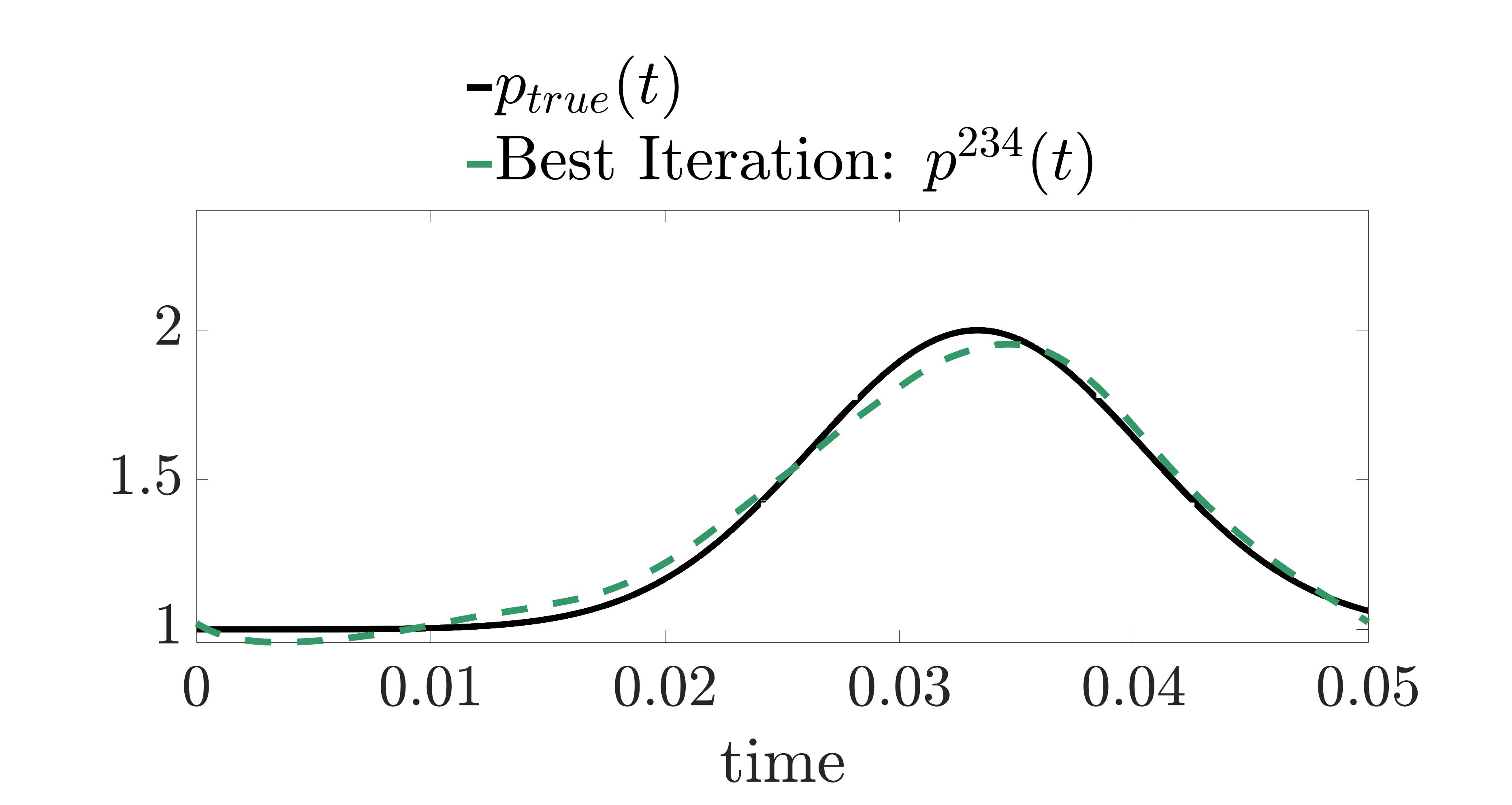}
    \end{subfigure}
    \begin{subfigure}[t]{0.32\textwidth}
        \centering
        \includegraphics[width=\linewidth]{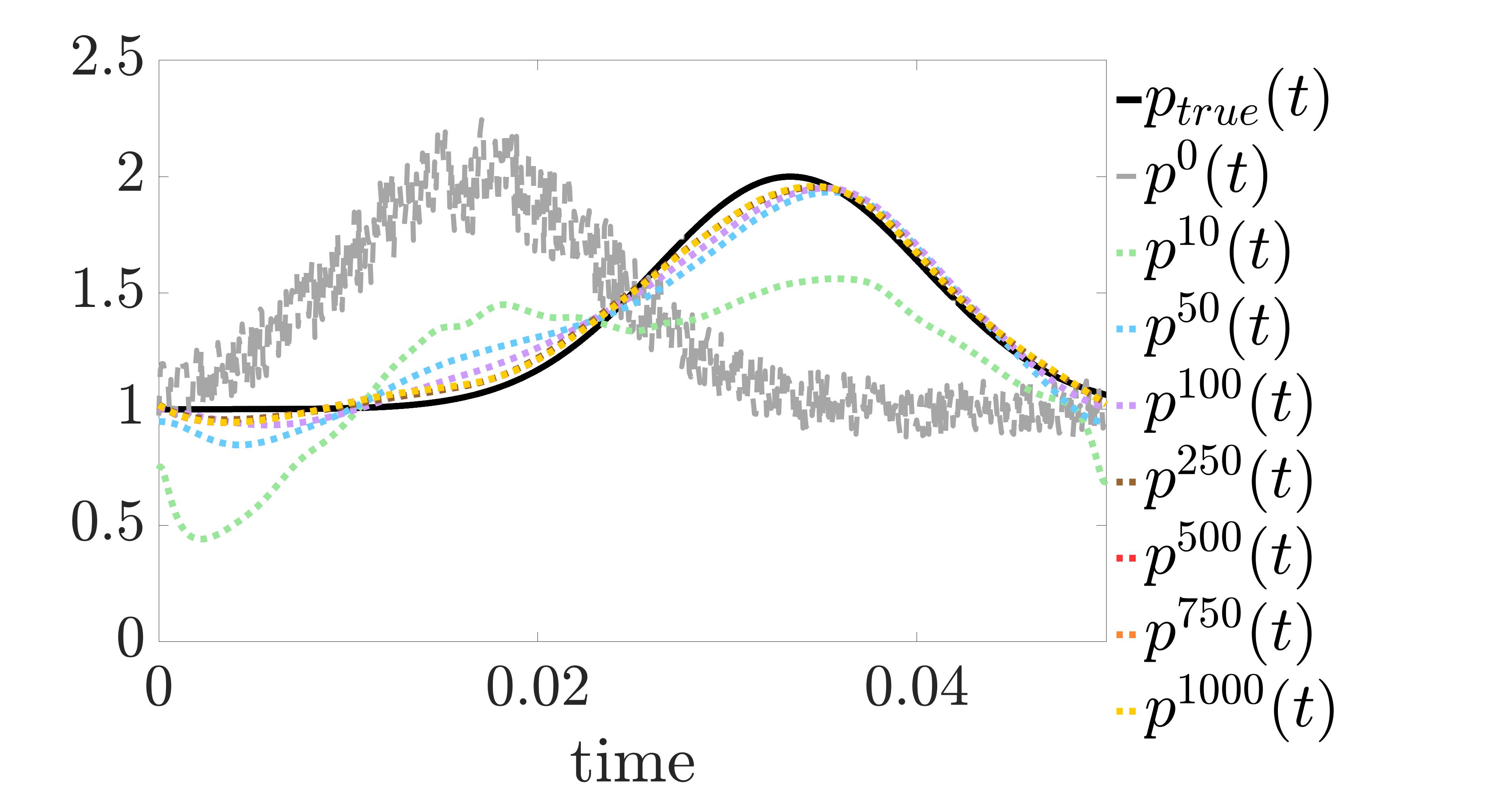}
    \end{subfigure}
    \begin{subfigure}[t]{0.32\textwidth}
        \includegraphics[width=\linewidth]{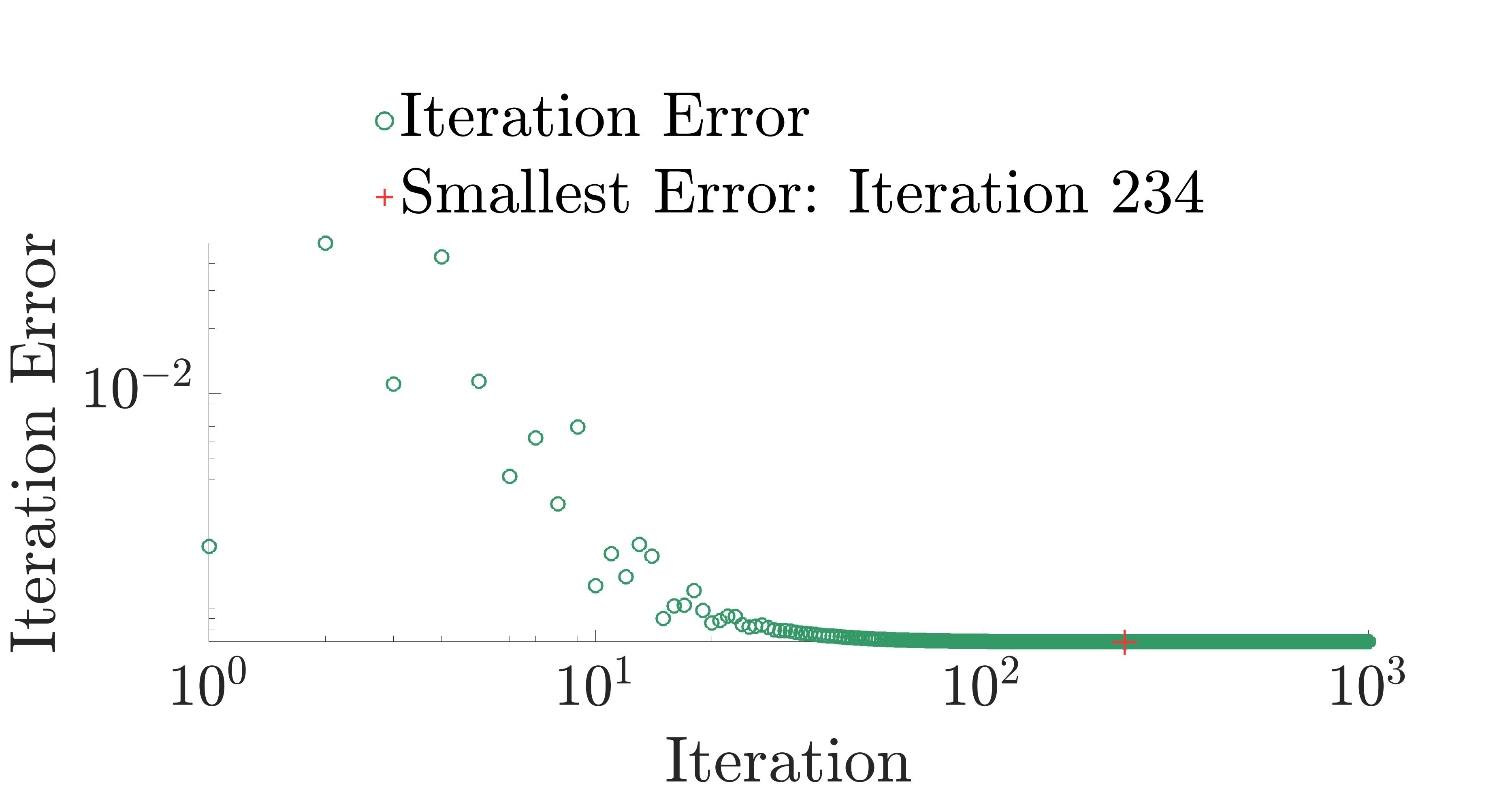}
    \end{subfigure}
    \caption{Results for Case \eqref{case_swe_2b}. 
    Left: plots of the true $p$ and the numerically recovered $p$ at iteration 234, corresponding with the smallest residue; Middle: plots of the true $p$, the noisy initial guess, and various iteration values for $p$;  Right: iteration errors on a log-log scale.
    }
    \label{fig:swe_2B}
\end{figure}

Cases \eqref{case_swe_2a} and \eqref{case_swe_2b} represent the situation in which the true value of $p(t)$ is a bump function that is non-centered with respect to the time interval and the corresponding initial guess is a noisy horizontal shift of $p_{true}(t)$. The numerical results are shown in Figures \ref{fig:swe_2A} and \ref{fig:swe_2B}. In both cases, the amplitude and shape of the true function and recovered numerical approximation are very close. The figures demonstrate that the scheme is robust even in the presence of multiplicative noise and the ill-posedness of the problem.

\begin{figure}[H]
    \begin{subfigure}[t]{0.32\textwidth}
        \centering
        \includegraphics[width=\linewidth]{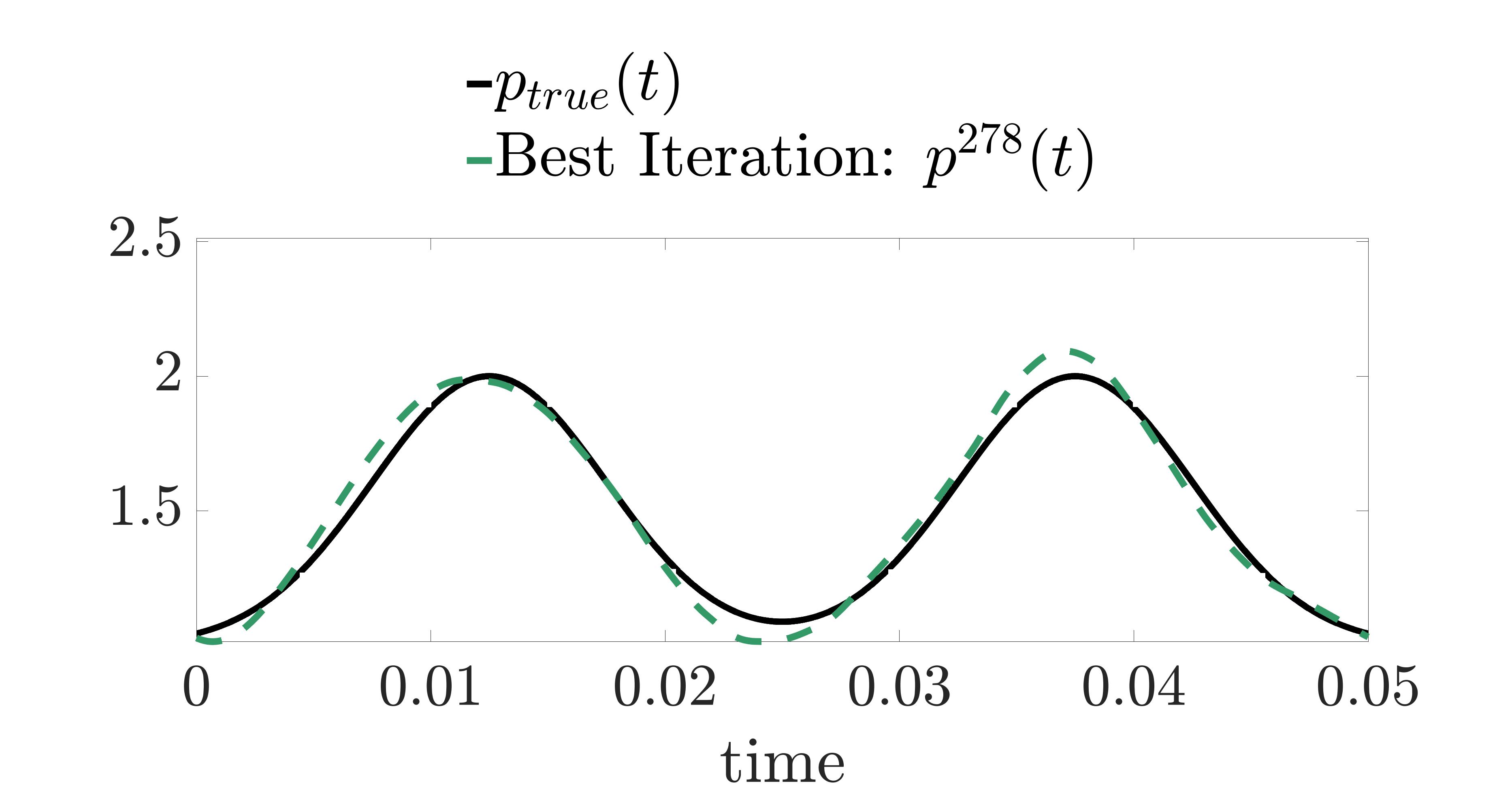}
    \end{subfigure}
    \begin{subfigure}[t]{0.32\textwidth}
        \centering
        \includegraphics[width=\linewidth]{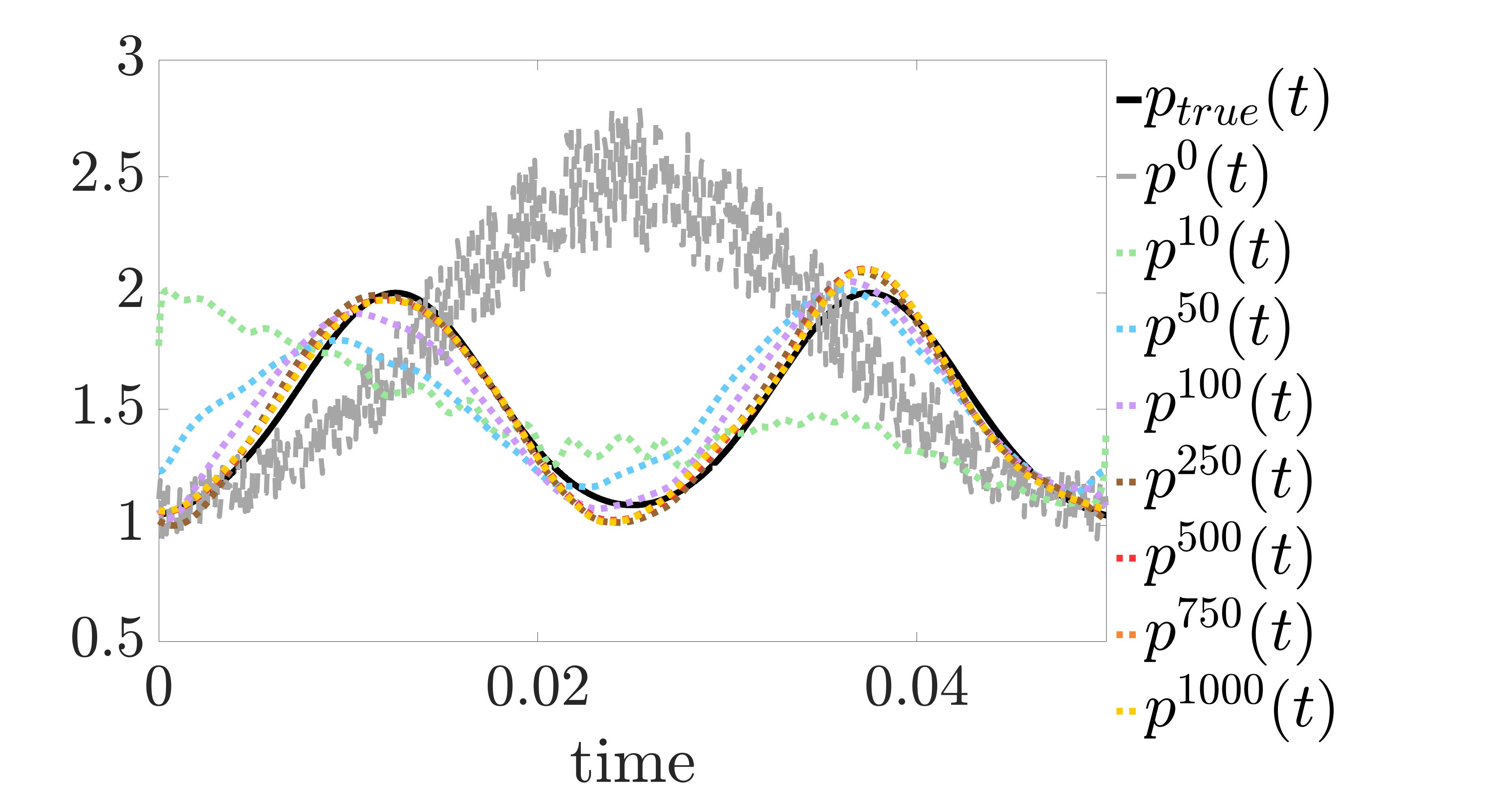}
    \end{subfigure}
    \begin{subfigure}[t]{0.32\textwidth}
        \includegraphics[width=\linewidth]{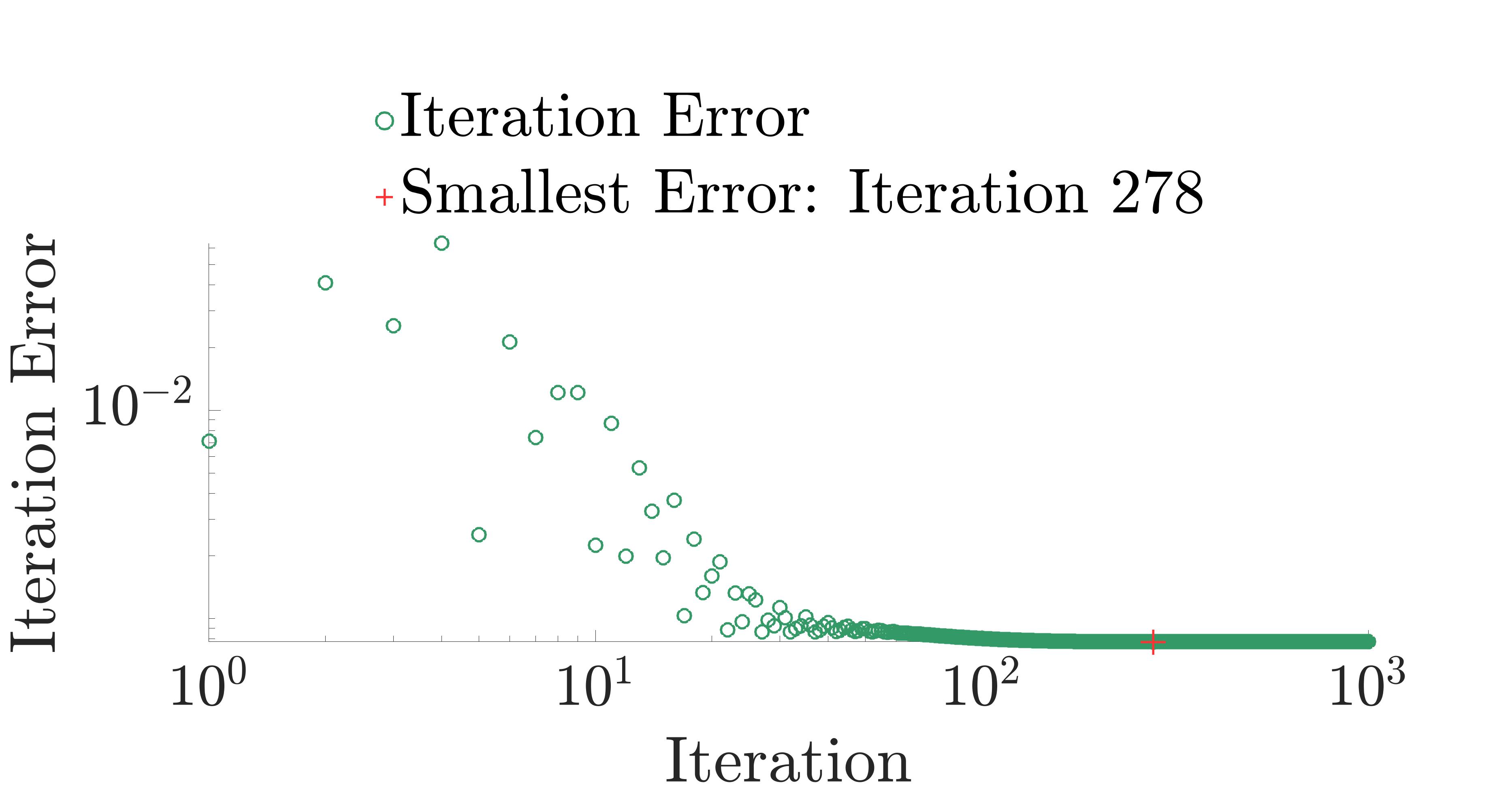}
    \end{subfigure}
    \caption{Results for Case \eqref{case_swe_2c}. Left: plots of the true $p$ and the numerically recovered $p$ at iteration 278, corresponding with the smallest residue; Middle: plots of the true $p$, the noisy initial guess, and various iteration values for $p$;  Right: iteration errors on a log-log scale.}
    \label{fig:swe_2C}
\end{figure}

\begin{figure}[h]
    \begin{subfigure}[t]{0.32\textwidth}
        \centering
        \includegraphics[width=\linewidth]{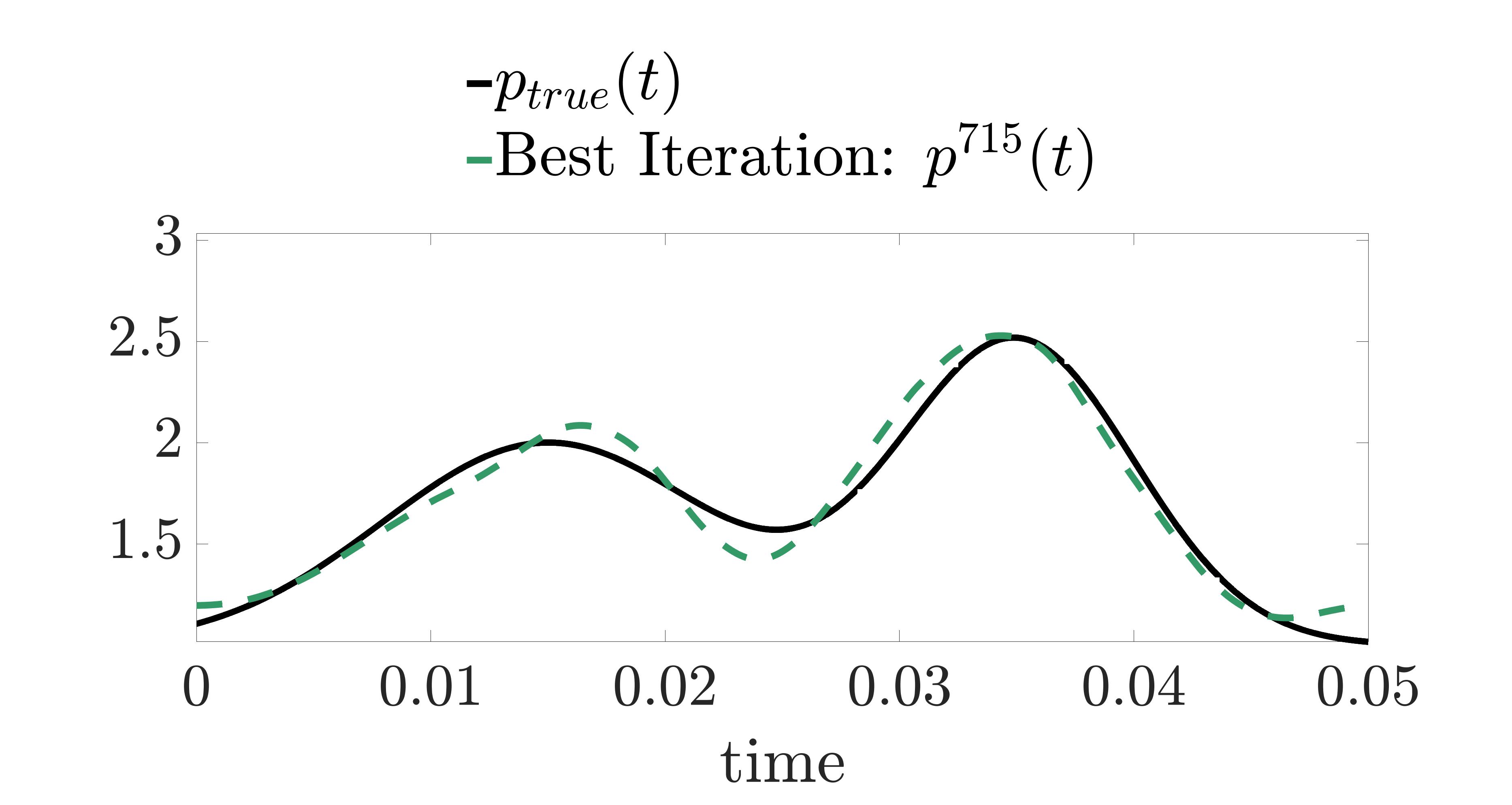}
    \end{subfigure}
    \begin{subfigure}[t]{0.32\textwidth}
        \centering
        \includegraphics[width=\linewidth]{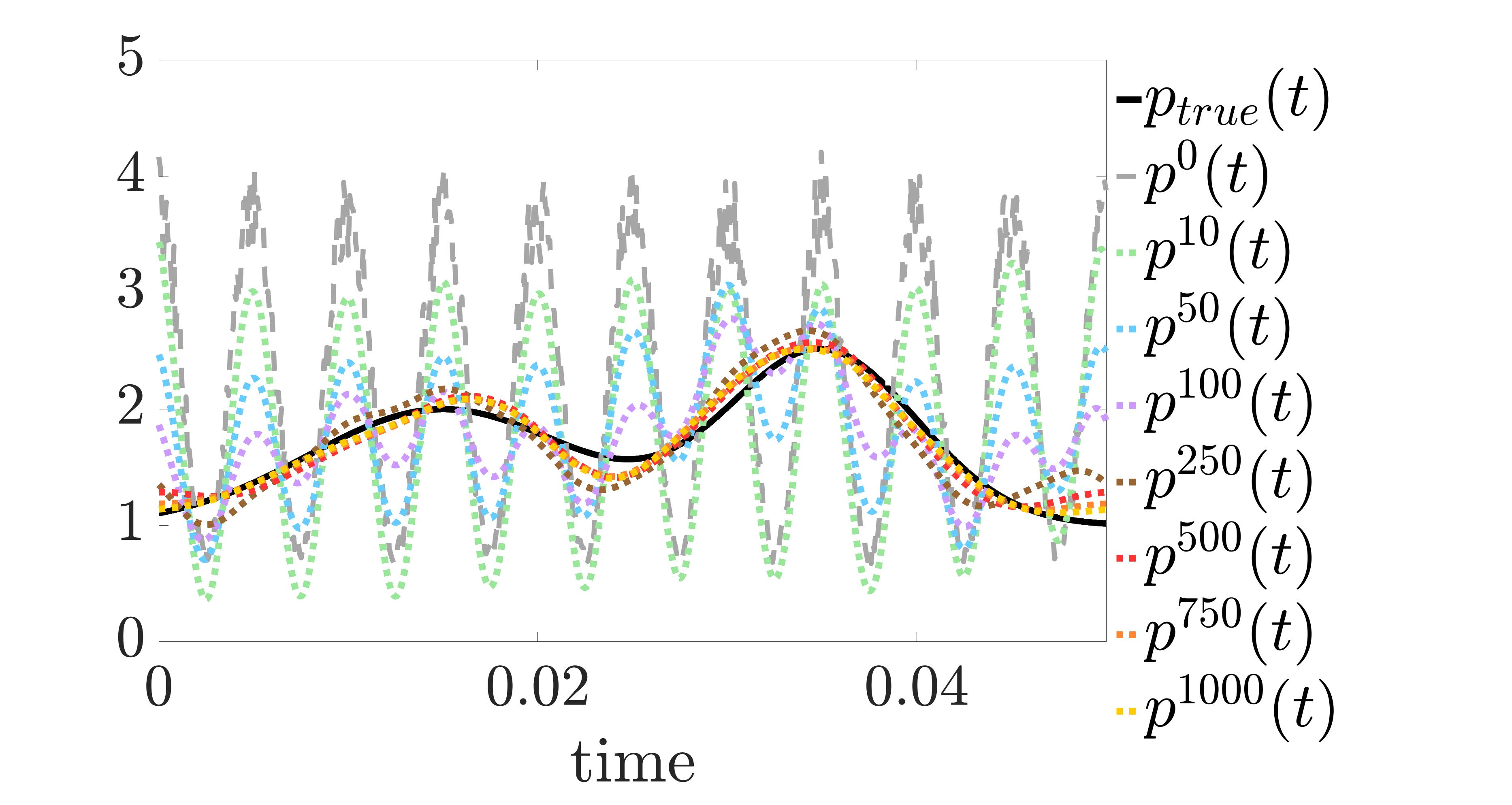}
        \caption{$\gamma_H=1\times 10^{-8}$}
    \end{subfigure}
    \begin{subfigure}[t]{0.32\textwidth}
        \includegraphics[width=\linewidth]{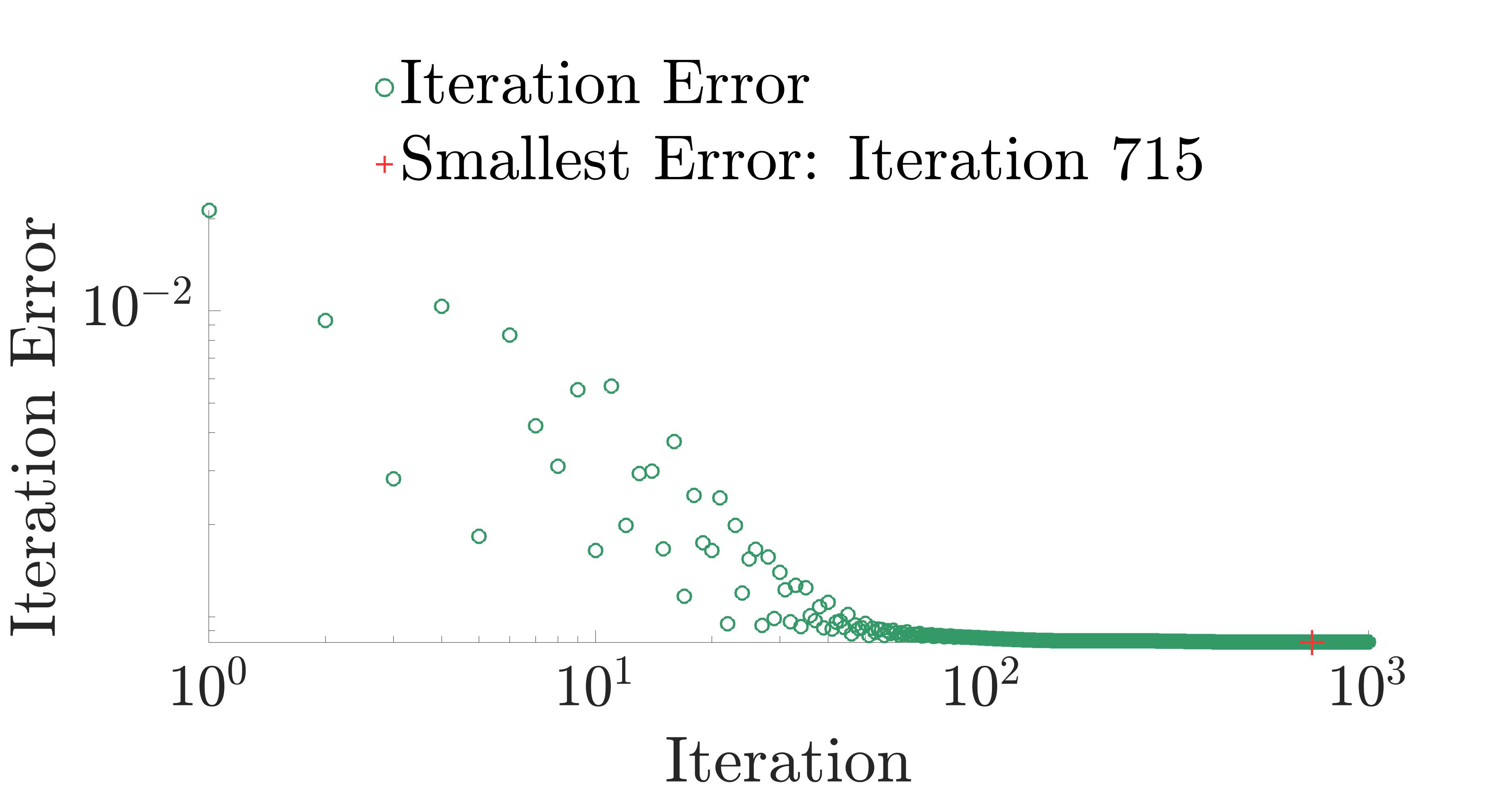}
    \end{subfigure} \\
    \begin{subfigure}[t]{0.32\textwidth}
        \centering
        \includegraphics[width=\linewidth]{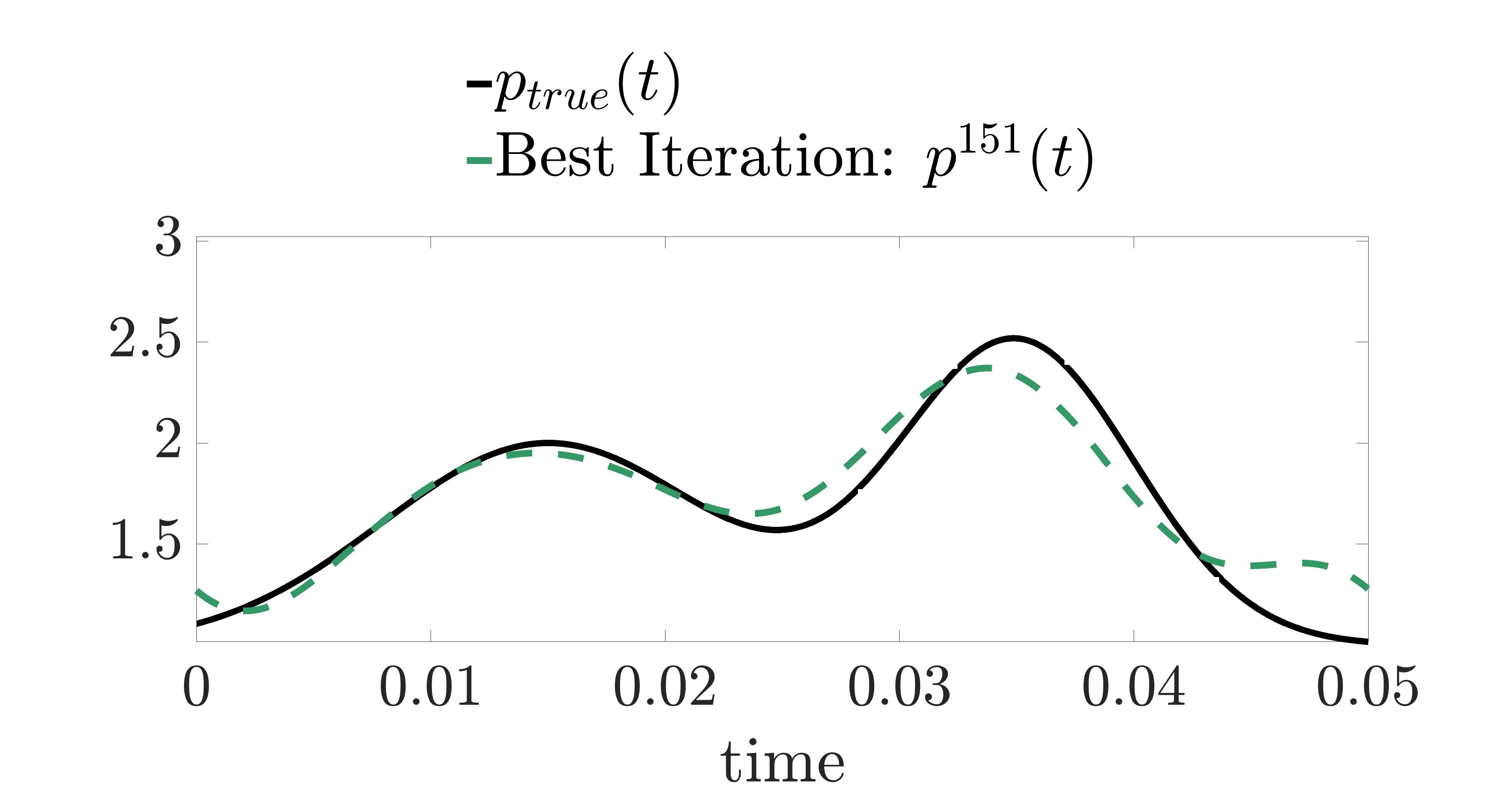}
    \end{subfigure}
    \begin{subfigure}[t]{0.32\textwidth}
        \centering
        \includegraphics[width=\linewidth]{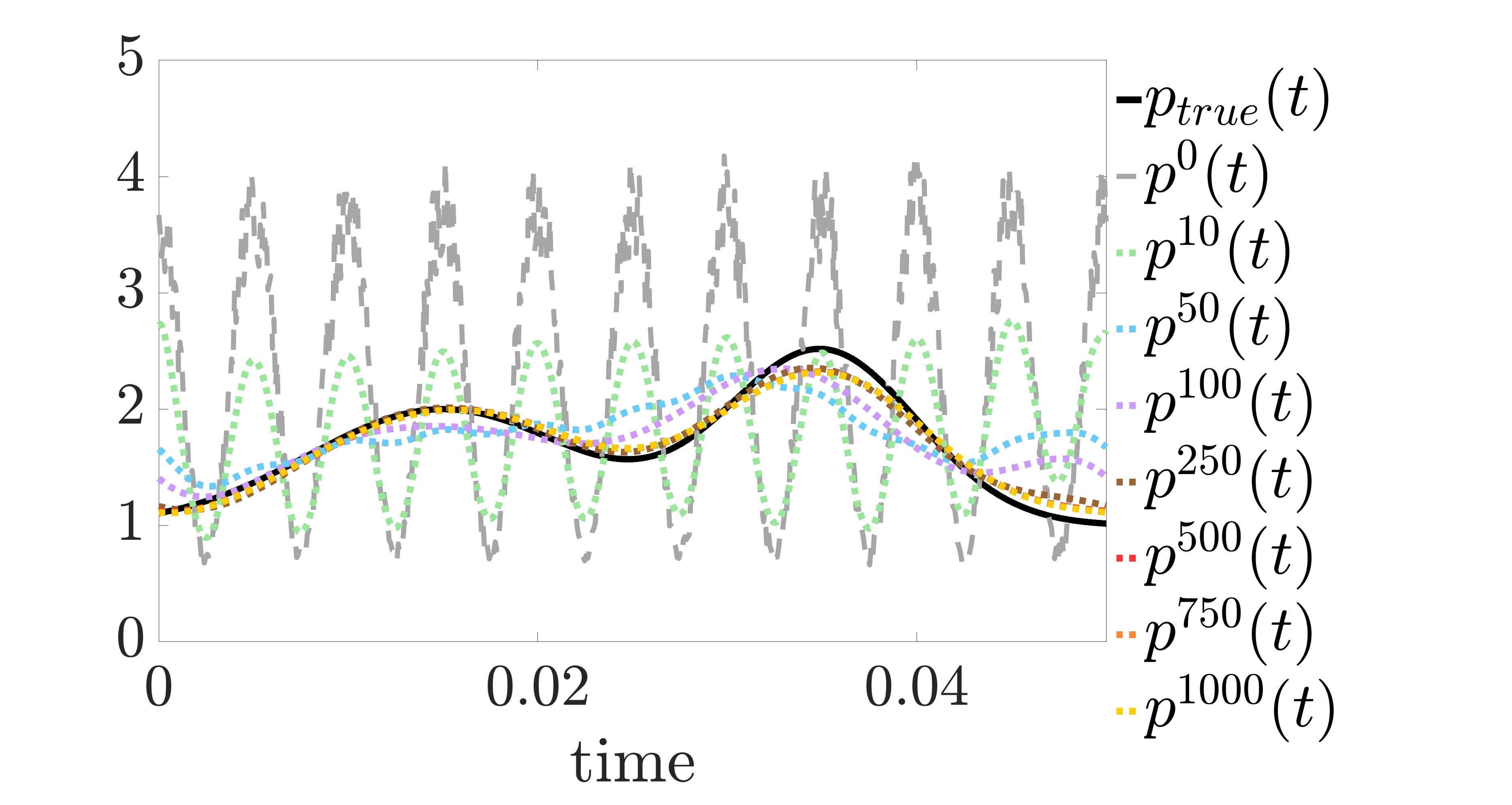}
        \caption{$\gamma_H=5\times 10^{-8}$}
        \label{fig:swe_2Db}
    \end{subfigure}
    \begin{subfigure}[t]{0.32\textwidth}
        \includegraphics[width=\linewidth]{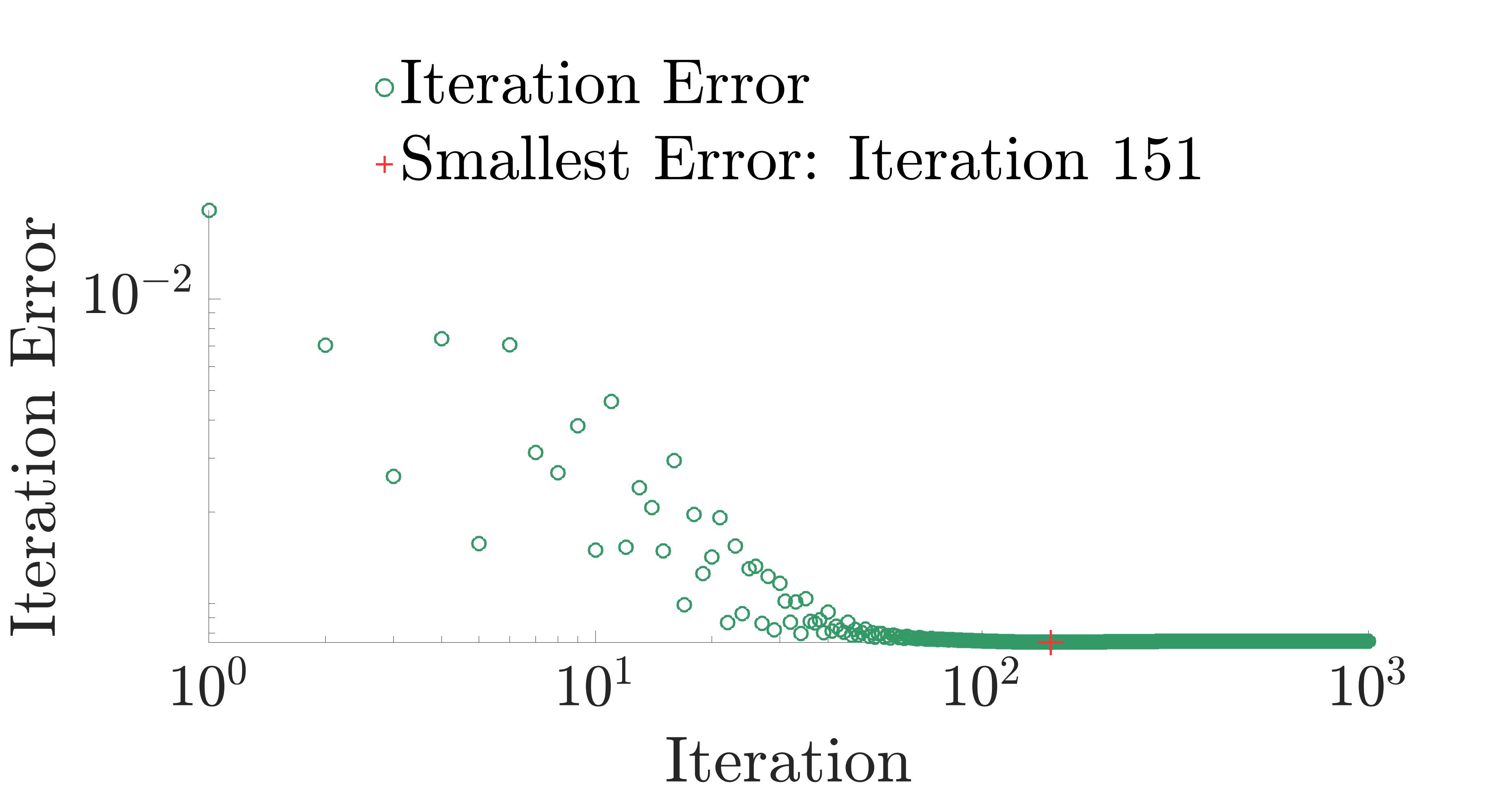}
    \end{subfigure}
    \caption{Results for Case \eqref{case_swe_2d}. 
    Left Column: plots of the true $p$ and the numerically recovered $p$ at iteration 715 (top row) and 151 (bottom row), corresponding with the smallest residue; Middle Column: plots of the true $p$, the noisy initial guess, and various iteration values for $p$; Right Column: iteration errors on a log-log scale.}
    \label{fig:swe_2D}
\end{figure}

\begin{figure}[H]
    \begin{subfigure}[t]{0.32\textwidth}
        \centering
        \includegraphics[width=\linewidth]{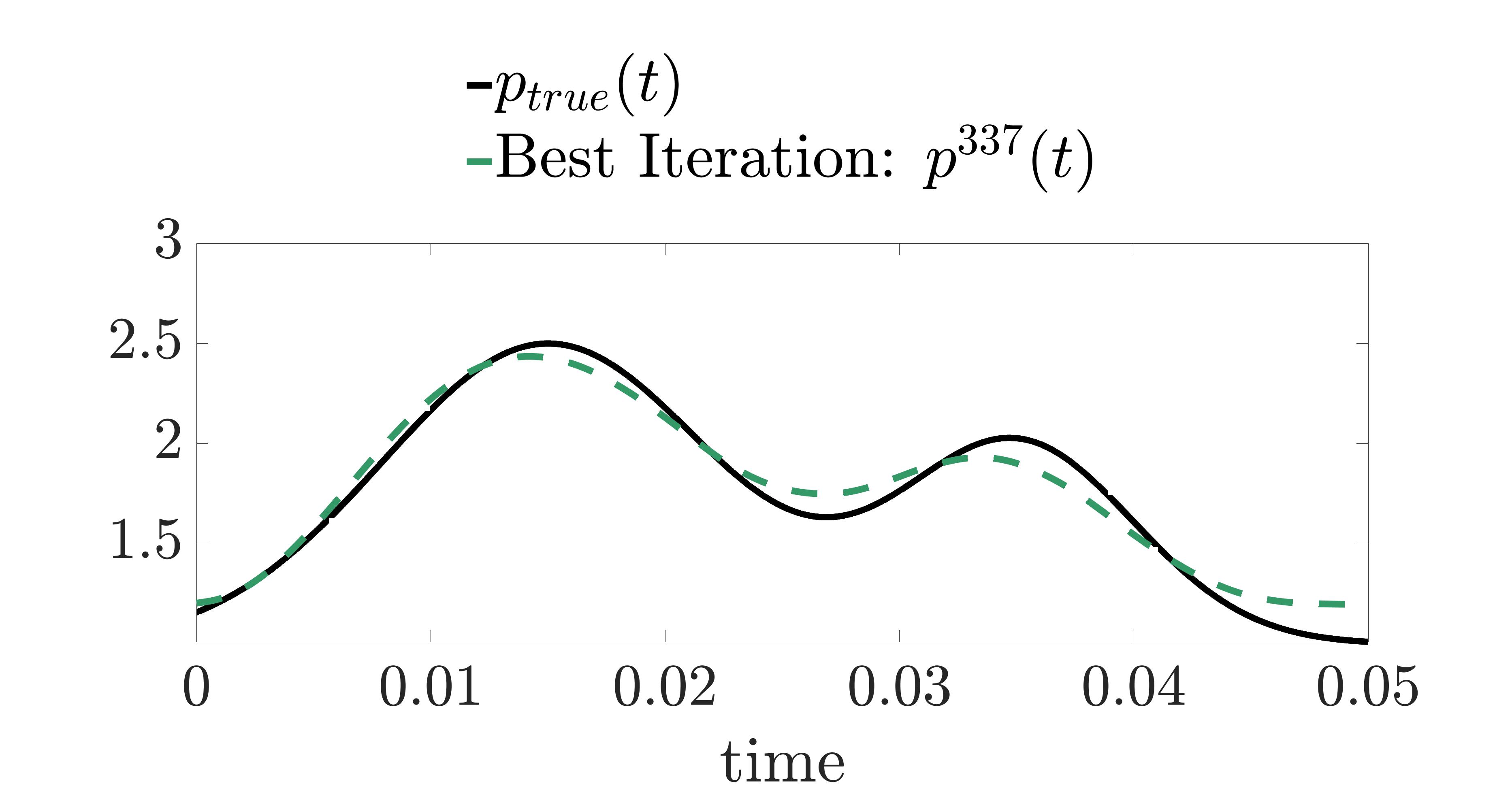}
    \end{subfigure}
    \begin{subfigure}[t]{0.32\textwidth}
        \centering
        \includegraphics[width=\linewidth]{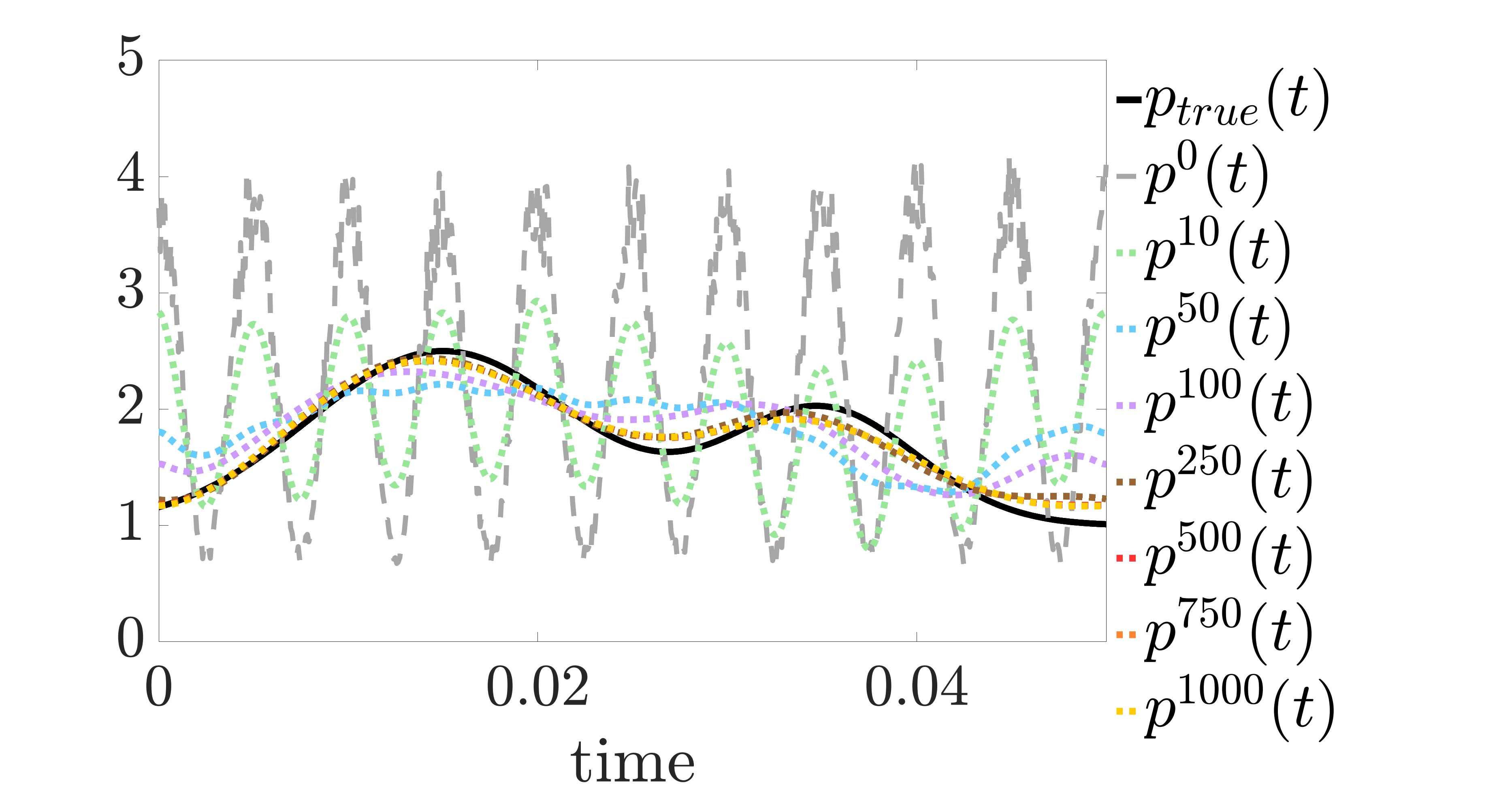}
    \end{subfigure}
    \begin{subfigure}[t]{0.32\textwidth}
        \includegraphics[width=\linewidth]{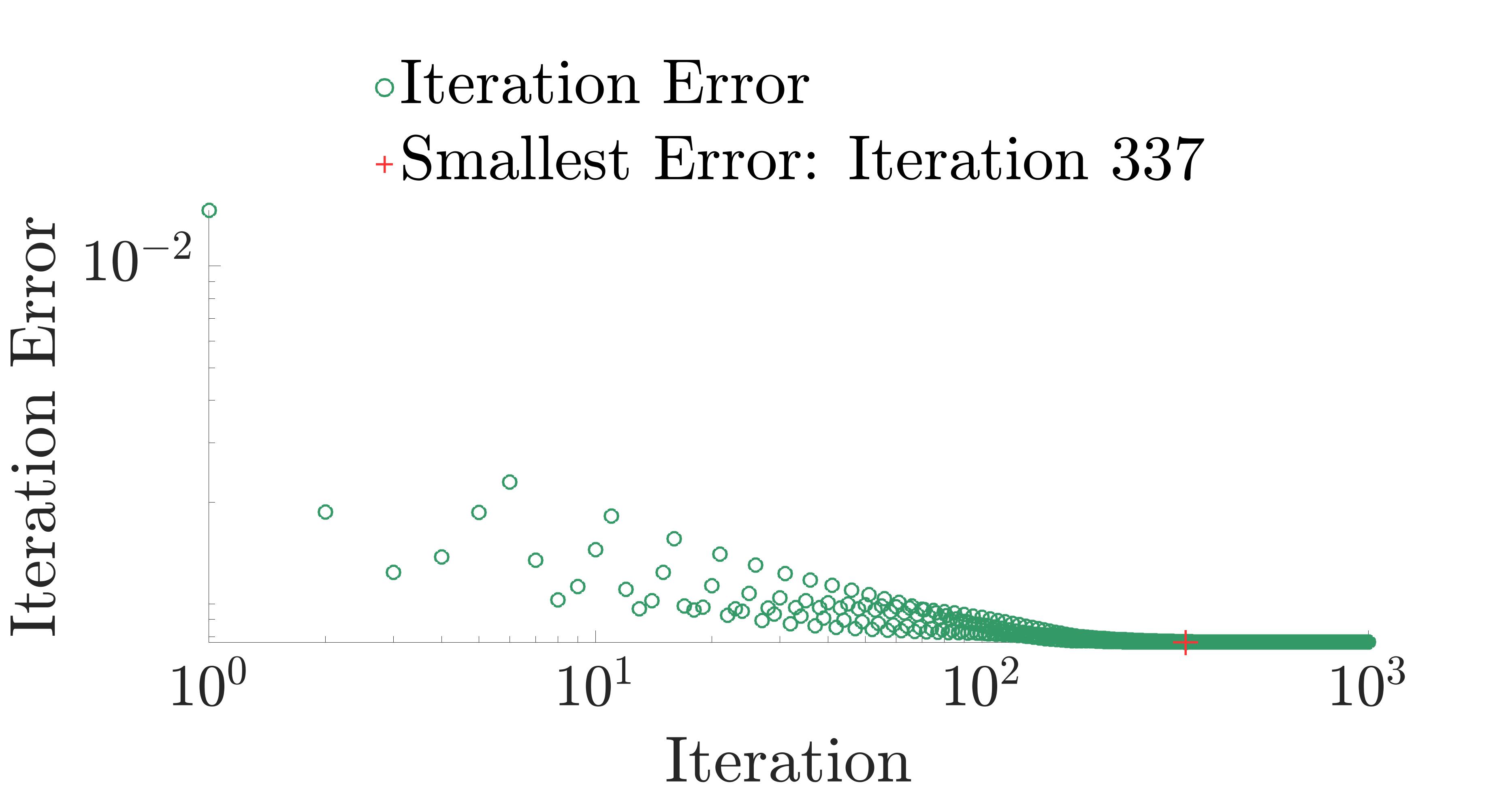}
    \end{subfigure}
    \caption{Results for Case \eqref{case_swe_2e}. 
    Left: plots of the true $p$ and the numerically recovered $p$ at iteration 337, corresponding with the smallest residue; Middle: plots of the true $p$, the noisy initial guess, and various iteration values for $p$; Right: iteration errors on a log-log scale.
    }
    \label{fig:swe_2E}
\end{figure}

In case \eqref{case_swe_2c}, we examine a true function $p$ that has two bumps of equal amplitude with an initial guess consisting of one bump with a larger amplitude. The results can be found in Figure \ref{fig:swe_2C}. Cases \eqref{case_swe_2d} and \eqref{case_swe_2e} also include a true $p$ function of two bumps, however they have different amplitudes and the corresponding initial guesses are highly oscillatory trigonometric functions. The corresponding results are shown in Figures \ref{fig:swe_2D} and \ref{fig:swe_2E}. In the examples with two bumps, the reconstructed function was also able to identify the two crests. The effect of the parameter $\gamma_H$ is explored in these cases. The value $\gamma_H = 1 \times 10^{-8}$ is used for Cases 
\eqref{case_swe_2c} and \eqref{case_swe_2d}. We can see that this smaller choice of $\gamma_H$ results in a less smooth solution in comparison with the results from Cases \eqref{case_swe_2d} and \eqref{case_swe_2e} when $\gamma_H = 5\times 10^{-8}$ is used. However, the plots corresponding to $\gamma_H = 5\times 10^{-8}$ while more smooth, are more flattened. 

\begin{figure}[h]
    \begin{subfigure}[t]{0.32\textwidth}
        \centering
        \includegraphics[width=\linewidth]{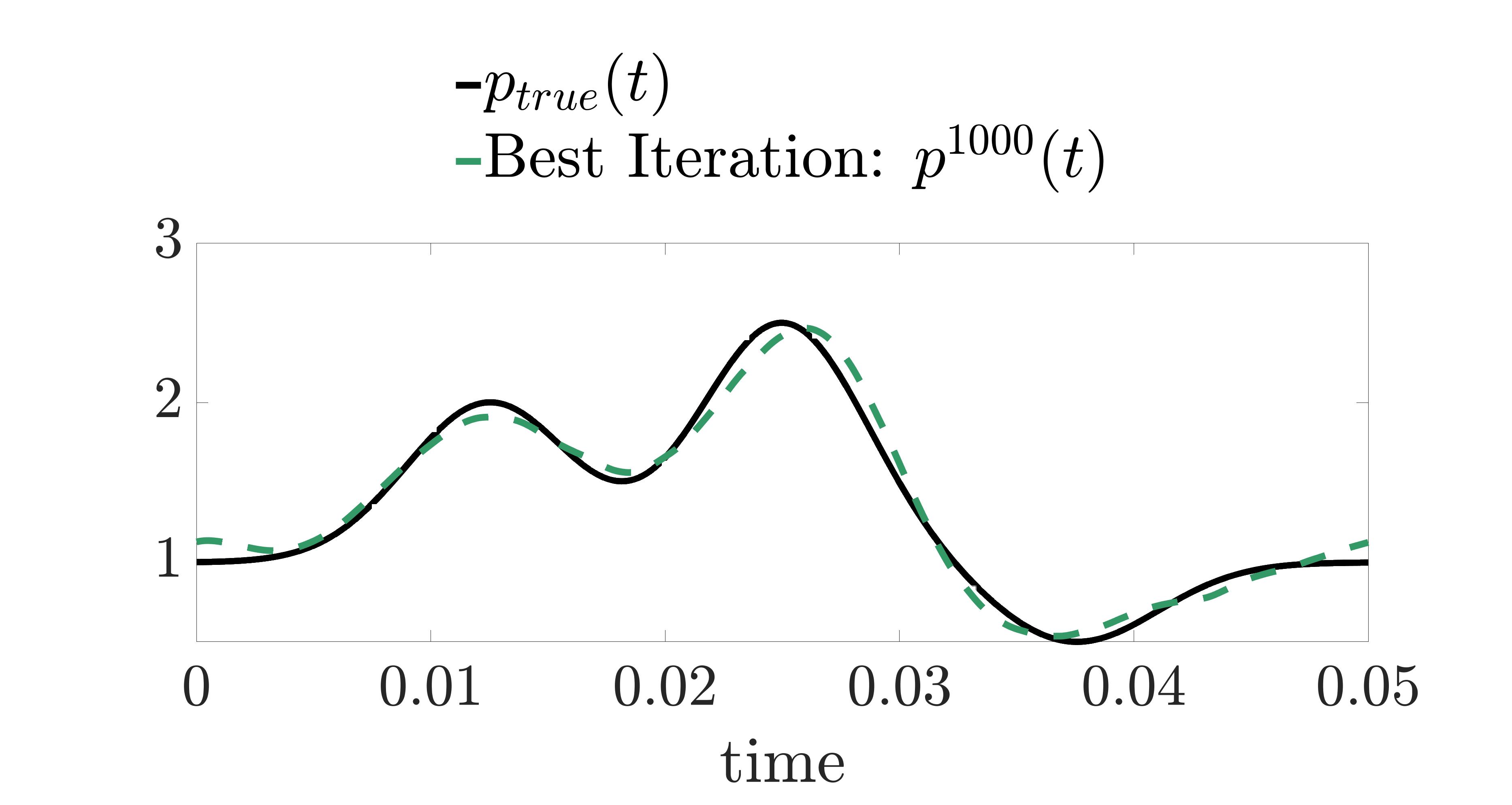}
    \end{subfigure}
    \begin{subfigure}[t]{0.32\textwidth}
        \centering
        \includegraphics[width=\linewidth]{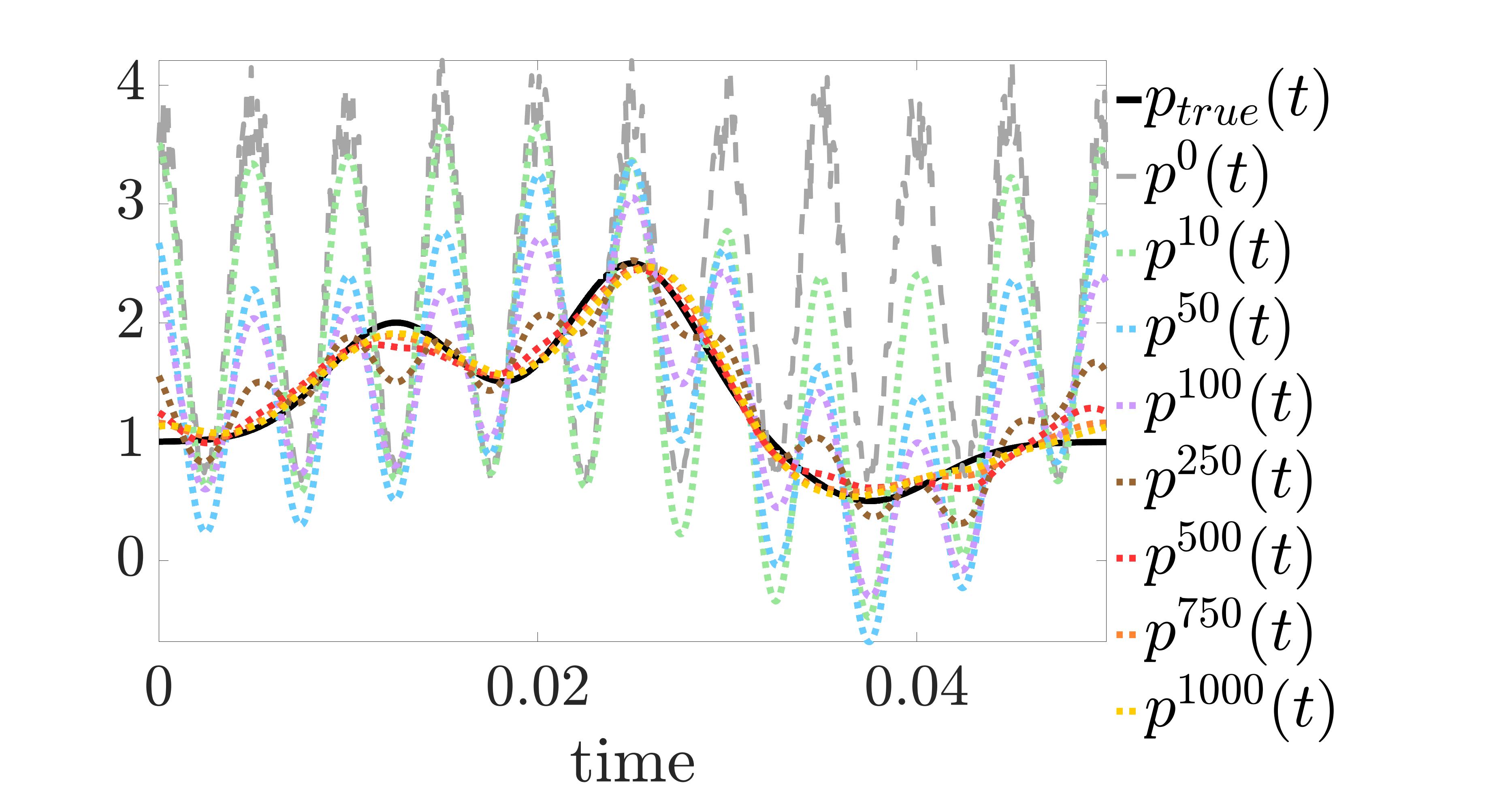}
    \end{subfigure}
    \begin{subfigure}[t]{0.32\textwidth}
        \includegraphics[width=\linewidth]{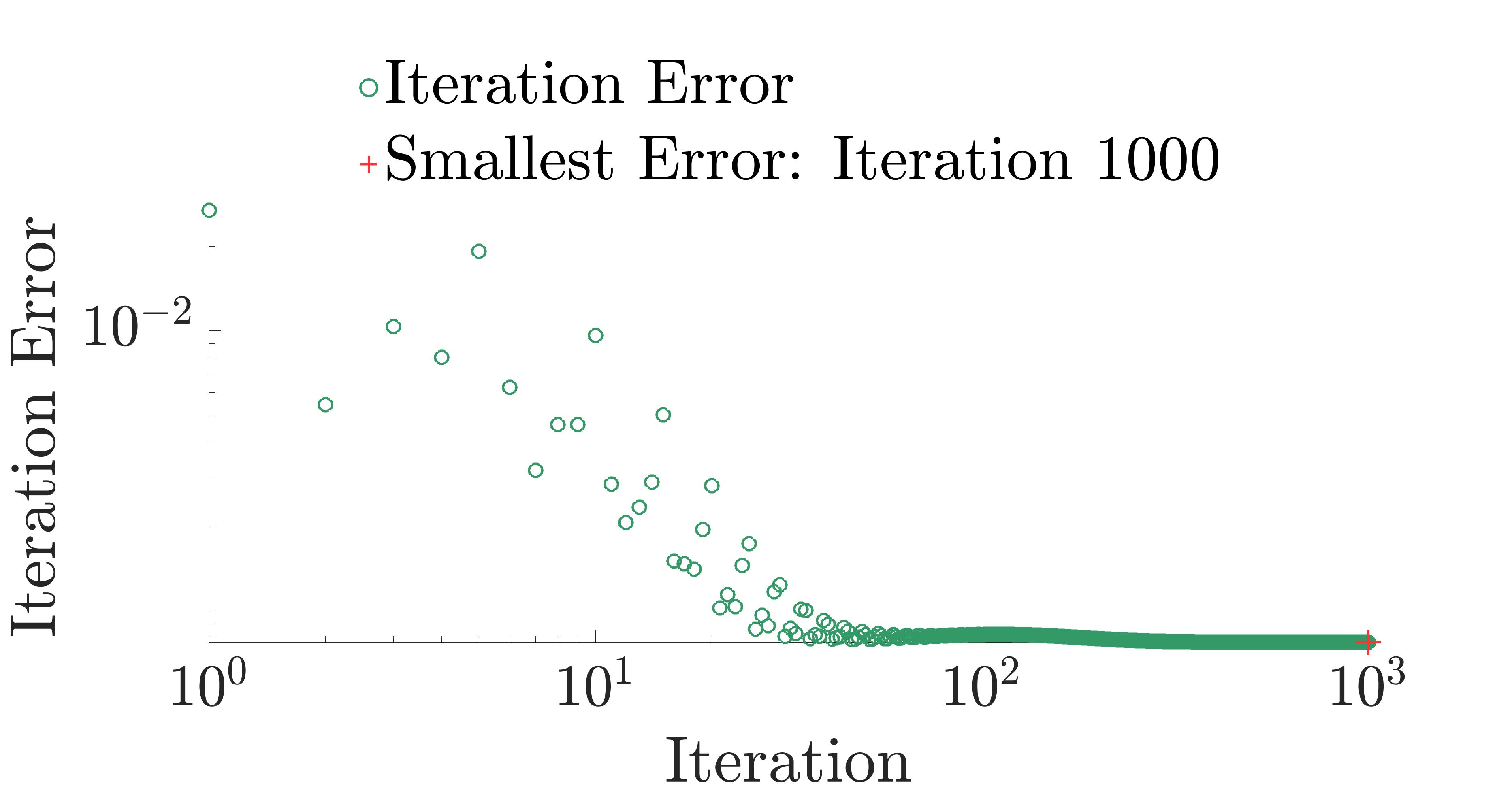}
    \end{subfigure}
    \caption{Results for Case \eqref{case_swe_2f}. Left: plots of the true $p$ and the numerically recovered $p$ at iteration 1000, corresponding with the smallest residue; Middle: plots of the true $p$, the noisy initial guess, and various iteration values for $p$; Right: iteration errors on a log-log scale.}
    \label{fig:swe_2F}
\end{figure}

\begin{figure}[h]
    \begin{subfigure}[t]{\textwidth}
        \centering
        \includegraphics[width=\linewidth]{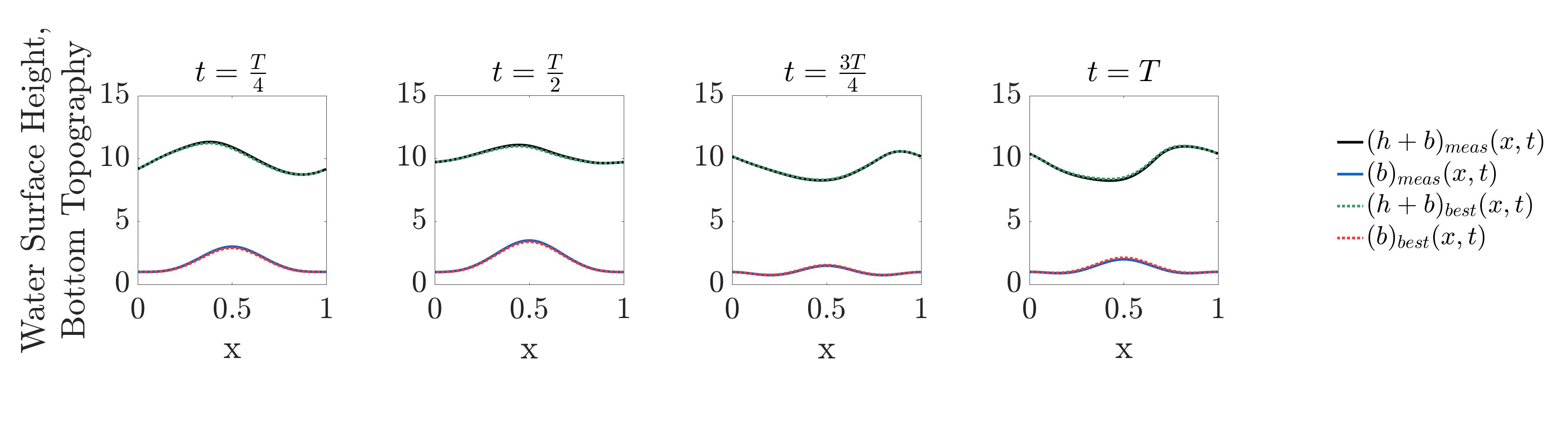}
    \end{subfigure}
    \begin{subfigure}[t]{\textwidth}
        \centering
        \includegraphics[width=\linewidth]{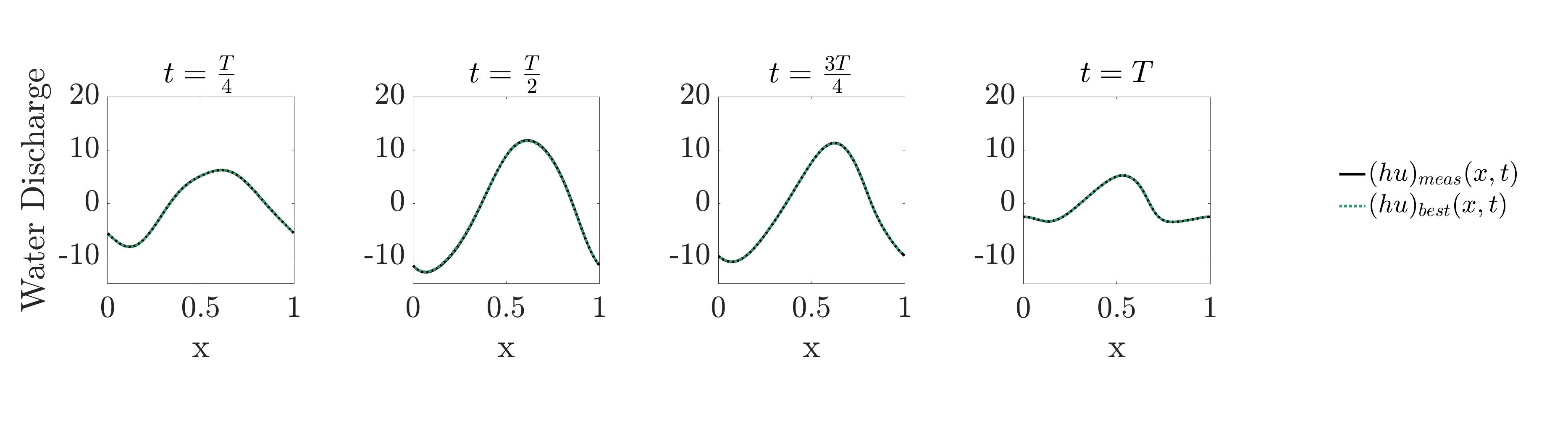}
    \end{subfigure}
    \caption{Comparison between the measured forward solutions and the results from the best iteration for SWEs Case \eqref{case_swe_2f} at times $t=\frac{T}{4},\frac{T}{2},\frac{3T}{4}, \text{ and } T$. In the top row the bottom topography function, $b$, and the water surface heights, $h+b$, are compared. In the bottom row, the water discharge, $hu$, is compared. The measured data and results from the iterative scheme are well matched for all functions in for each of the selected time snapshots.}
    \label{fig:swe_2F_2}
\end{figure}

The case \eqref{case_swe_2f} contains two crests of different amplitudes and a trough for the true function with a highly oscillatory trigonometric function as the initial guess. Plots of the results corresponding to this case can be found in Figure \ref{fig:swe_2F}. Additionally, for case \eqref{case_swe_2f} we show the solutions of the forward problem in Figure \ref{fig:swe_2F_2} at different times ($t=\frac{T}{4}, \frac{T}{2}, \frac{3T}{4}$, and $T$). The water surface height, bottom topography, and water discharge of the measured data and the numerical solution at the iteration with the smallest residue are compared. The recovered bottom topography along with the recovered state variables match the true functions well, even for some more complicated choices of $p(t)$.


Figures \ref{fig:swe_2A}-\ref{fig:swe_2F} each contain a plot of the residues $\J_0$, defined by \eqref{cost_funct} and \eqref{error_funct}, at the endpoints of the spatial domain for each iteration on a log-log scale. We see in each case a similar behavior occurs in which an `elbow'-like shape appears. The portion of this residue curve with a steeper slope corresponds to the situation when the term $\J_0$ has a greater impact on the update of the function $p$, which happens for the beginning iterations. The flat portion of the residue curve corresponds to the situation when the iteration starts to enter a small neighborhood where the regularization term $\R(p)$ convexifies the optimization problem and dominates the update. 

\subsection{Tests for Recovering \texorpdfstring{$p(t)$}{p(t)} from Different Initial Guesses} \label{sec_test_swe_2}
In this subsection, we run simulations with different initial guesses of $p^0(t)$ to recover the same $p_{true}(t)$. The goal is to demonstrate that the ability of our algorithm in recovering $p_{true}(t)$ does not depend on the initial guess. 

We consider the forward problem with the initial conditions in \eqref{ics_swe_smooth0} and the spatial bottom topography functions described in \eqref{bvals_swe_smooth0}. The true time component of the bottom topography function is fixed to be 
\begin{equation}\label{pcal_swe_2bumps}
    p(t) = \exp\left(\beta(t-0.3T)^2\right)+\frac{3}{2}\exp\left(2\beta(t-0.7T)^2\right)+1,
\end{equation}
with $\beta = -10,000$. Four different representative initial guesses, $p^0(t)$, listed in Table \ref{tab:swe_p_2}, will be tested. In all cases the final time is $T=0.05$ (while the solution is still smooth) and periodic boundary conditions are used. 
The remaining hyperparameters include a learning rate of $\ell=0.6$, $\gamma_L = 1 \times 10^{-6}$, and $\gamma_H = 1 \times 10^{-8}$.

\begin{table}[h]
    \centering
    \begin{tabular}{c c}
    \toprule
        Case & $p^0(t)$ \\
        \midrule
        \stepcounter{subsecval}\refstepcounter{subsecvalnext}(\thesubsecvalnext)\label{case_swe_3a} &  $1$\\
        \refstepcounter{subsecvalnext}(\thesubsecvalnext)\label{case_swe_3b} &  $4\sin^2\left(\frac{\pi}{T}t\right)$\\
        \refstepcounter{subsecvalnext}(\thesubsecvalnext)\label{case_swe_3c} &  $-2\sin^2\left(\frac{\pi}{T}t\right)+2$\\
        \refstepcounter{subsecvalnext}(\thesubsecvalnext)\label{case_swe_3d} &  $3\cos^2\left(\frac{10\pi}{T}t\right)+0.75$\\
    \bottomrule
    \end{tabular}
    \caption{The corresponding initial guesses used, $p^0$ with $T=0.05$. Multiplicative noise is applied to $p^0$ in the simulations.}
    \label{tab:swe_p_2}
\end{table}

\begin{figure}[h]
    \begin{subfigure}[t]{0.32\textwidth}
        \centering
        \includegraphics[width=\linewidth]{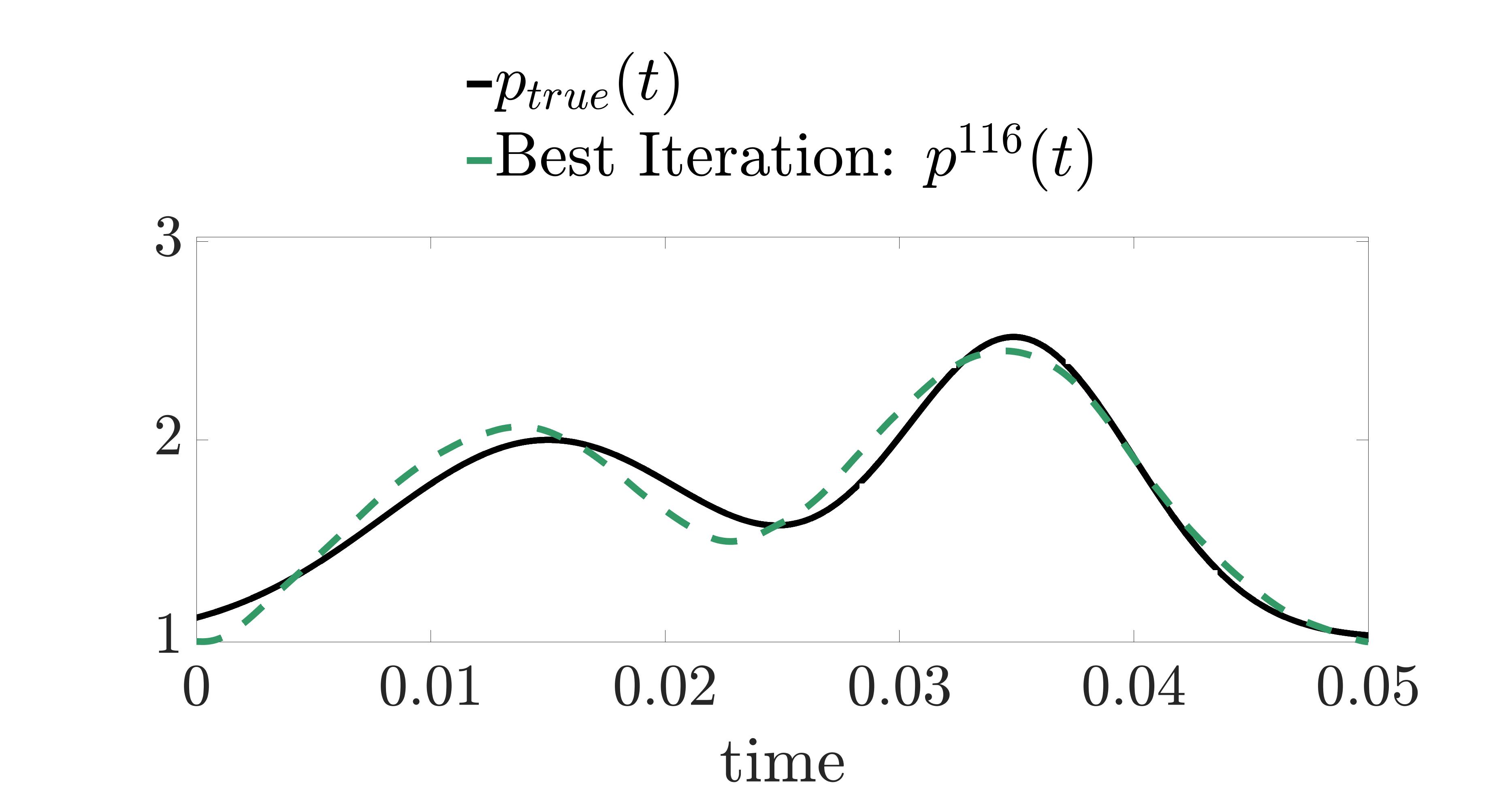}
    \end{subfigure}
    \begin{subfigure}[t]{0.32\textwidth}
        \centering
        \includegraphics[width=\linewidth]{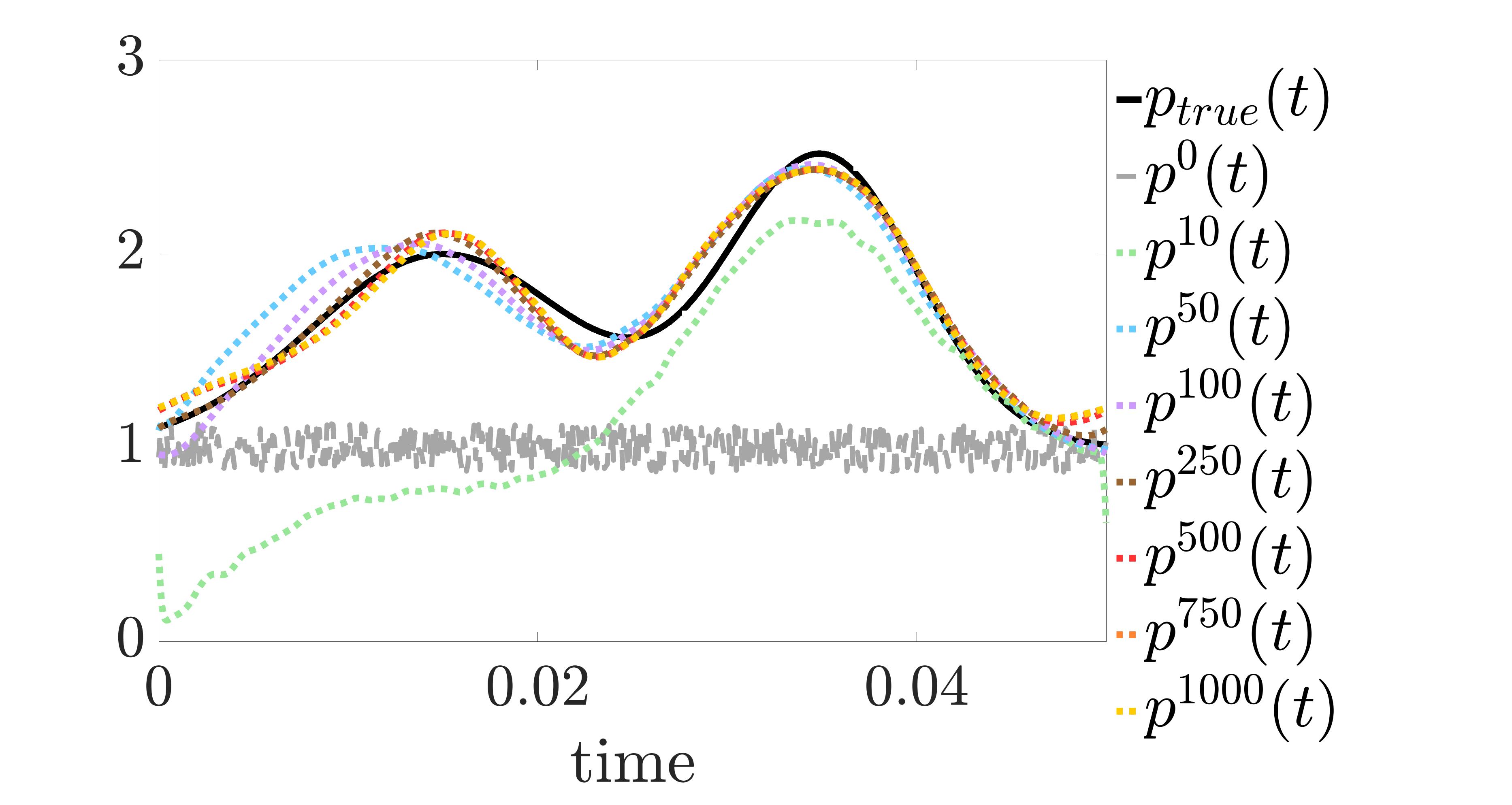}
    \end{subfigure}
    \begin{subfigure}[t]{0.32\textwidth}
        \includegraphics[width=\linewidth]{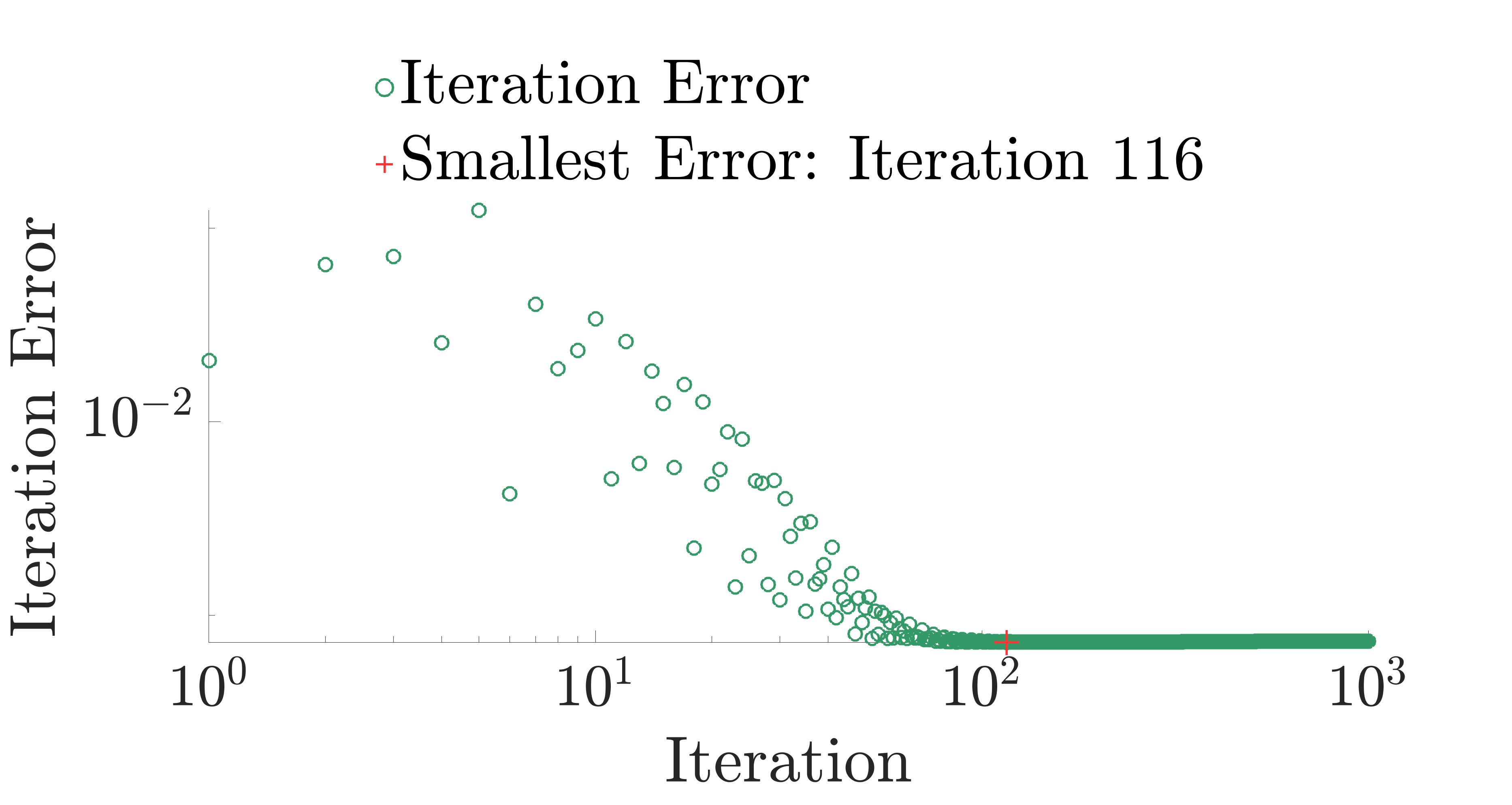}
    \end{subfigure}
    \caption{Results for Case \eqref{case_swe_3a}. Left: plots of the true $p$ and the $p$ corresponding with the smallest residue error, at iteration 116; Middle: plots of the true $p$, the noisy initial guess, and various iteration values for $p$; Right: iteration errors on a log-log scale.}
    \label{fig:swe_3A}
\end{figure}

\begin{figure}[h]
    \begin{subfigure}[t]{0.32\textwidth}
        \centering
        \includegraphics[width=\linewidth]{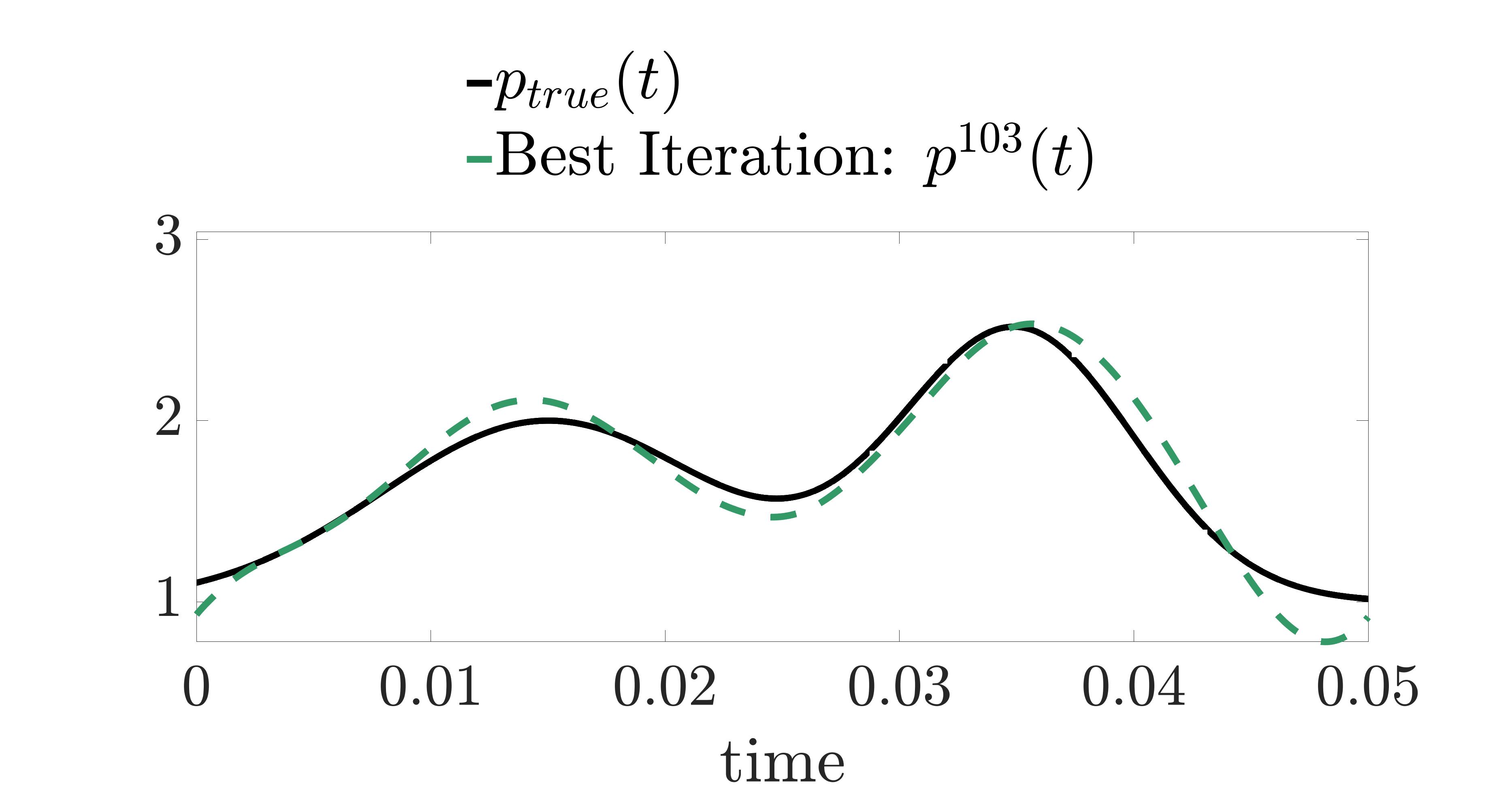}
    \end{subfigure}
    \begin{subfigure}[t]{0.32\textwidth}
        \centering
        \includegraphics[width=\linewidth]{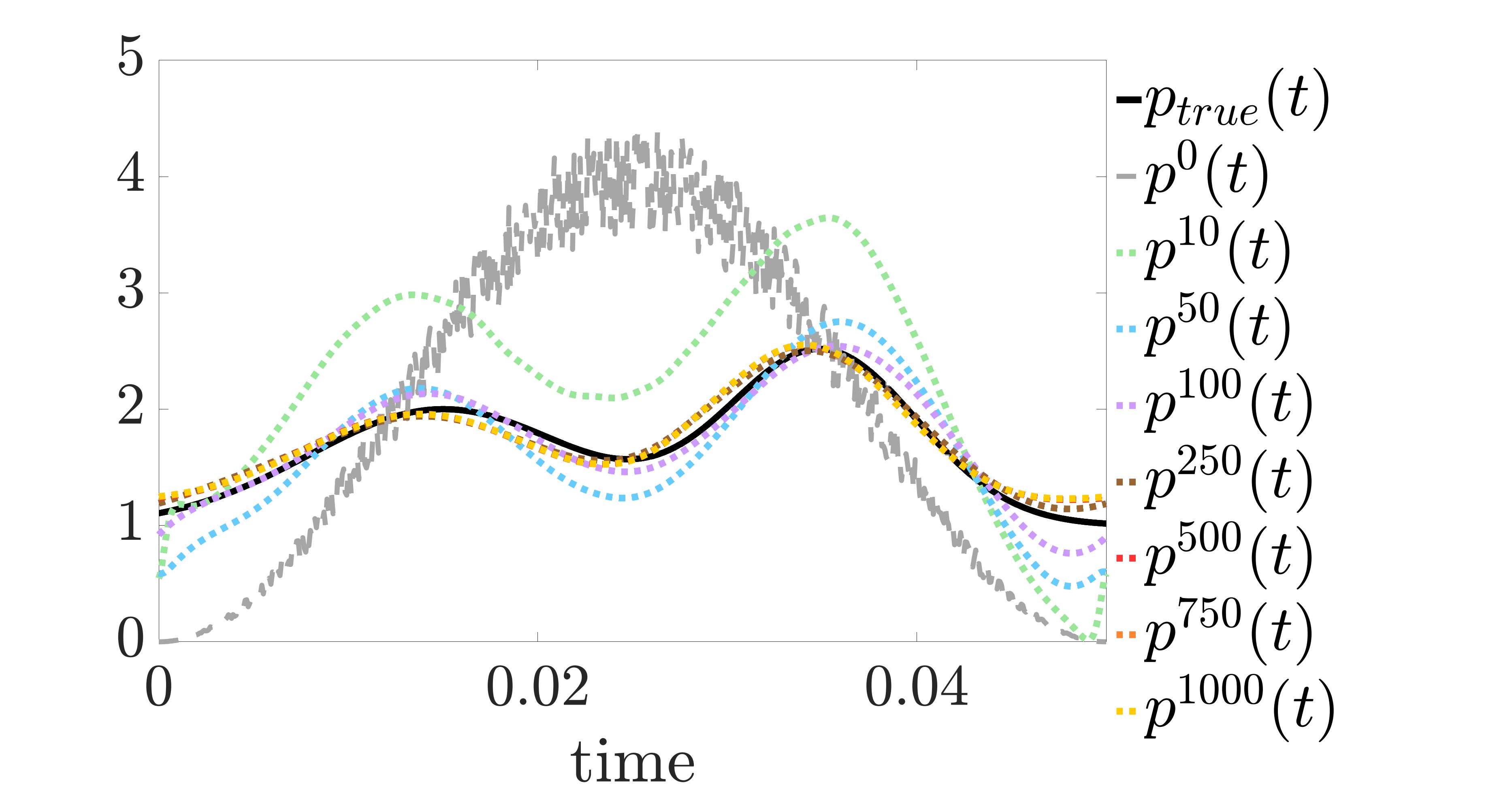}
    \end{subfigure}
    \begin{subfigure}[t]{0.32\textwidth}
        \includegraphics[width=\linewidth]{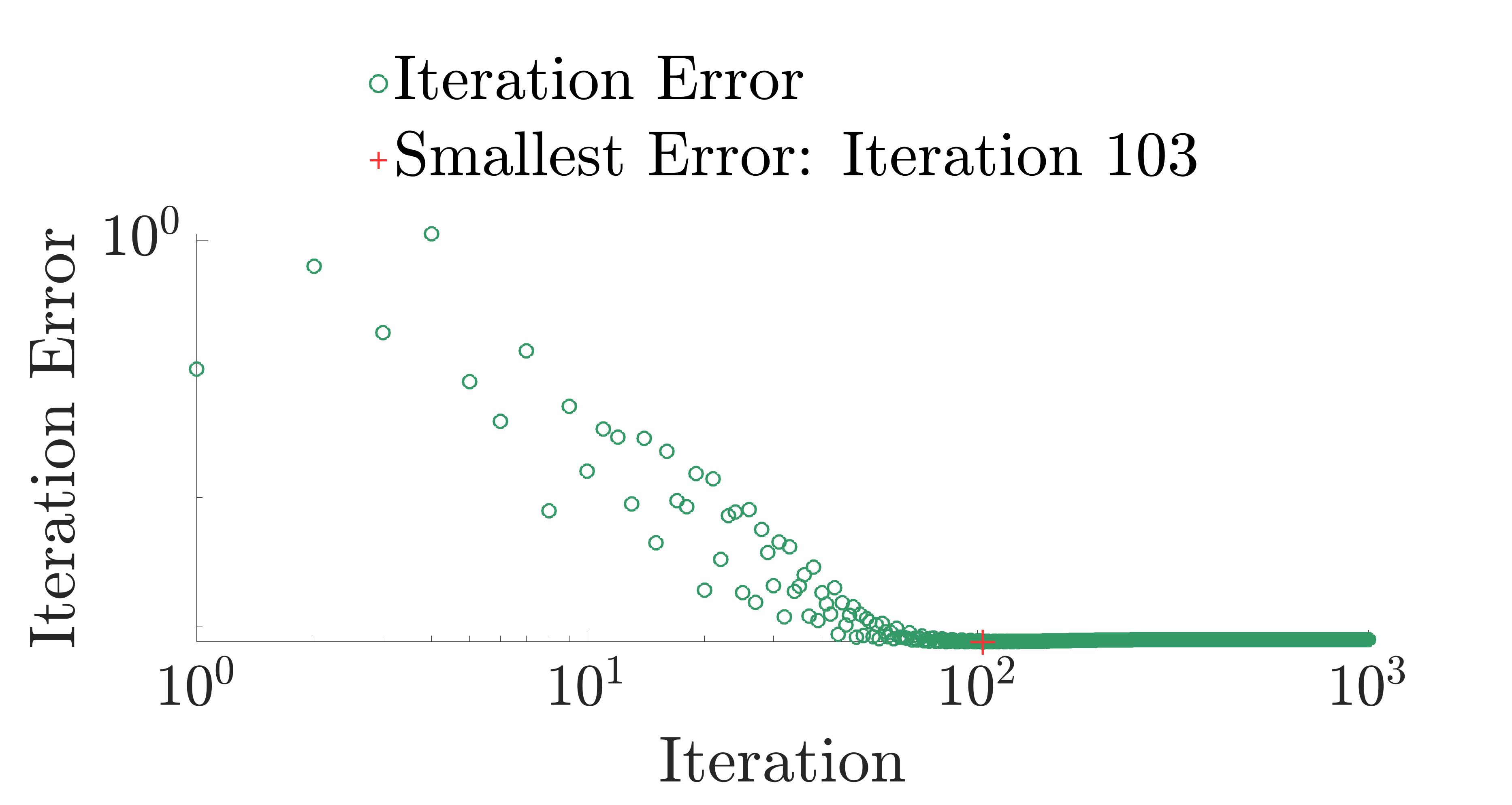}
    \end{subfigure}
    \caption{Results for Case \eqref{case_swe_3b}. Left: plots of the true $p$ and the $p$ corresponding with the smallest residue error, at iteration 103; Middle: plots of the true $p$, the noisy initial guess, and various iteration values for $p$; Right: iteration errors on a log-log scale.}
    \label{fig:swe_3B}
\end{figure}

\begin{figure}[H]
    \begin{subfigure}[t]{0.32\textwidth}
        \centering
        \includegraphics[width=\linewidth]{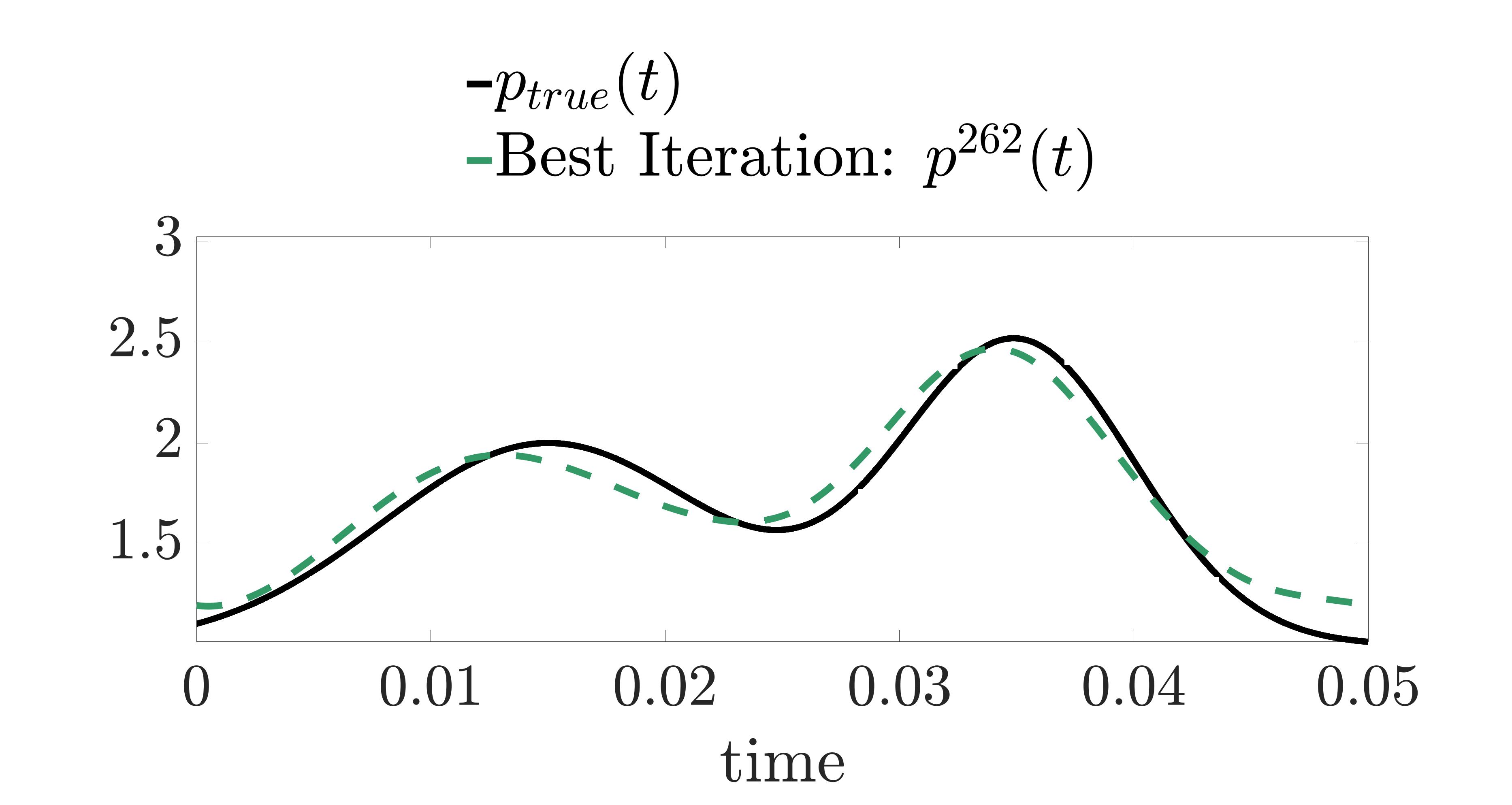}
    \end{subfigure}
    \begin{subfigure}[t]{0.32\textwidth}
        \centering
        \includegraphics[width=\linewidth]{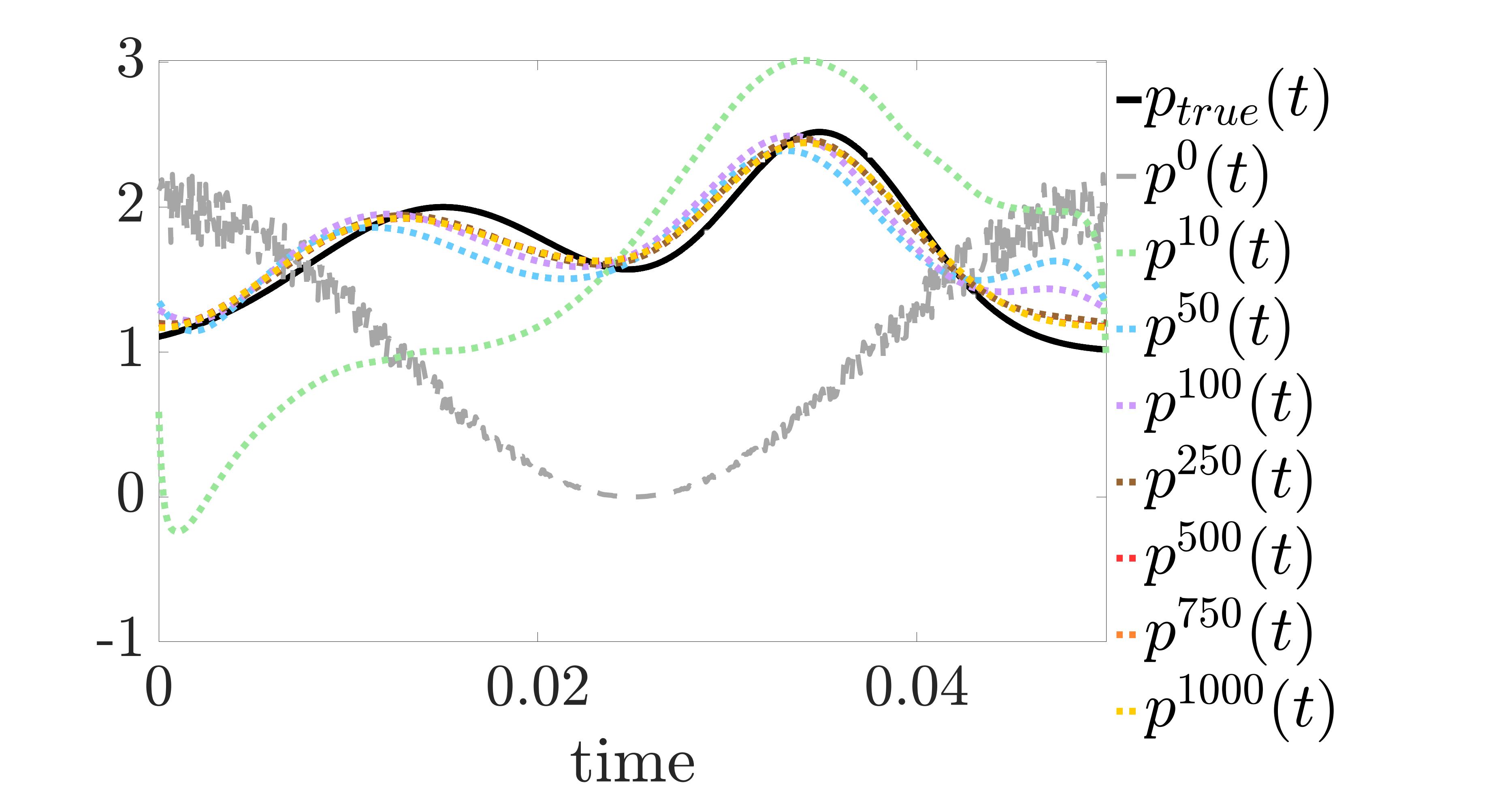}
    \end{subfigure}
    \begin{subfigure}[t]{0.32\textwidth}
        \includegraphics[width=\linewidth]{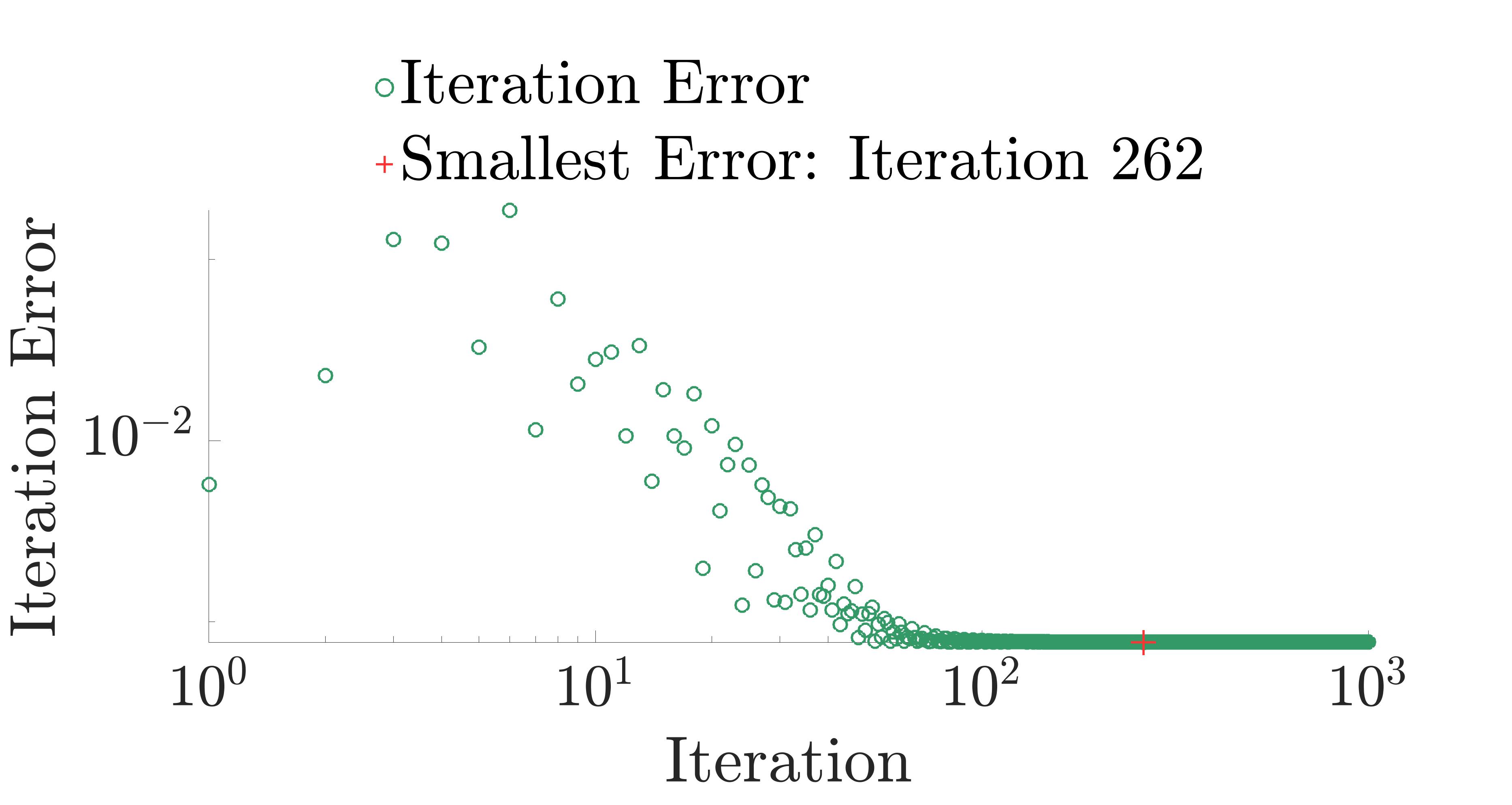}
    \end{subfigure}
    \caption{Results for Case \eqref{case_swe_3c}. Left: plots of the true $p$ and the $p$ corresponding with the smallest residue error, at iteration 262; Middle: plots of the true $p$, the noisy initial guess, and various iteration values for $p$; Right: iteration errors on a log-log scale.}
    \label{fig:swe_3C}
\end{figure}

The numerical results for Cases \eqref{case_swe_3a}, \eqref{case_swe_3b}, and \eqref{case_swe_3c} are shown in Figures \ref{fig:swe_3A} - \ref{fig:swe_3C}, while the results for Case \eqref{case_swe_3d} can be found in Figure \ref{fig:swe_2D}. In all four cases, the scheme was able to identify that $p_{true}(t)$ was a function consisting of two bumps, with the left bump (occurring earlier in time) having a smaller amplitude than the right bump (occurring later in time). This indicates the true function $p(t)$ can be recovered with the initial condition chosen from a wide range of functions. The numerical performances are similar in appearance and convergence rate. Cases \eqref{case_swe_3a} - \eqref{case_swe_3c} all achieve their best guess in less than 300 iterations. The true $p$ and corresponding $p$ with smallest residue error tend to have the some discrepancy near the final time $T$. The exception is Case \eqref{case_swe_3a} in which the initial guess for $p$ at time $T$ is near to $p_{true}(T)$.

\subsection{Convergence and Accuracy}
\hly{We discuss the convergence and accuracy of our numerical schemes in this section.} All numerical simulations presented in this section consist of the same problem set up as was used in Case \eqref{case_swe_2d}.

First, we discuss the local convergence of the optimization scheme. Since our problem is ill-posed, we structure our algorithm such that it optimizes the sum of the residual errors and the regularizers, that is $\J = \J_0(p) + \R_{L^1}(p) + \R_{H^1}(p)$. We apply our method to $\J$, which is locally-convex, and therefore the algorithm converges locally in first order to a local minimum of $J$. 

\begin{center}
\begin{figure}[H]
    \begin{subfigure}[t]{0.49\textwidth}
        \centering
        \includegraphics[width=\textwidth]{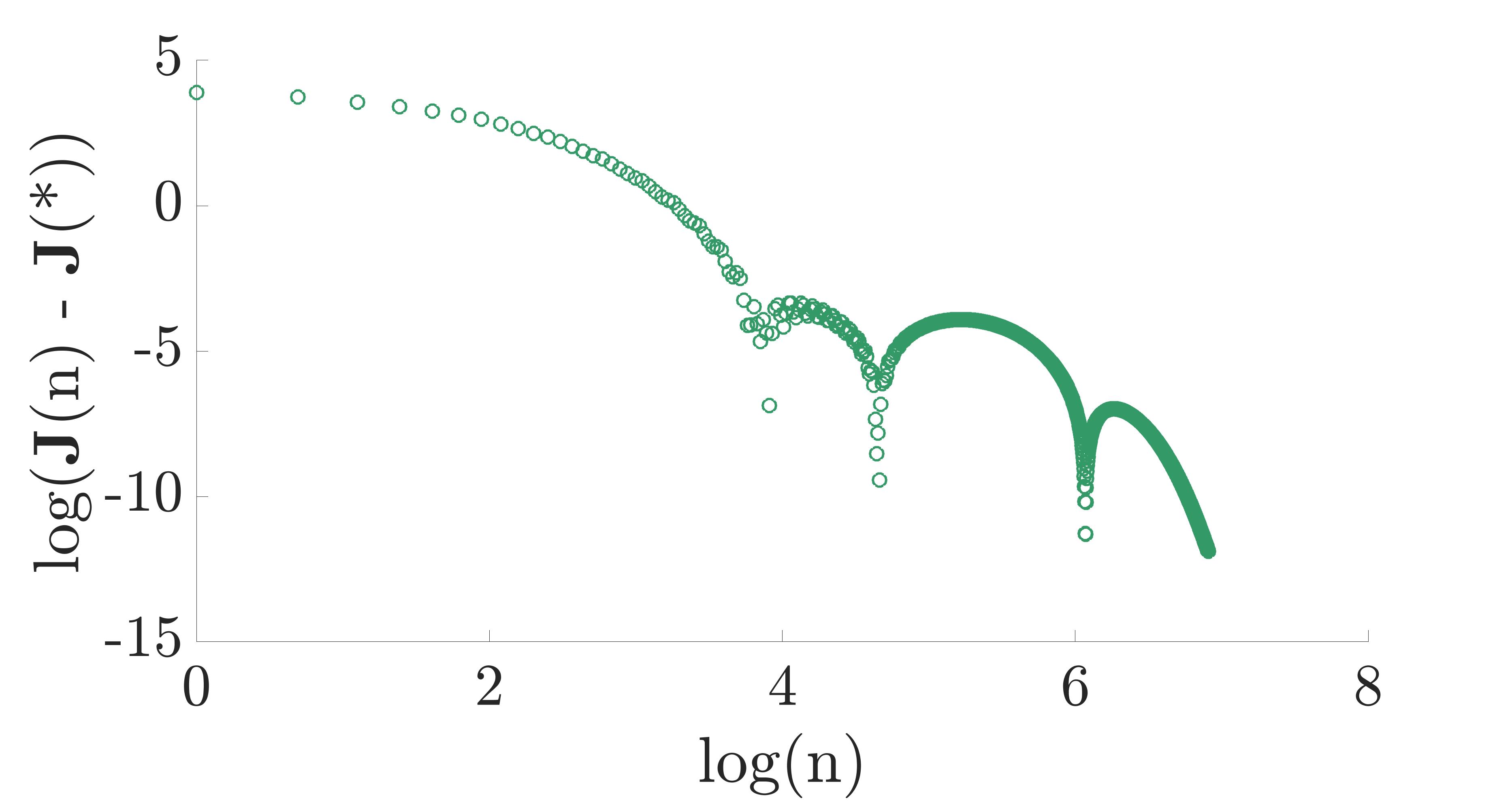}
        \caption{Global convergence behavior.}
    \end{subfigure}
    \begin{subfigure}[t]{0.49\textwidth}
        \centering
        \includegraphics[width=\textwidth]{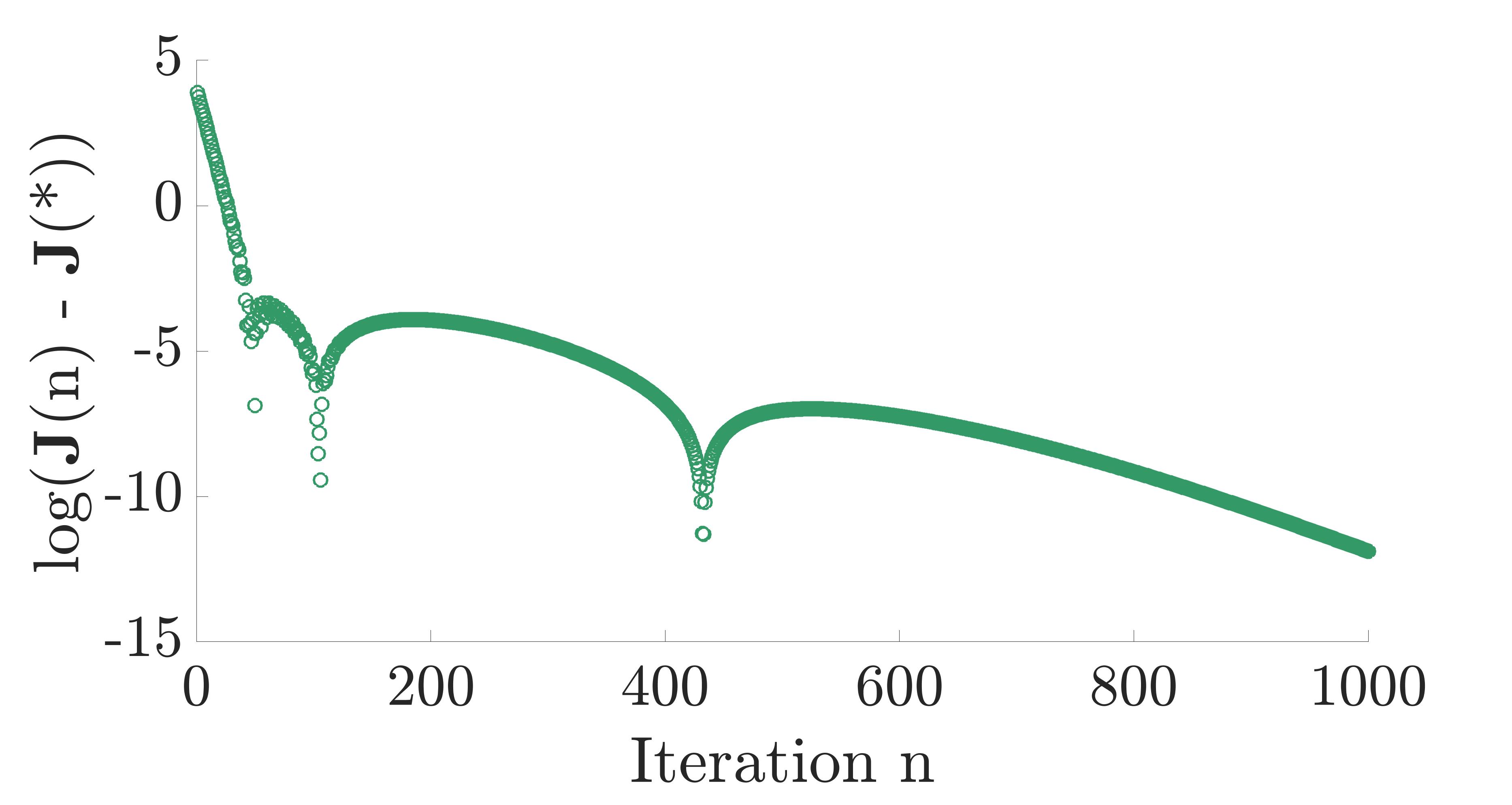}
        \caption{Global convergence behavior.}
    \end{subfigure}\\
    \begin{subfigure}[t]{0.49\textwidth}
        \includegraphics[width=\textwidth]{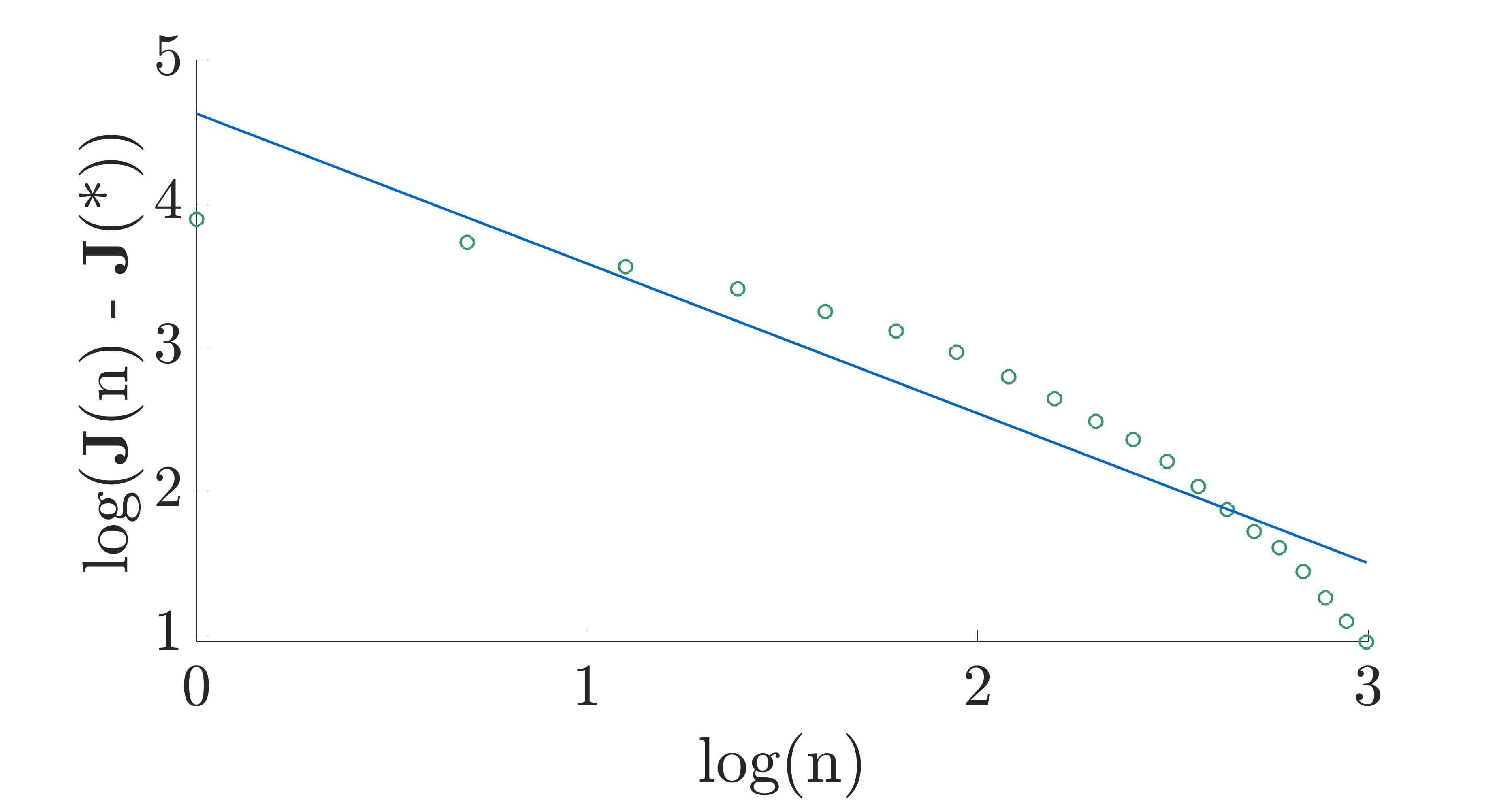}
        \caption{Front end convergence behavior.}
    \end{subfigure}
    \begin{subfigure}[t]{0.49\textwidth}
        \includegraphics[width=\textwidth]{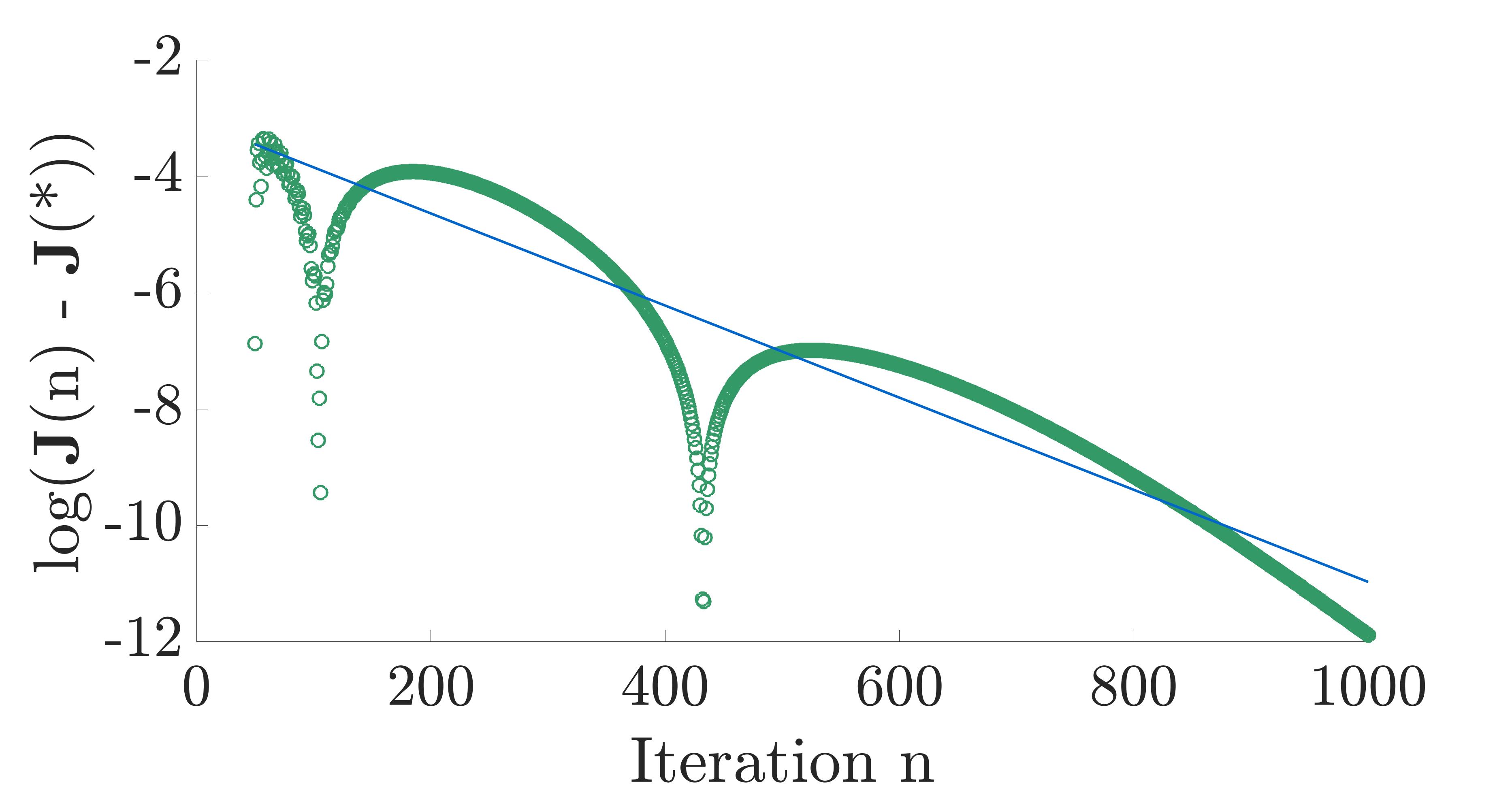}
        \caption{Tail end convergence behavior.}
    \end{subfigure}
    \caption{\hly{Plots for convergence behavior} corresponding to regularization parameters $\gamma_L = 1 \times 10^{-6}$ and $\gamma_H = 5 \times 10^{-8}$. Here $\J(n)$ is the error at iteration $n$ including regularization and $\J(*)$ is the error at the local minimum. Plots in the top row demonstrate the global convergence behavior with (a) a log scale on the x-axis and (b) a standard scale on the x-axis. The bottom row of plots show local convergence behavior overlaid with the best fit linear line. Specifically, (c) demonstrates the first order convergence for the front end iterations (1 to 20) on a log-log scale with a fit line slope of -1.0417 and (d) demonstrates the first order convergence for the tail iterations (50 to 1000) with a log scale on the y-axis and a fit line slope of -0.0079.}
    \label{fig:err_convg_1}
\end{figure}
\end{center}

\begin{figure}[H]
    \begin{subfigure}[t]{0.49\textwidth}
        \centering
        \includegraphics[width=\textwidth]{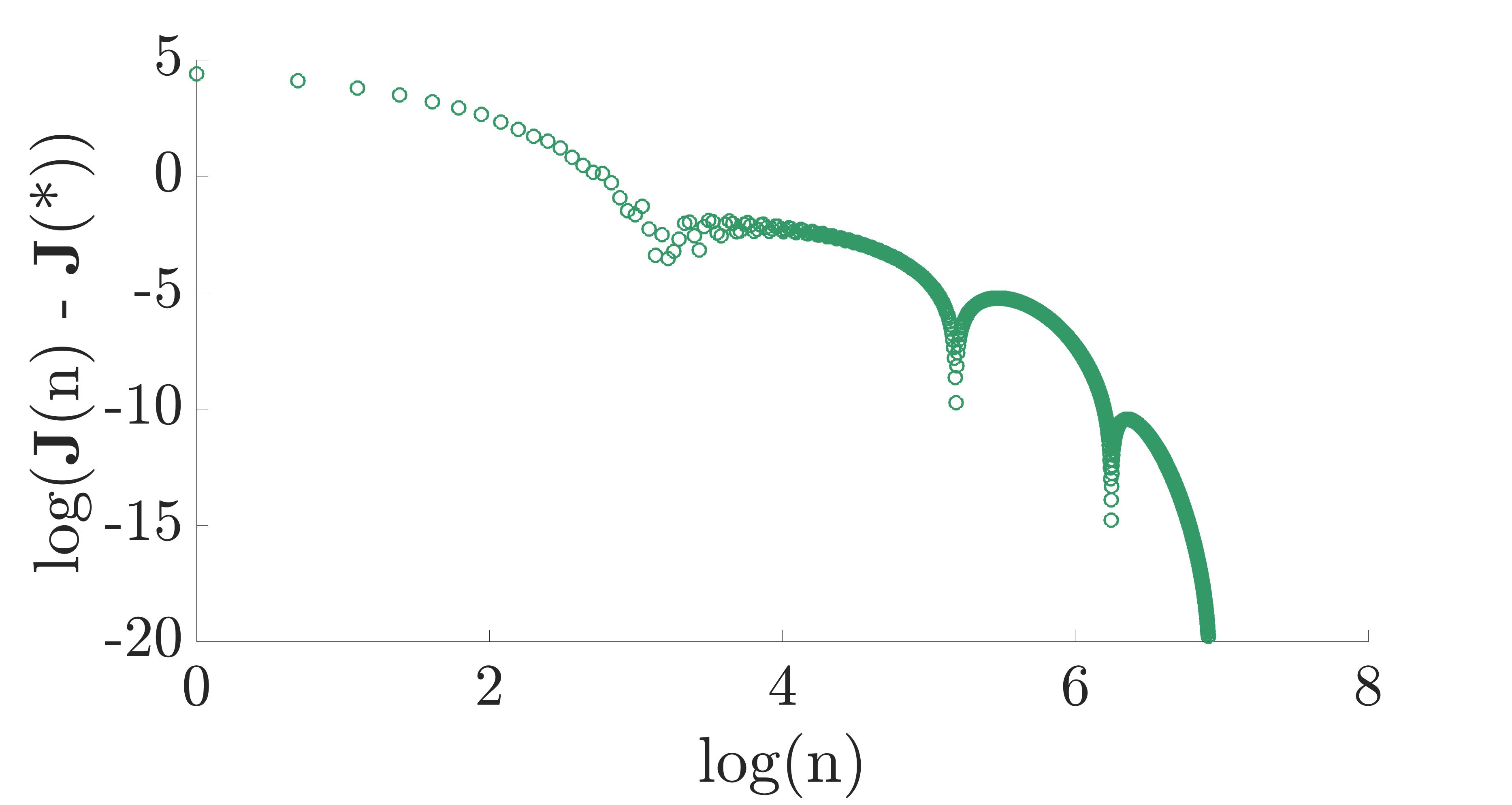}
        \caption{Global convergence behavior.}
    \end{subfigure}
    \begin{subfigure}[t]{0.49\textwidth}
        \centering
        \includegraphics[width=\textwidth]{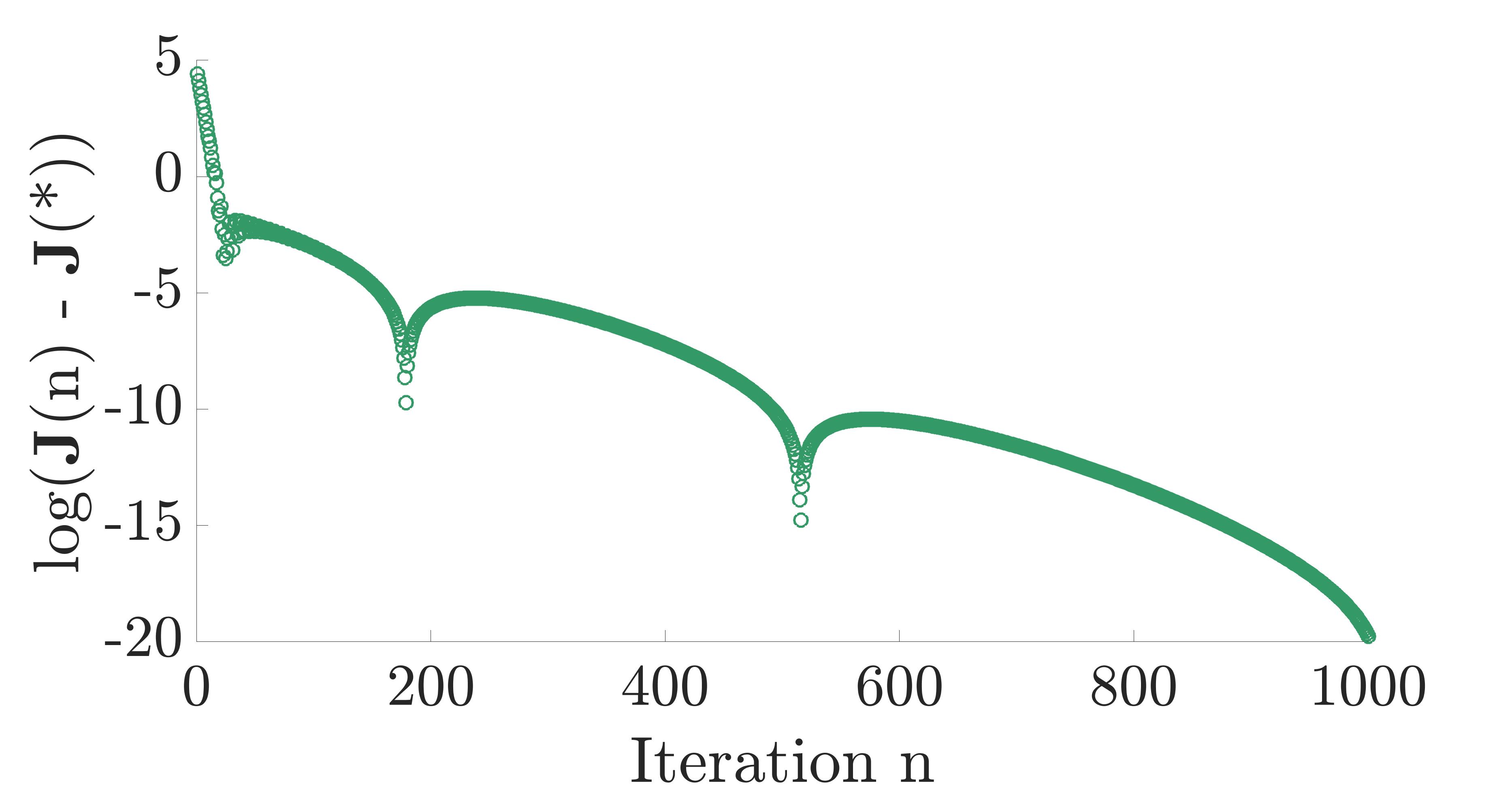}
        \caption{Global convergence behavior.}
    \end{subfigure}
    \begin{subfigure}[t]{0.49\textwidth}
        \includegraphics[width=\textwidth]{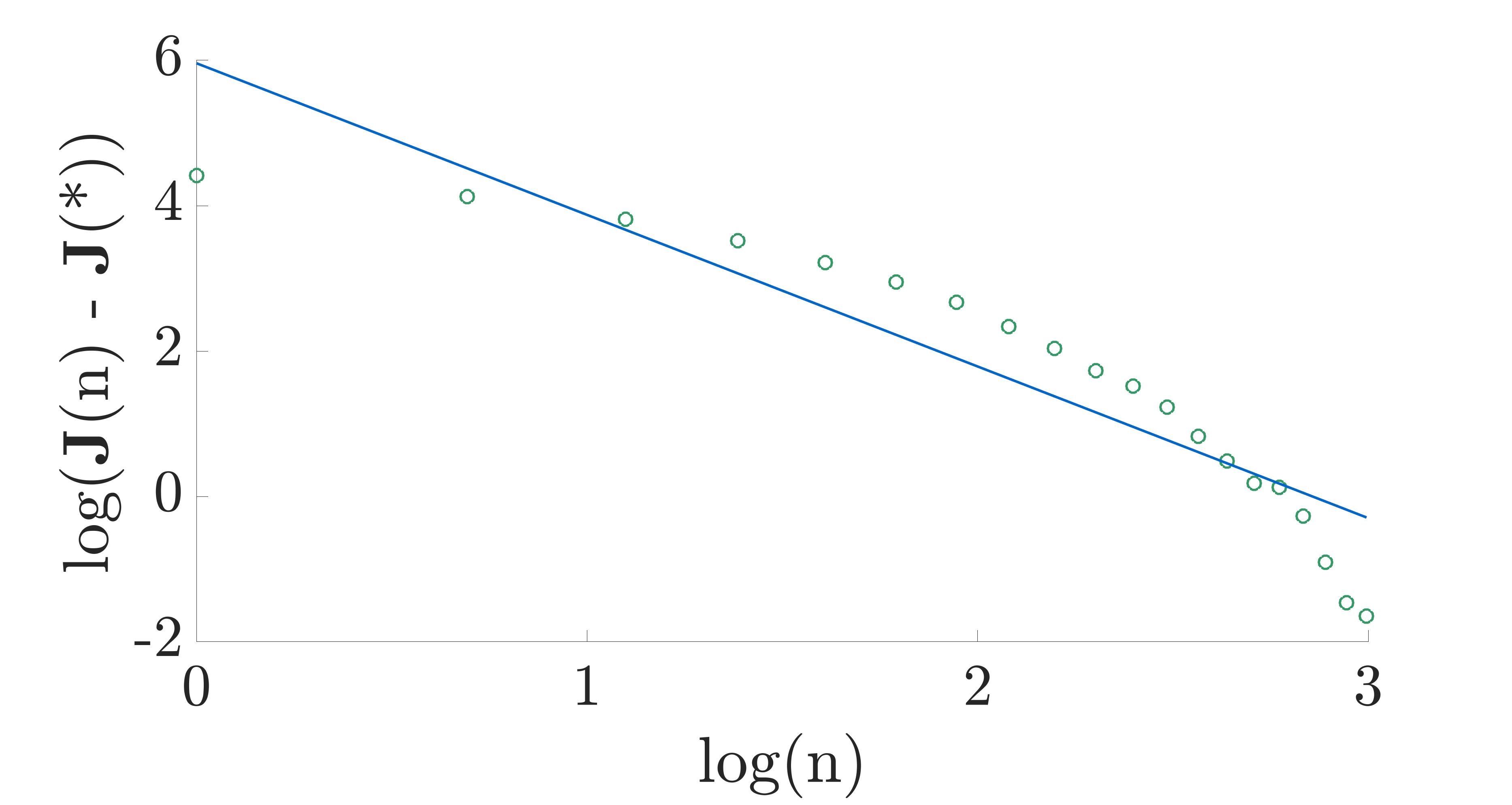}
        \caption{Front end convergence behavior.}
    \end{subfigure}
    \begin{subfigure}[t]{0.49\textwidth}
        \includegraphics[width=\textwidth]{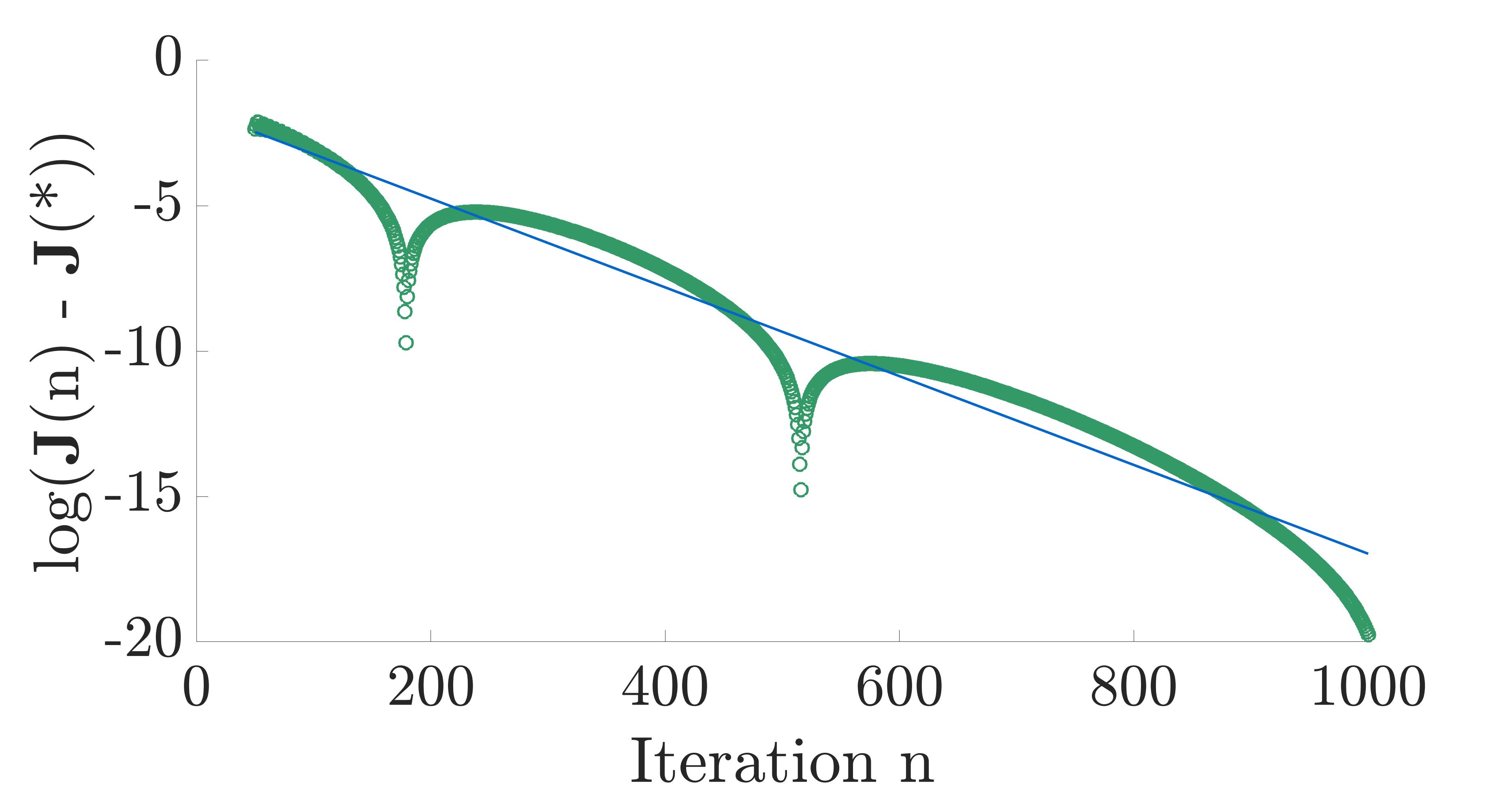}
        \caption{Tail end convergence behavior.}
    \end{subfigure}
    \caption{\hly{Plots for convergence behavior} corresponding to regularization parameters $\gamma_L = 5 \times 10^{-6}$ and $\gamma_H = 1 \times 10^{-7}$. Here $\J(n)$ is the error at iteration $n$ including regularization and $\J(*)$ is the error at the local minimum. Plots in the top row demonstrate the global convergence behavior with (a) a log scale on the x-axis and (b) a standard scale on the x-axis. The bottom row of plots show local convergence behavior overlaid with the best fit linear line. Specifically, (c) demonstrates the first order convergence for the front end iterations (1 to 20) on a log-log scale with a fit line slope of -2.0861 and (d) demonstrates the first order convergence for the tail iterations (50 to 1000) with a log scale on the y-axis and a fit line slope of -0.0153.}
    \label{fig:err_convg_2}
\end{figure}

For illustrative purpose, Figures \ref{fig:err_convg_1} and \ref{fig:err_convg_2} demonstrate the convergence behavior for two different choices of regularization. The results in Figure \ref{fig:err_convg_1} correspond to regularization parameters $\gamma_L = 1 \times 10^{-6}$ and $\gamma_H = 5 \times 10^{-8}$ while the results in Figure \ref{fig:err_convg_2} correspond to slightly higher regularization parameters $\gamma_L = 5 \times 10^{-6}$ and $\gamma_H = 1 \times 10^{-7}$. In both cases we examine the global and local convergence behaviors. The difference between the iteration error and the local minimum error on a log scale is plotted against the log of the iteration number in Figures \ref{fig:err_convg_1}(a) and \ref{fig:err_convg_2}(a) and against the iteration number in Figures \ref{fig:err_convg_1}(b) and \ref{fig:err_convg_2}(b). These plots show different convergence behavior at the front end and the tail end of the iterative process. Hence, local convergence behavior is shown as well. For both cases we plot the errors against the first 20 iterations on a log-log scale with a best fit line to show the linear convergence. The results in Figure \ref{fig:err_convg_1}(c) is matched to a best fit line with a slope of $-1.0417$ showing first order convergence. The best fit line in Figure \ref{fig:err_convg_2}(c) has a slope of $-2.0861$ showing that an increased rate of convergence occurs with increased regulation rates. On the other hand, we also plot the errors corresponding to iterations 50 to 1000 on a log scale against the iterations with a best fit line. Since the x-axis, in this case, is not on a log scale, the slope of the best fit line now corresponds to exponential rate of convergence. The results in Figure \ref{fig:err_convg_1}(d) correspond with a best fit line consisting of a slope of $-0.0079$ while the results for the case with higher regularization parameters demonstrated in Figure \ref{fig:err_convg_2}(d) have a slope of $-0.0153$, which is again twice as large as the case with smaller regularization parameters. 

Second, we discuss the accuracy of the scheme. The solution of the forward problem is not analytically available, hence the order of accuracy is determined by comparing the state variables computed on uniform mesh of cell width $\Delta x$ that is repeatedly halved. We again consider the problem set up that was used in Case \eqref{case_swe_2d} with regularization parameters $\gamma_L = 1 \times 10^{-6}$ and $\gamma_H = 5 \times 10^{-8}$ and no noise. The $L^1$ errors and corresponding orders of accuracy for the forward problem computed with $P^0, P^1,$ and $P^2$ piecewise polynomials are found in Table \ref{tab:acc}. We see that an order of $k+1$ accuracy is achieved for each polynomial degree $k$ for both state variables $h$ and $hu$.

\begin{center}
\begin{small}
    \begin{table}[H]
        \centering
        \begin{tabular}{c c c c c c c c }
        \toprule
            & & \multicolumn{2}{c}{$k=0$} & \multicolumn{2}{c}{$k=1$}&\multicolumn{2}{c}{$k=2$} \\
        \midrule
              & N & $L^1$ Error & Order & $L^1$ Error & Order & $L^1$ Error & Order \\
        \midrule
              & 25  & 0.0490 & -      & $2.9640 \times 10^{-3}$ & -      & $2.4032 \times 10^{-3}$ & - \\  
              & 50  & 0.0260 & 0.9143 & $4.9877 \times 10^{-4}$ & 2.5711 & $1.1128 \times 10^{-5}$ & 7.7547 \\ 
            $h$ & 100 & 0.0133 & 0.9725 & $1.2434 \times 10^{-4}$ & 2.0041 & $1.4679 \times 10^{-6}$ & 2.9223 \\ 
              & 200 & 0.0067 & 0.9858 & $3.0999 \times 10^{-5}$ & 2.0040 & $1.7730 \times 10^{-7}$ & 3.0494 \\ 
              & 400 & 0.0034 & 0.9927 & $7.7396 \times 10^{-6}$ & 2.0019 & $2.2152 \times 10^{-8}$ & 3.0008\\ 
        \midrule
              & 25  & 0.4101 & -      & $5.3311 \times 10^{-2}$ & -      & $5.1219 \times 10^{-2}$ & - \\
              & 50  & 0.1980 & 1.0504 & $5.1953 \times 10^{-3}$ & 3.3592 & $8.5791 \times 10^{-5}$ & 9.2216\\
            $hu$ & 100 & 0.1008 & 0.9745 & $1.4146 \times 10^{-3}$ & 1.8768 & $1.1960 \times 10^{-5}$ & 2.8425\\
              & 200 & 0.0509 & 0.9857 & $3.7772 \times 10^{-4}$ & 1.9050 & $1.7950 \times 10^{-6}$ & 2.7362\\
              & 400 & 0.0256 & 0.9928 & $9.6943 \times 10^{-5}$ & 1.9621 & $2.3099 \times 10^{-7}$ & 2.9581\\
        \bottomrule
        \end{tabular}        
        \caption{\hly{$L^1$ errors and convergence orders of the accuracy test for the forward problem}  using $P^0$, $P^1$ and $P^2$ piecewise polynomials and a uniform mesh of N cells. In each case, $k+1$ order of accuracy is achieved. 
        }
        \label{tab:acc}
    \end{table}
\end{small}
\end{center}

\subsection{Impact of Inconsistent Discretization Schemes}
\begin{figure}[h]
    \begin{subfigure}[t]{0.32\textwidth}
        \centering
        \includegraphics[width=\linewidth]{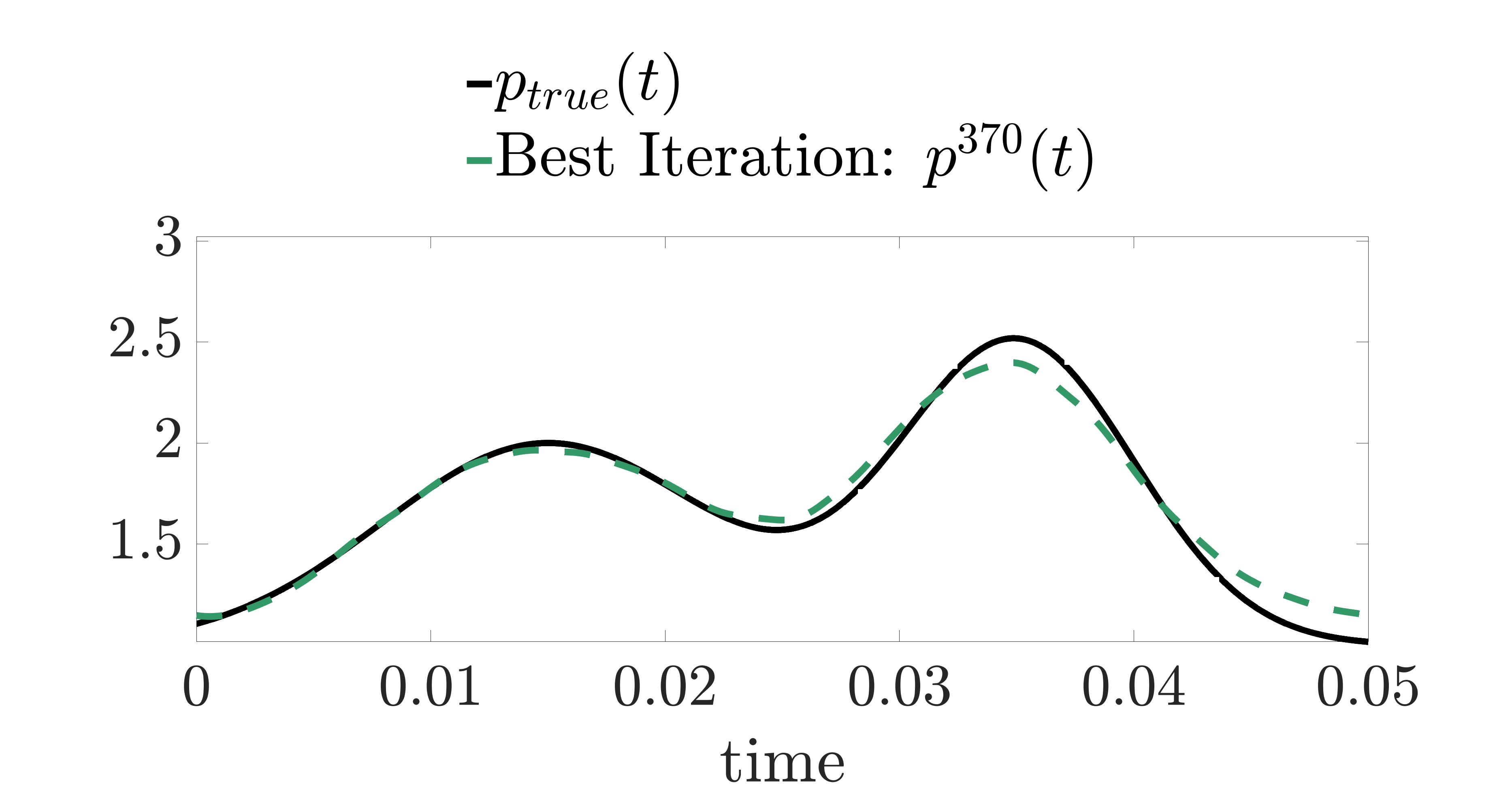}
    \end{subfigure}
    \begin{subfigure}[t]{0.32\textwidth}
        \centering
        \includegraphics[width=\linewidth]{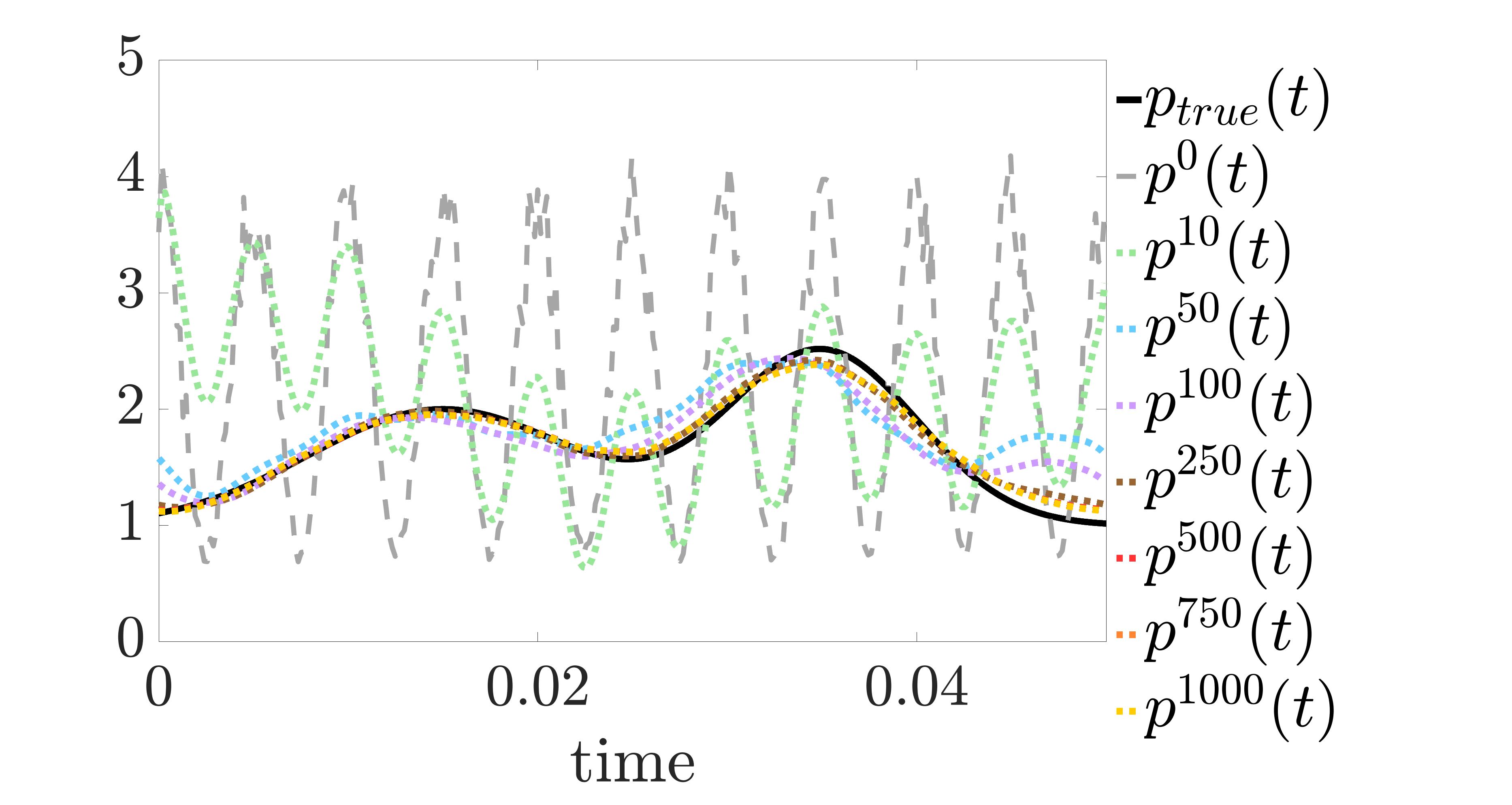}
    \end{subfigure}
    \begin{subfigure}[t]{0.32\textwidth}
        \includegraphics[width=\linewidth]{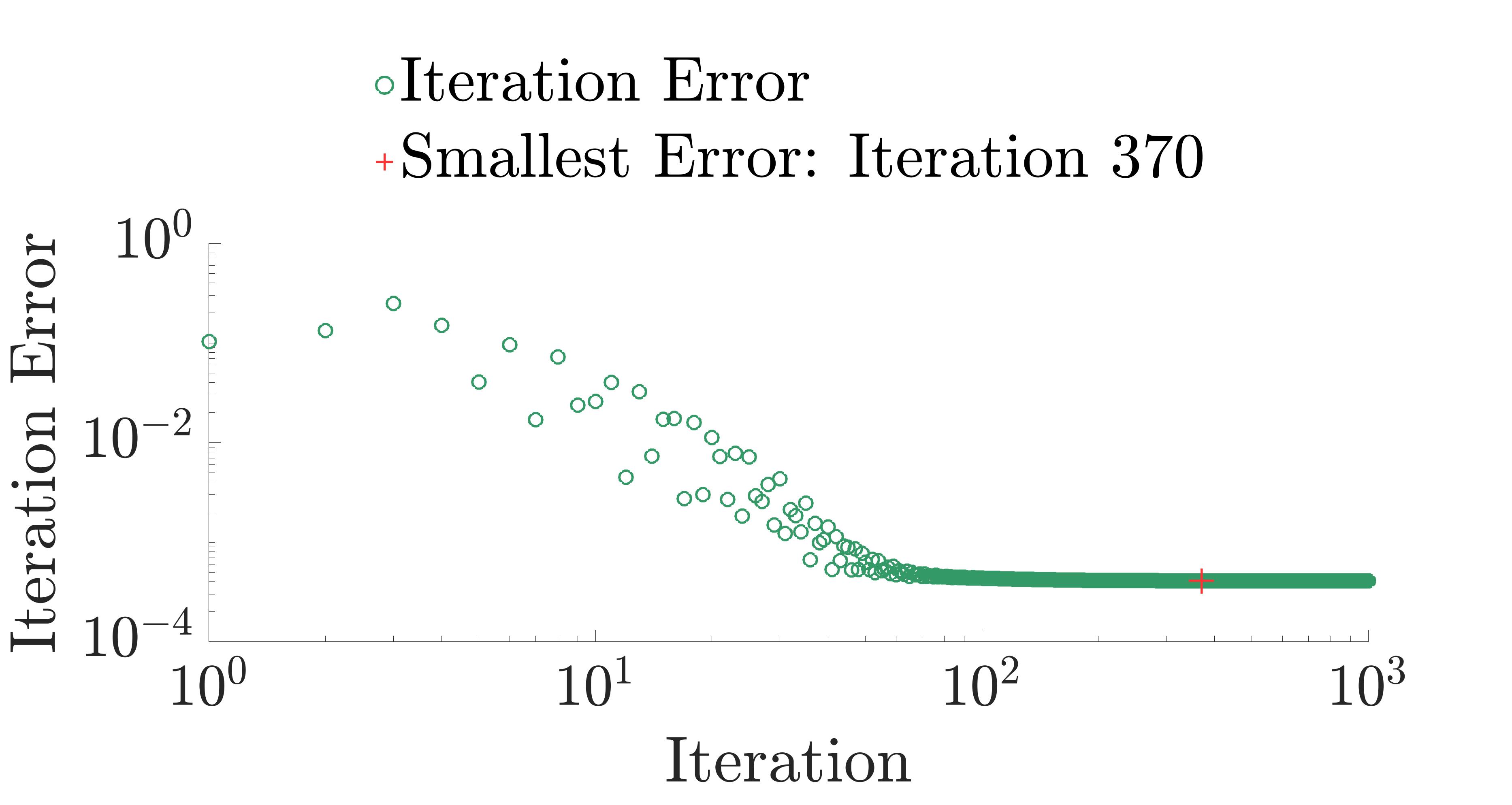}
    \end{subfigure}
    \caption{Results for Case \eqref{case_swe_2d} under consistent conditions. Left: plots of the true $p$ and the $p$ corresponding with the smallest residue error, at iteration 370; Middle: plots of the true $p$, the noisy initial guess, and various iteration values for $p$; Right: iteration errors on a log-log scale.}
    \label{fig:swe_6A}
\end{figure}

\hly{In this section we further justify the decision to employ different discretization schemes for the forward and (linearized) adjoint schemes}. The primary reason for this choice is to avoid any speculation of inverse crime. Here we compare the numerical results for case \eqref{case_swe_2d} under two conditions. First, Figure \ref{fig:swe_6A} shows the numerical results under consistent conditions in which the forward and adjoint problems are both discretized using $P^2$ polynomials and a uniform mesh of 25 cells. Second, Figure \ref{fig:swe_2Db} displays the numerical results under inconsistent conditions. The forward problem was discretized with $P^2$ piecewise polynomials and a mesh of $50$ uniform cells, while the adjoint problem was solved with $P^1$ piecewise polynomials and a uniform mesh of $25$ cells. The inconsistency of the discretized gradient, derived from the discrete adjoint solution, with the discretized forward model does not contaminate the behavior of the inverse problem algorithm as evidence by the results in Figure \ref{fig:swe_2Db}. 

\subsection{Impact of Regularization Terms}
\begin{figure}[h]
    \begin{subfigure}[t]{0.32\textwidth}
        \centering
        \includegraphics[width=\linewidth]{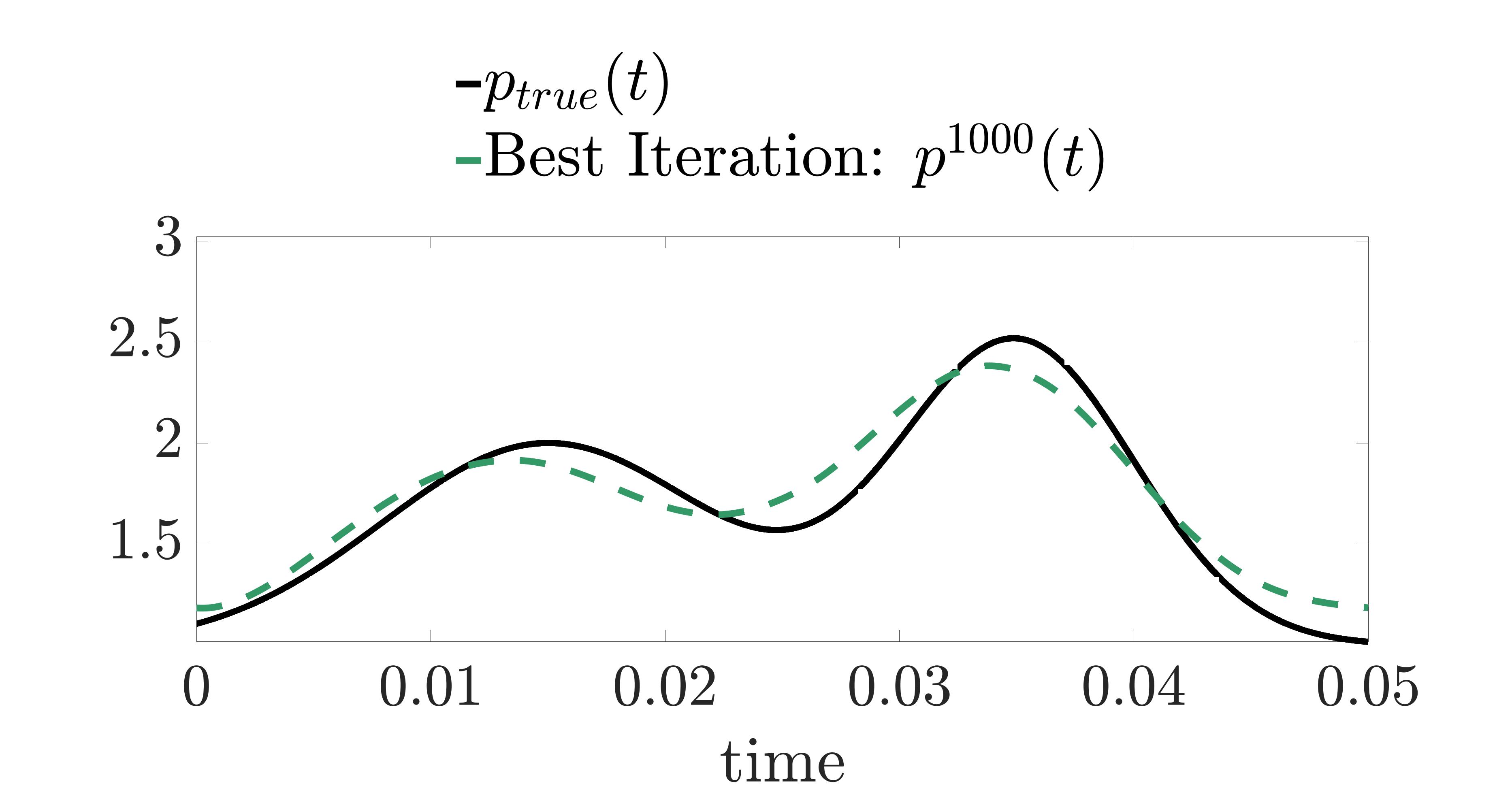}
    \end{subfigure}
    \begin{subfigure}[t]{0.32\textwidth}
        \centering
        \includegraphics[width=\linewidth]{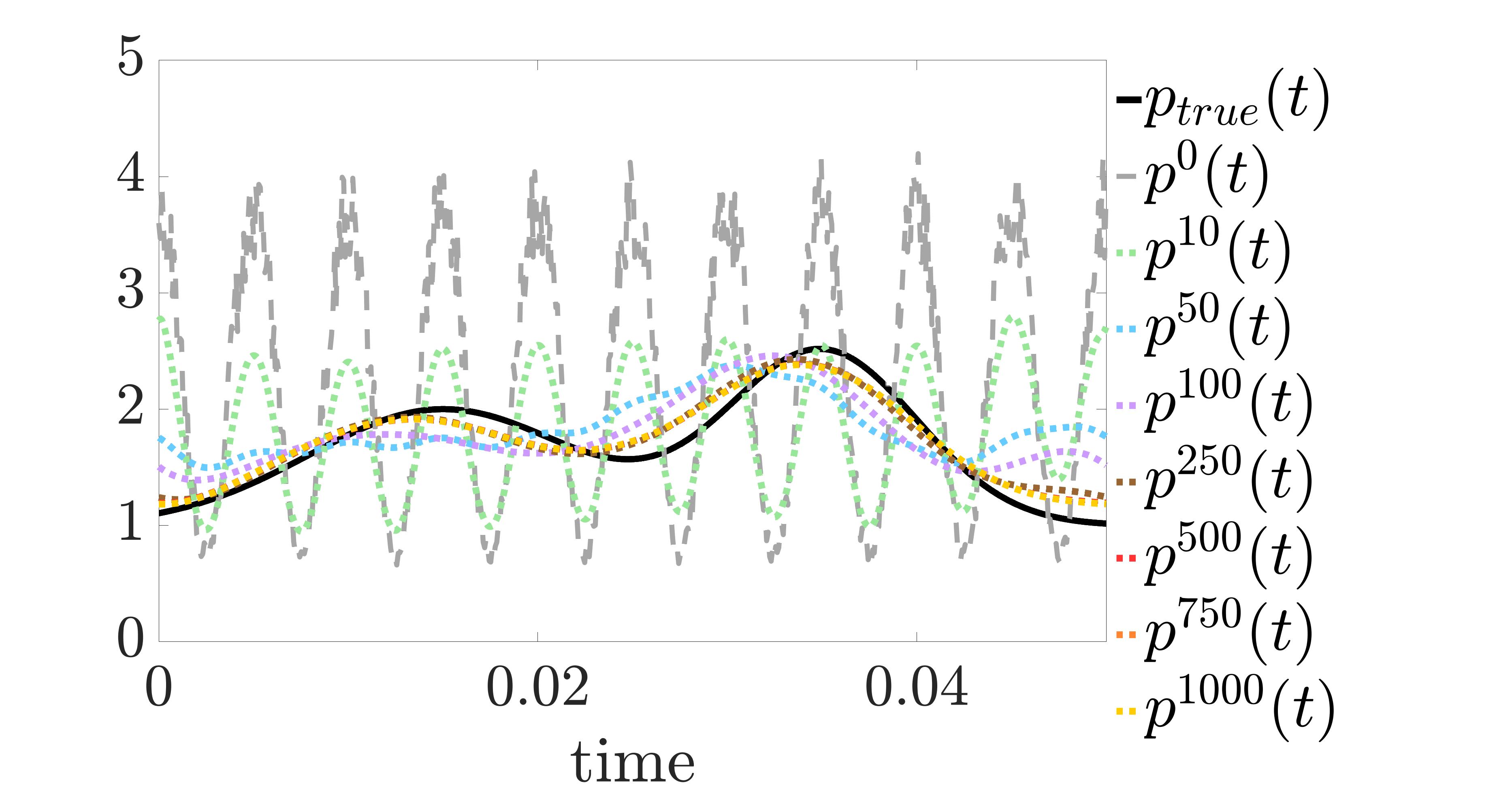}
    \end{subfigure}
    \begin{subfigure}[t]{0.32\textwidth}
        \includegraphics[width=\linewidth]{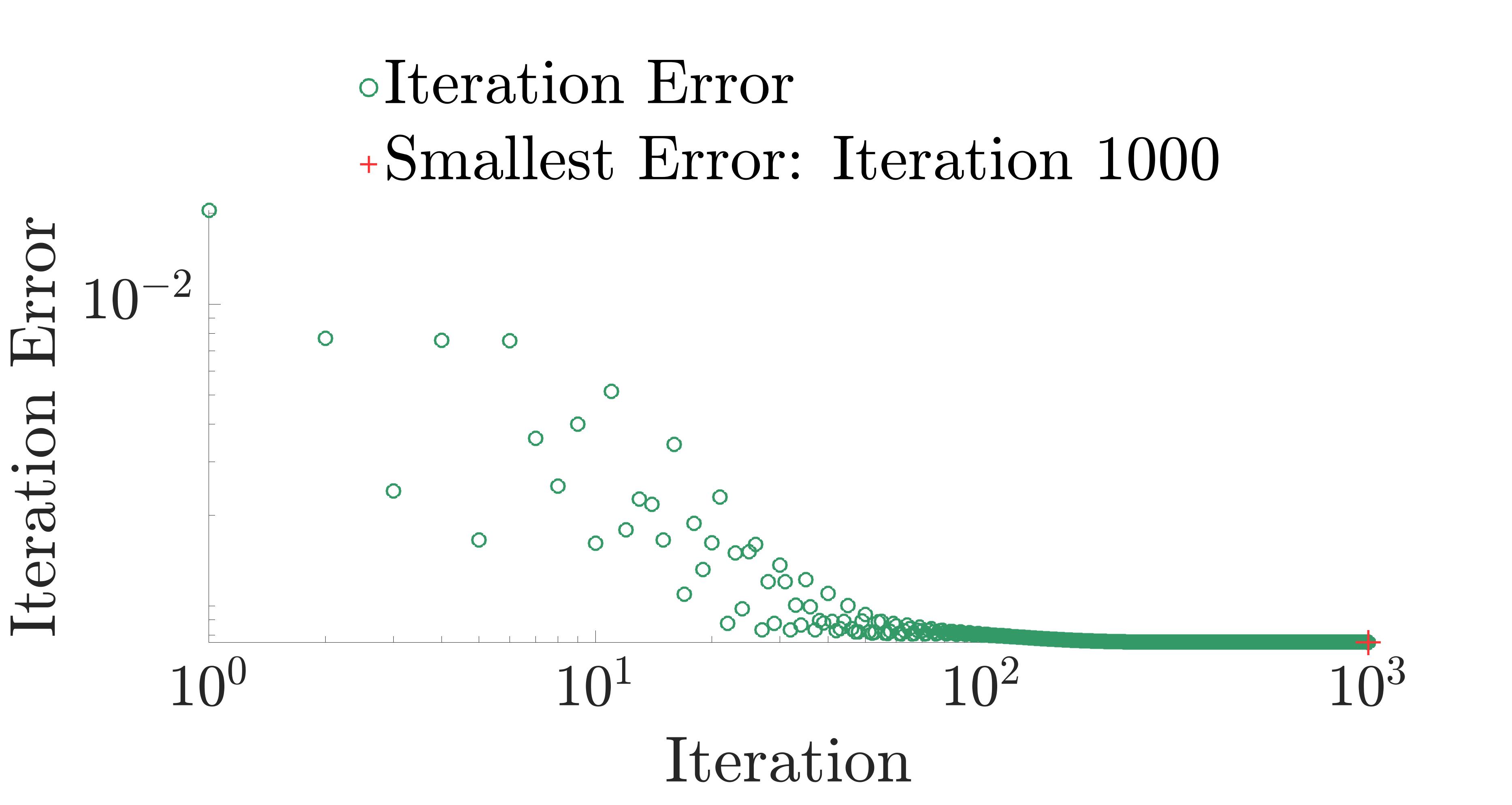}
    \end{subfigure}
    \caption{Results for Case \eqref{case_swe_5a}. 
    Left: plots of the true $p$ and the numerically recovered $p$ at iteration 1000 corresponding with the smallest residue; Middle: plots of the true $p$, the noisy initial guess, and various iteration values for $p$; Right: iteration errors on a log-log scale.}
    \label{fig:swe_5A}
\end{figure}

\begin{figure}[h]
    \begin{subfigure}[t]{0.32\textwidth}
        \centering
        \includegraphics[width=\linewidth]{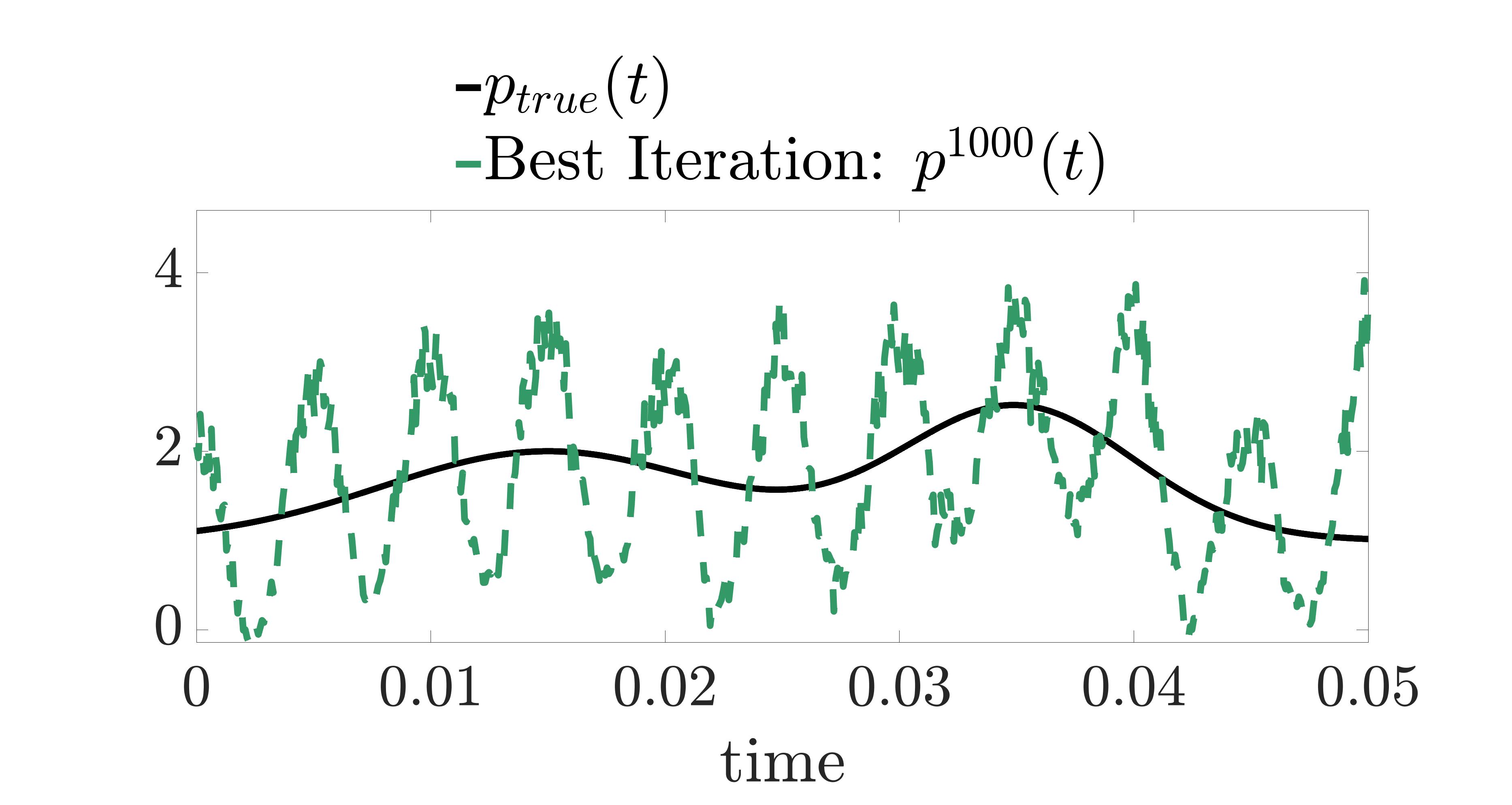}
    \end{subfigure}
    \begin{subfigure}[t]{0.32\textwidth}
        \centering
        \includegraphics[width=\linewidth]{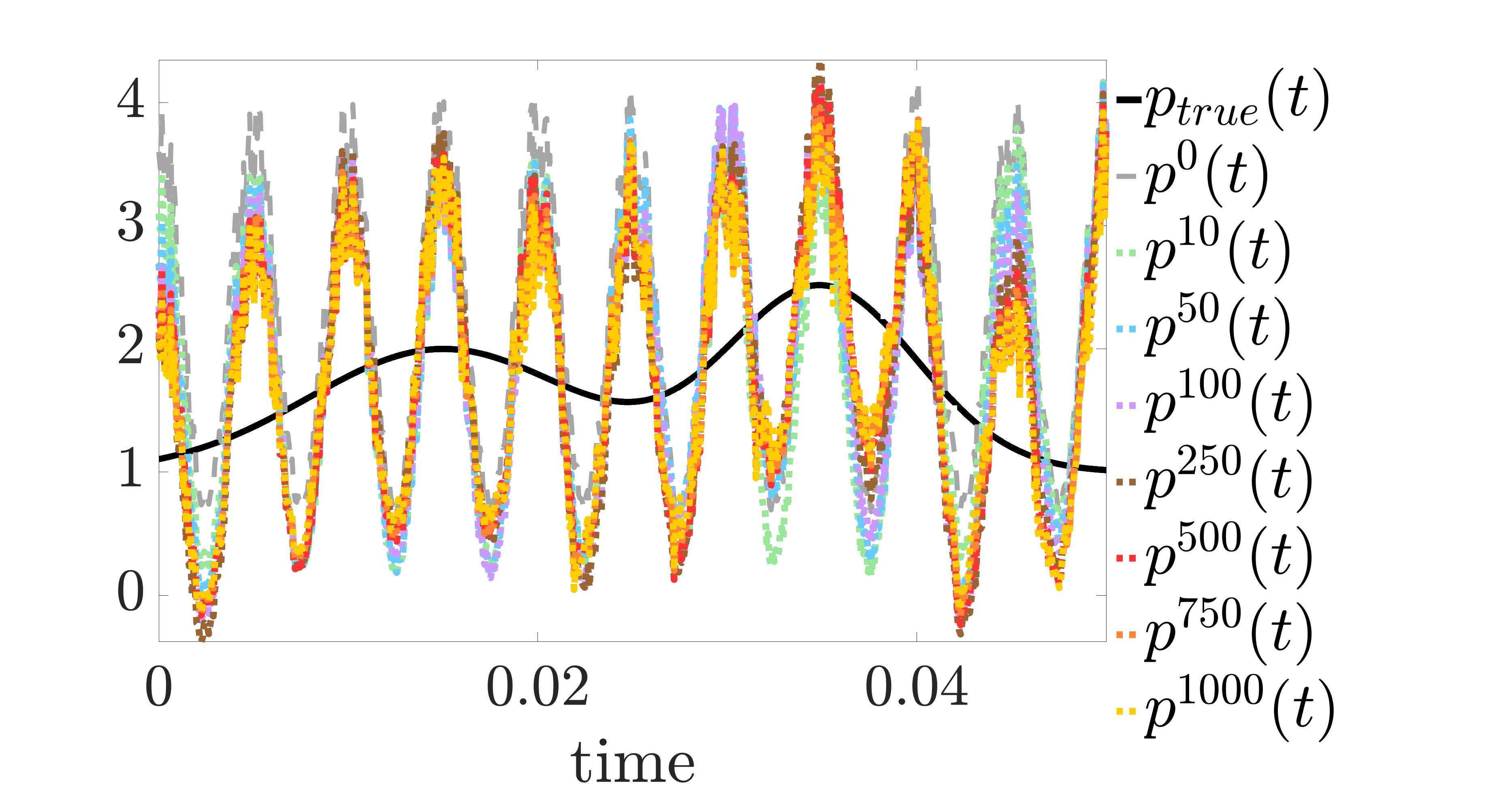}
    \end{subfigure}
    \begin{subfigure}[t]{0.32\textwidth}
        \includegraphics[width=\linewidth]{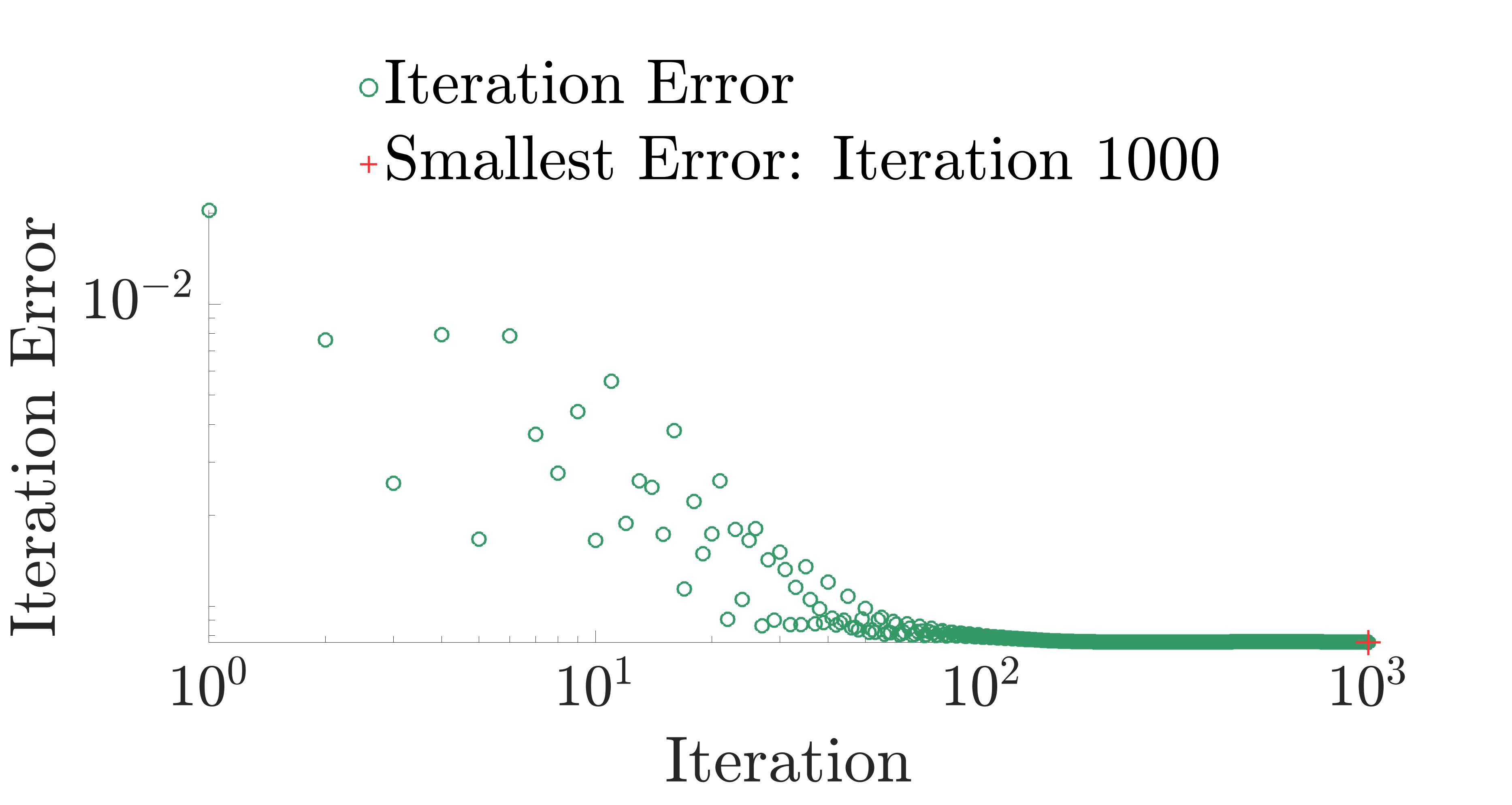}
    \end{subfigure}
    \caption{Results for Case \eqref{case_swe_5b}. 
    Left: plots of the true $p$ and the numerically recovered $p$ at iteration 1000 corresponding with the smallest residue; Middle: plots of the true $p$, the noisy initial guess, and various iteration values for $p$; Right: iteration errors on a log-log scale.}
    \label{fig:swe_5B}
\end{figure}

\begin{figure}[h]
    \begin{subfigure}[t]{0.32\textwidth}
        \centering
        \includegraphics[width=\linewidth]{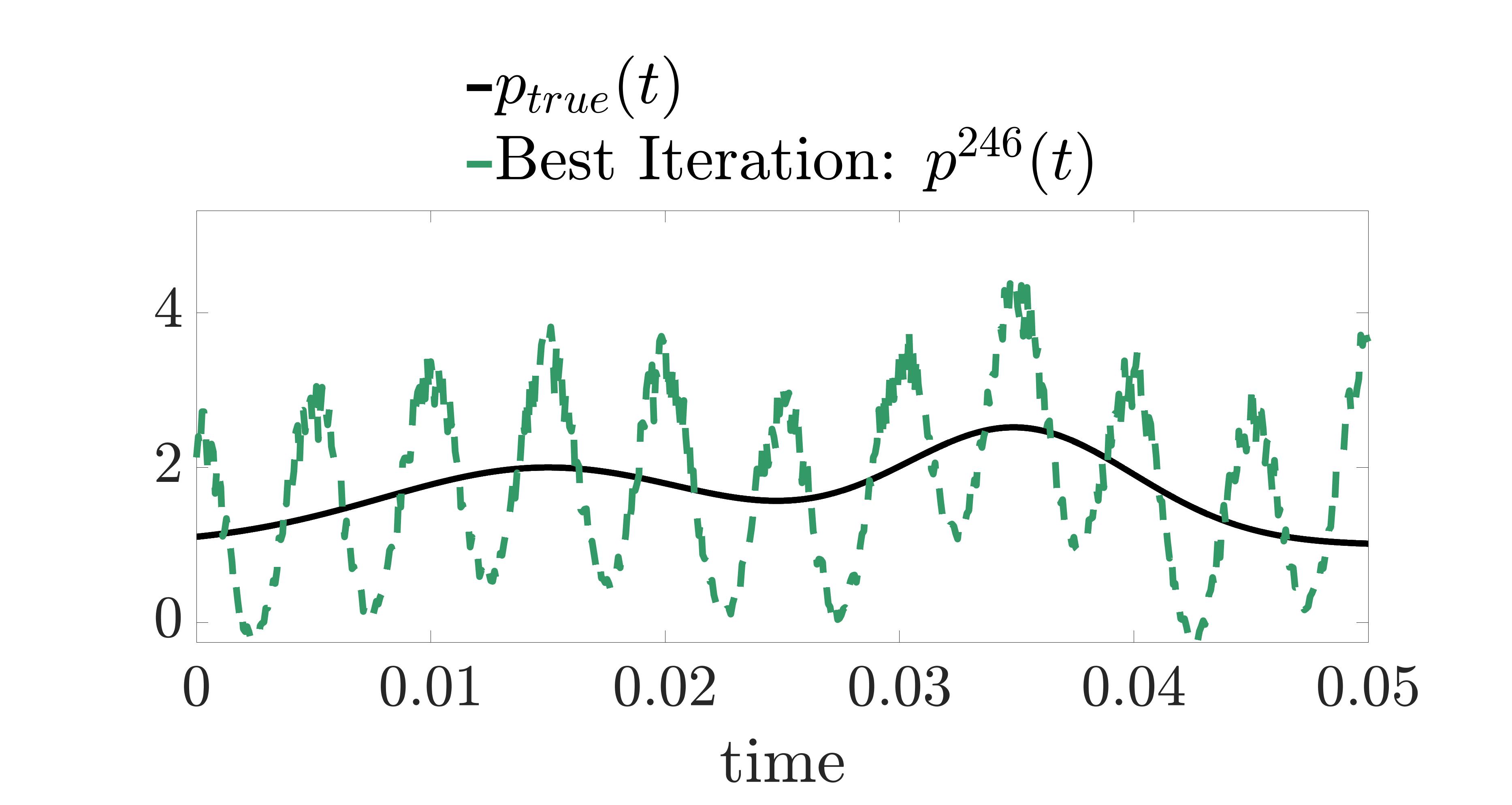}
    \end{subfigure}
    \begin{subfigure}[t]{0.32\textwidth}
        \centering
        \includegraphics[width=\linewidth]{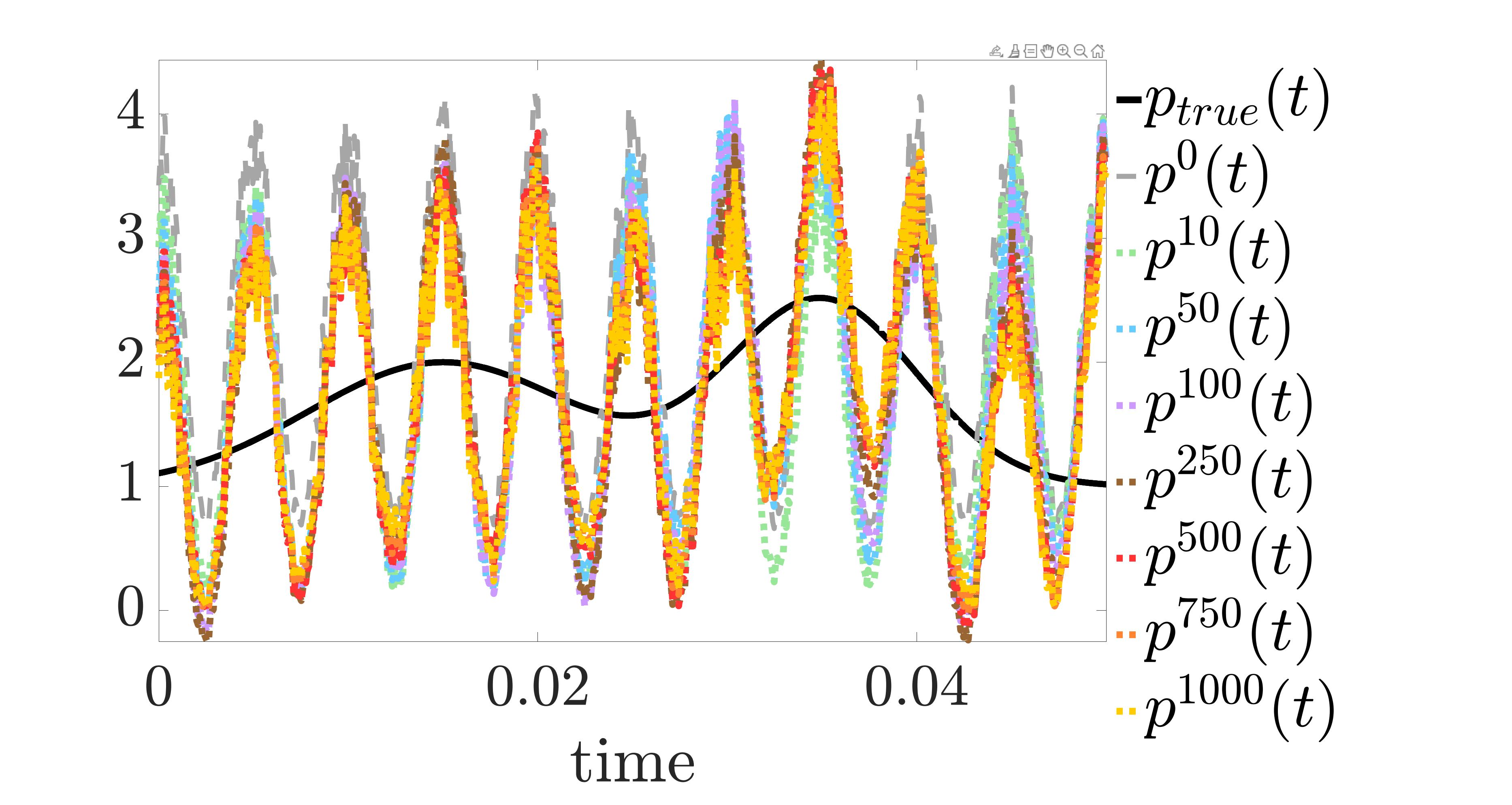}
    \end{subfigure}
    \begin{subfigure}[t]{0.32\textwidth}
        \includegraphics[width=\linewidth]{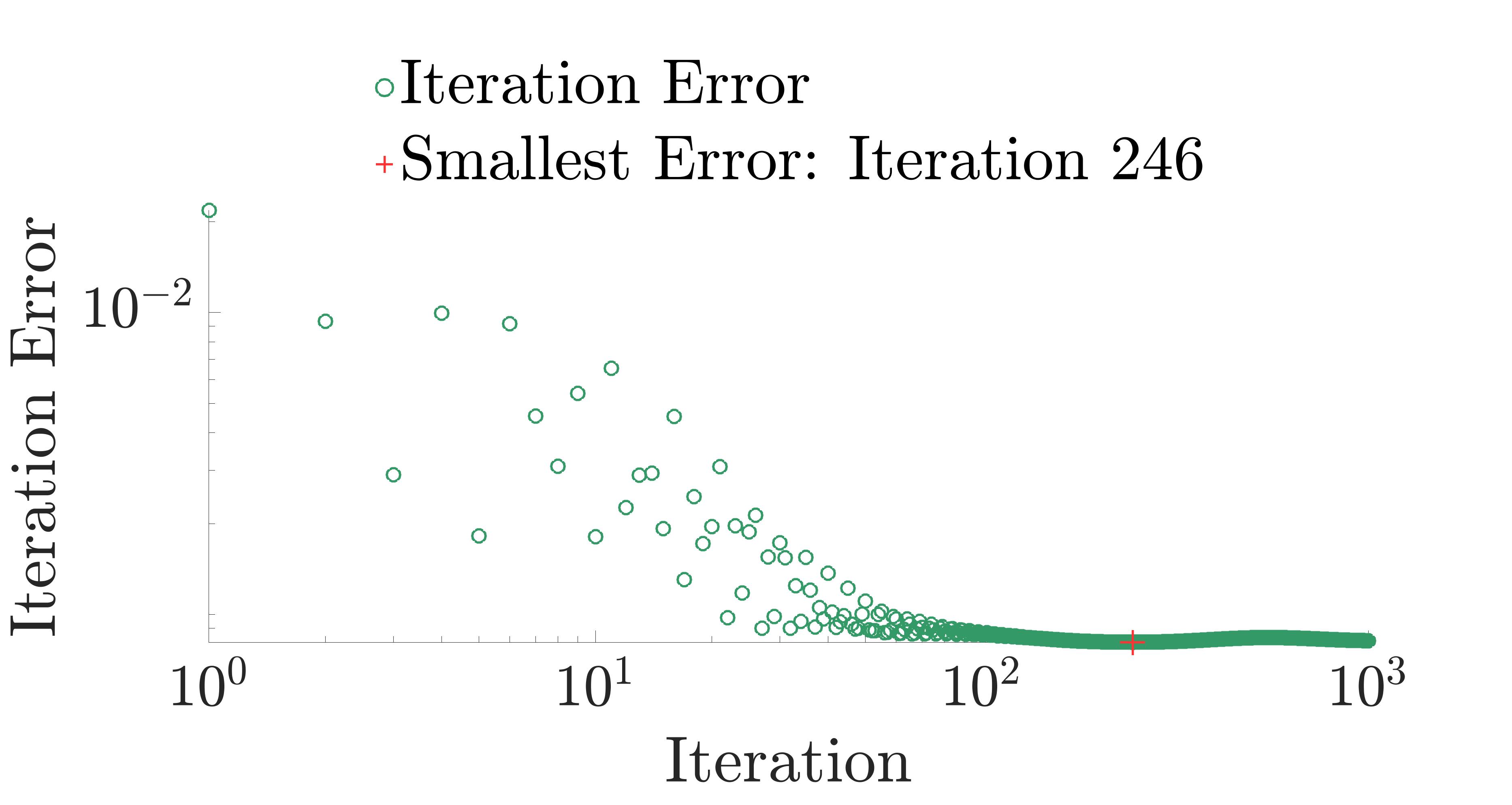}
    \end{subfigure}
    \caption{Results for Case \eqref{case_swe_5c}. 
    Left: plots of the true $p$ and the numerically recovered $p$ at iteration 246 corresponding with the smallest residue; Middle: plots of the true $p$, the noisy initial guess, and various iteration values for $p$; Right: iteration errors on a log-log scale.}
    \label{fig:swe_5C}
\end{figure}

\hlc{In this subsection, we demonstrate the necessity of the regularization terms} discussed in Section \ref{sec_regularization}. First, we demonstrate what happens when one or more of the regularization terms are completely removed from the algorithm. Table \ref{tab:swe_p_5} provides the regularization parameters used for the three test cases examined which are visualized in Figures \ref{fig:swe_5A} - \ref{fig:swe_5C}.

\begin{table}[h]
    \centering
    \begin{tabular}{c c c}
    \toprule
        Case & $\gamma_L$ & $\gamma_H$ \\
        \midrule
        \stepcounter{subsecval}\refstepcounter{subsecvalnext}(\thesubsecvalnext)\label{case_swe_5a} &  $0$ & $5 \times 10^{-8}$\\
        \refstepcounter{subsecvalnext}(\thesubsecvalnext)\label{case_swe_5b} &  $1 \times 10^{-6}$ & $0$\\
        \refstepcounter{subsecvalnext}(\thesubsecvalnext)\label{case_swe_5c} &  $0$ & $0$\\
    \bottomrule
    \end{tabular}
    \caption{The regularization parameters tested for Case \eqref{case_swe_2d}.}
    \label{tab:swe_p_5}
\end{table}

Results for Case \eqref{case_swe_5a}, representing the case with no $L^1$ regularization, are displayed in Figure \ref{fig:swe_5A}. While in this case, two peaks of differing heights are recovered, the recovered peak heights are lower than the true peak locations and occur slightly earlier in time. Furthermore, the iteration with the smallest error (iteration 1000) occurs in Case \eqref{case_swe_5a} much later than in Case \eqref{case_swe_2d} (iteration 151).
Case \eqref{case_swe_5b} represents the case with no $H^1$ regularization. The corresponding results can be found in Figure \ref{fig:swe_5B}. Not only does the lack of $H^1$ regularization result in a noisy solution, but also the two peaks are not even recovered after 1000 iterations, demonstrating slowed convergence of the iterative scheme.
The situation in which neither $L^1$ nor $H^1$ regularization is implemented, Case \eqref{case_swe_5c}, is shown in Figure \ref{fig:swe_5C}. Results in this case are nearly identical to those found in Case \eqref{case_swe_5b} and Figure \ref{fig:swe_5B}.

\begin{figure}[h]
    \centering
    \includegraphics[width=0.5\linewidth]{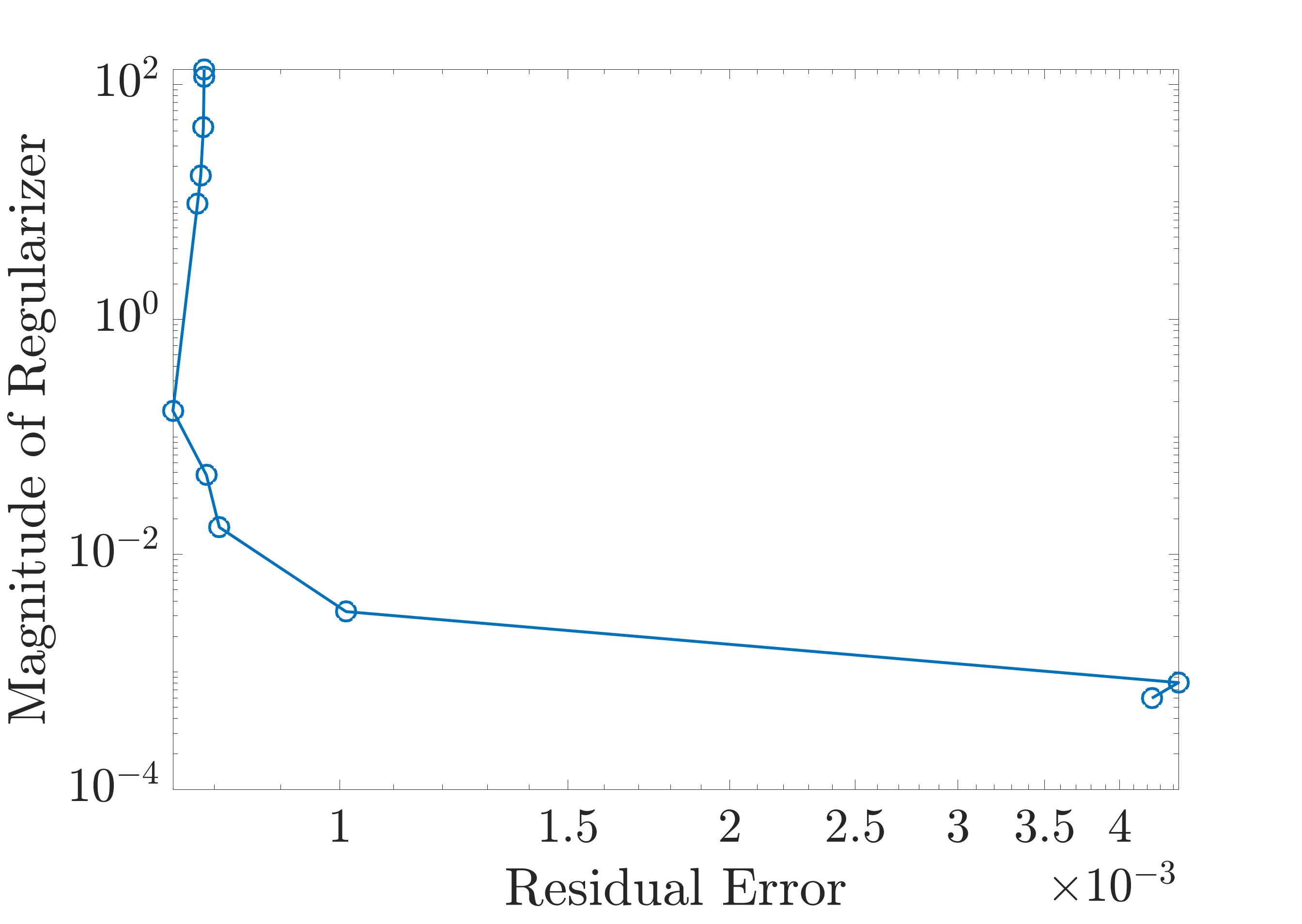}
    \caption{L-curve corresponding to Case \eqref{case_swe_2d} for regularization parameter $\hat{\gamma} = 10^i$, $i \in [-5,5]$.}
    \label{fig:lcurve}
\end{figure}

The second way we demonstrate the impact of regularization is by way of a modification to the minimization problem presented in \eqref{min_prob}. The modification involves an additional parameter and takes on the following form
\begin{equation} \label{min_prob_new}
    \begin{split}
        & \text{minimize } \J(p) := \J_0(p) + \hat{\gamma}\R(p) \quad \text{subject to } \eqref{swe_jacobiansys}.
    \end{split}
\end{equation}
This formulation allows us to perform an L-curve test \cite{Hansen00thel-curve} when running an example for multiple values of $\hat{\gamma}$. The L-curve test is a visualization tool used for examining the impact of regularization and to find a balance with the residual errors. Here we consider the same problem setup as in Case \eqref{case_swe_2d} and set the regularization parameter to $\hat{\gamma} = 10^i$ for the integer $i$ ranging between $-5$ and $5$. Note that the choice for $\gamma_L$ and $\gamma_H$ stay fixed at $1\times 10^{-6}$ and $5\times 10^{-8}$, respectively, for each choice of $\hat{\gamma}$. Figure \ref{fig:lcurve} shows the smallest residual error over all iterations on the x-axis and the magnitude of the regularizer from the same iteration on the y-axis, both on a log-scale. The values of $\hat{\gamma}$ occurring near the `elbow' of the L-curve correspond to $\hat{\gamma} = 1, 10, 100, 1000$, and are reasonable choices of regularization parameter to use.

\subsection{Recovering \texorpdfstring{$p(t)$}{p(t)} with Spatial Discontinuities in the Solutions of the Forward Problem}
In the previous two subsections, the final stopping time is chosen to ensure that the solutions are smooth over the entire computational region. One well-known fact of hyperbolic conservation laws is that discontinuities may appear even when the initial conditions are smooth. In this subsection, we explore the capability of the proposed inverse algorithms when the solutions contain spatial discontinuities.  
We use the same initial conditions as in \eqref{ics_swe_smooth0} and spatial components of the bottom topography functions as in \eqref{bvals_swe_smooth0}. We run the forward problem for a longer time so that discontinuities form in the solutions. Here, we set the final time as $T=0.2$. We consider two choices for $p_{true}(t)$ in this subsection, which are listed in Table \ref{tab:swe_p_4} along with their corresponding initial guesses.

\begin{table}[h]
    \centering
    \begin{tabular}{c c c}
    \toprule
        Case & $p_{true}(t)$ & $p^0(t)$ \\
        \midrule
        \stepcounter{subsecval}\refstepcounter{subsecvalnext}(\thesubsecvalnext)\label{case_swe_4a} 
        & $e^{\beta\left(t-\frac{1}{2}T\right)^2}+1$ 
        & $1$ \\
        \refstepcounter{subsecvalnext}(\thesubsecvalnext)\label{case_swe_4b}  
        & $e^{4\beta\left(t-\frac{1}{4}T\right)^2}
        +\frac{3}{2}e^{4\beta\left(t-\frac{1}{2}T\right)^2}
        -\frac{1}{2}e^{4\beta\left(t-\frac{3}{4}T\right)^2}+1$
        & $1$\\
    \bottomrule
    \end{tabular}
    \caption{The true function for $p(t)$ denoted as $p_{true}$ and the corresponding initial guess used, $p^0$ with $\beta = -700$. Multiplicative noise is applied to $p^0$ in the simulations.}
    \label{tab:swe_p_4}
\end{table}

To capture the discontinuities well and remove the possible oscillations, a slope limiter is often employed in the DG method. We implement two different slope limiters for generating the measured data, as well as for solving the forward problem in the inverse scheme. The simple minmod limiter \cite{Cockburn1989TVBFramework} is employed along the characteristic direction to generate the measured data. On the other hand, the WENO limiter, introduced by Qiu and Shu in \cite{Qiu2005Runge-KuttaLimiters} is used for the forward solver within the iterative inverse scheme. This limiter is known to be robust and it is able to capture the sharp transition of the discontinuities. The implementation of two different slope limiters was in an effort to avoid `inverse crime'.

\begin{figure}
\begin{subfigure}[t]{0.32\textwidth}
    \centering
    \includegraphics[width=\linewidth]{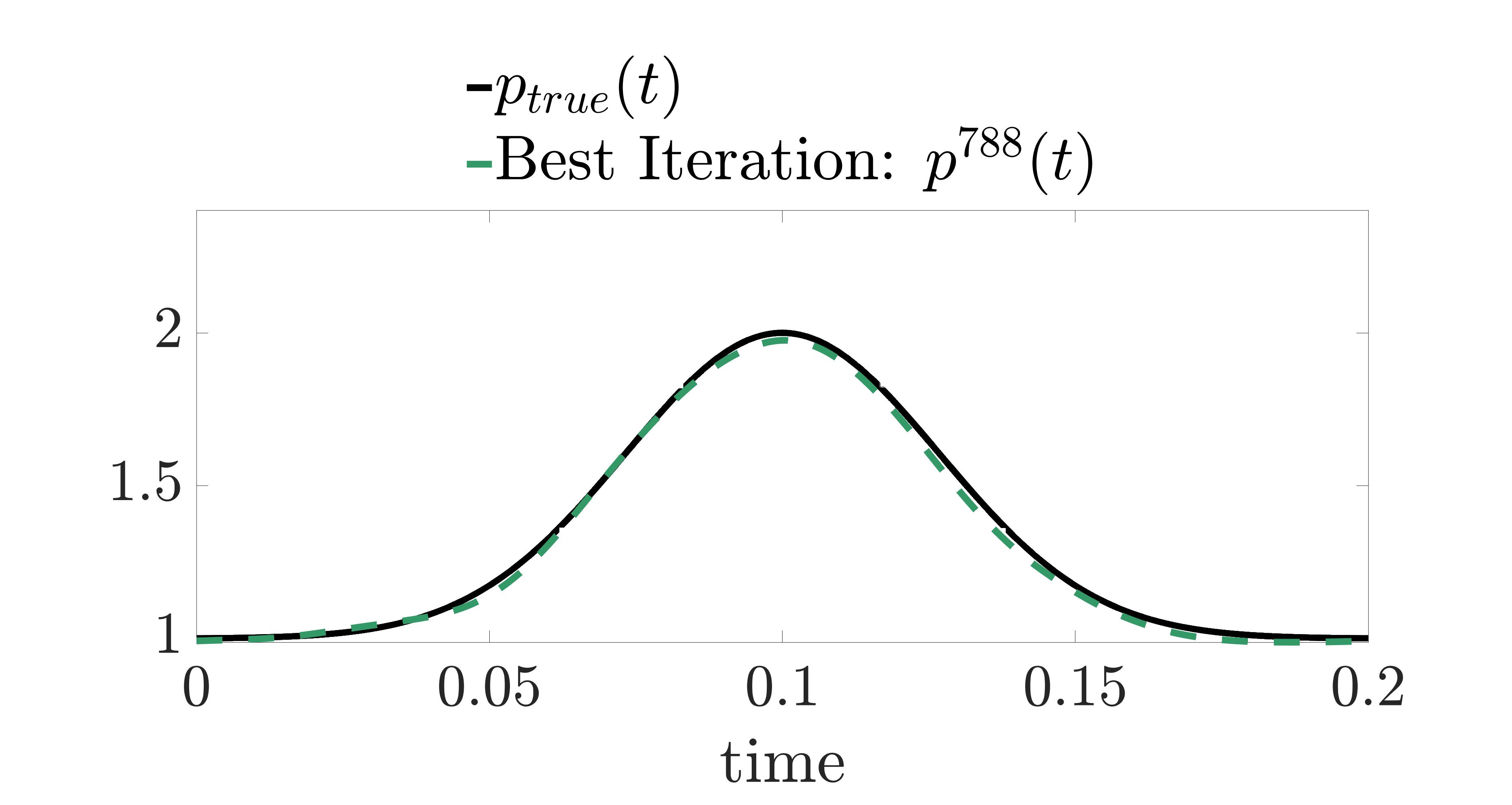}
\end{subfigure}
\begin{subfigure}[t]{0.32\textwidth}
    \centering
    \includegraphics[width=\linewidth]{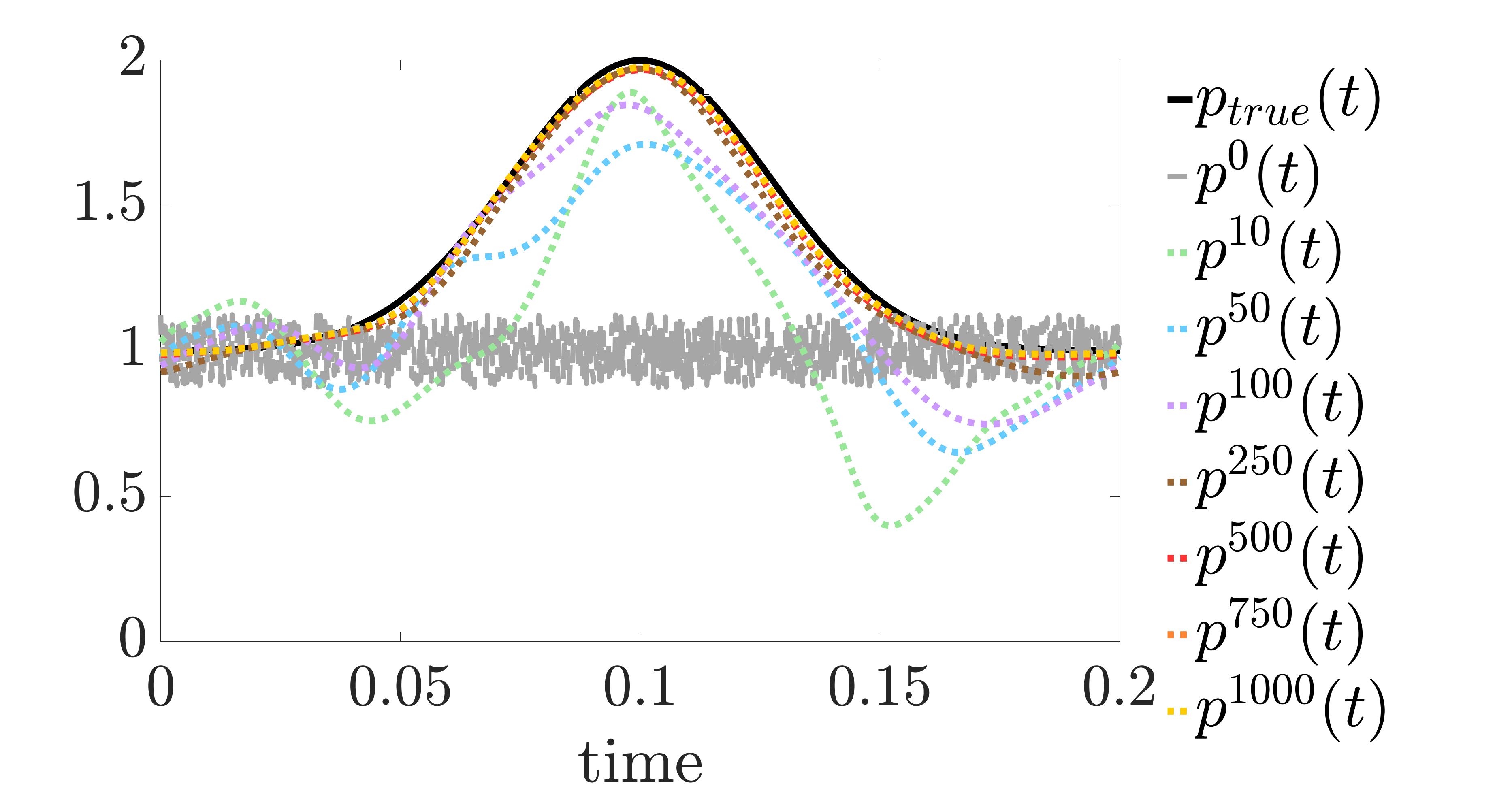}
\end{subfigure}
\begin{subfigure}[t]{0.32\textwidth}
    \centering
    \includegraphics[width=\linewidth]{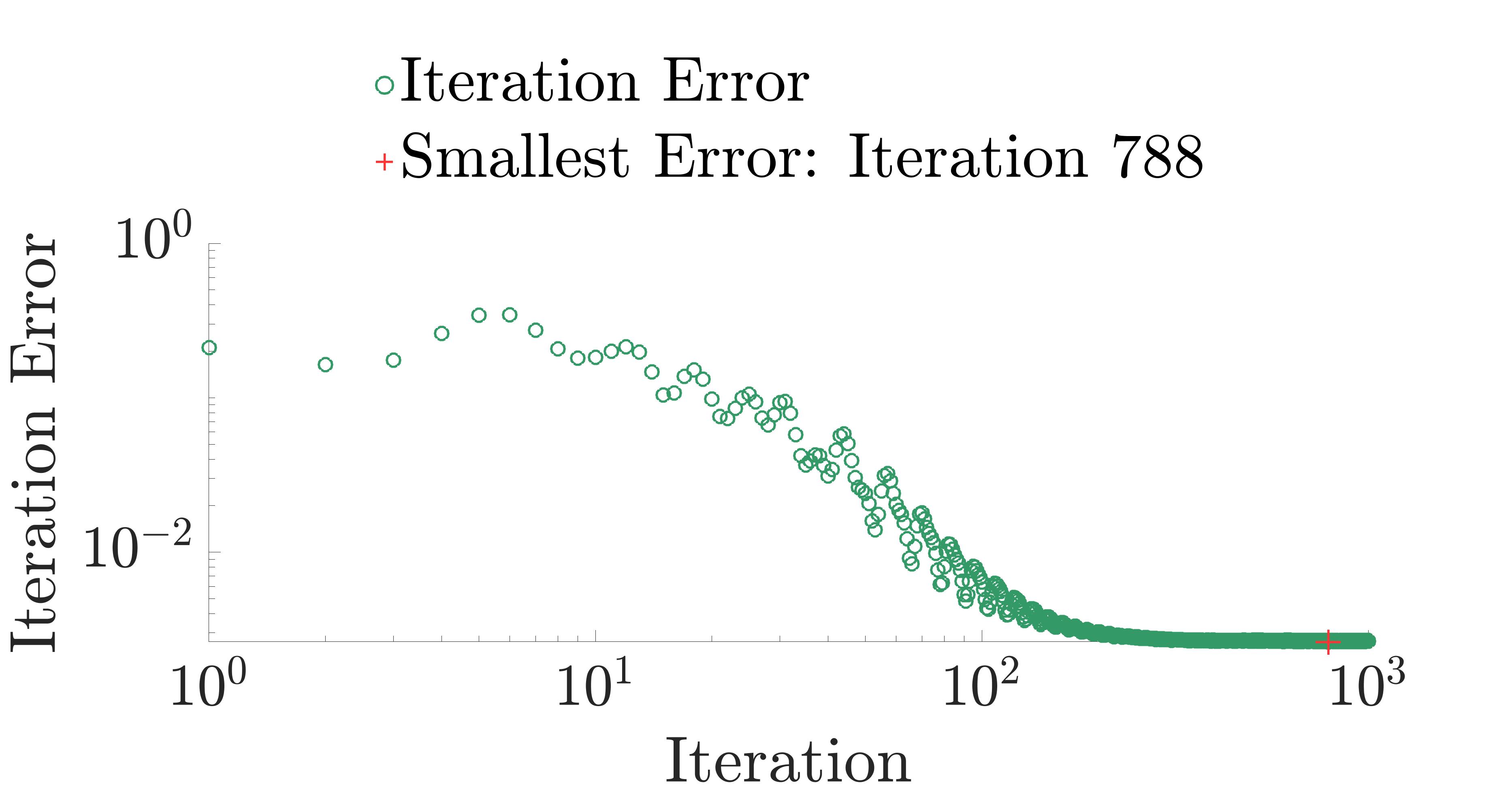}
\end{subfigure}
\caption{Results for Case \eqref{case_swe_4a}. Left: plots of the true $p$ and the $p$ corresponding with the smallest residue error, at iteration 788; Middle: plots of the true $p$, the noisy initial guess, and various iteration values for $p$; Right: iteration errors on a log-log scale.}
    \label{fig:swe_4A}
\end{figure}

\begin{figure}
    \centering
    \includegraphics[width=\linewidth]{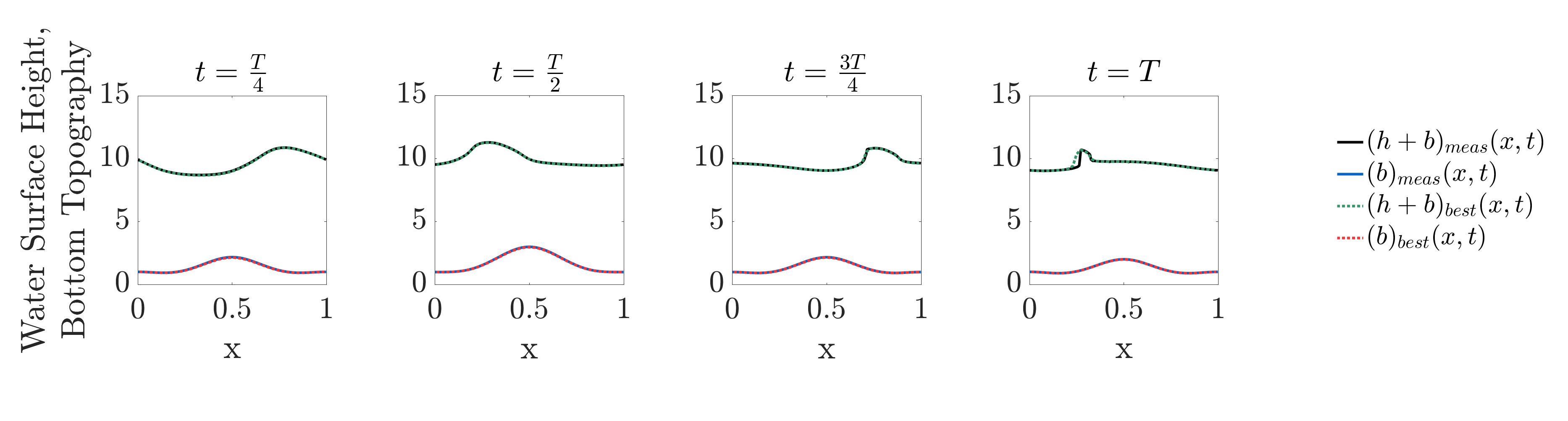}
    \includegraphics[width=\linewidth]{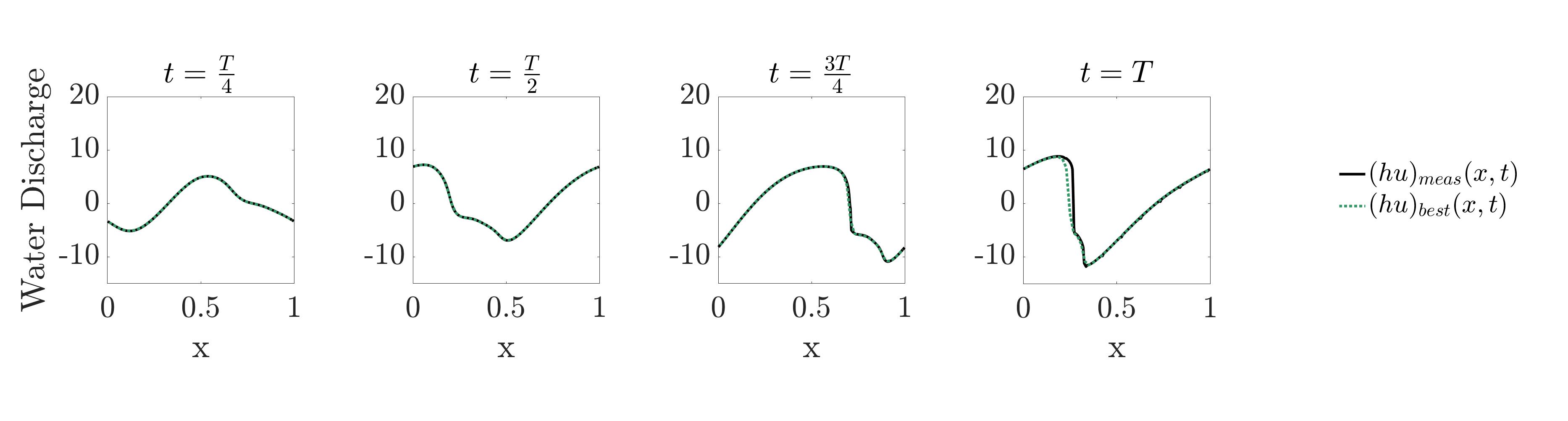}
    \caption{Comparison between the measured forward solutions and the numerical results from the iteration with the smallest residue error for Case \eqref{case_swe_4a}. The results shown are for 4 different time snapshots. The bottom topography function, $b$, water surface height, $h+b$ (top row), and the water discharge, $hu$ (bottom row), are compared.}
    \label{fig:swe_4A_2}
\end{figure}

\begin{figure}
\begin{subfigure}[t]{0.32\textwidth}
    \centering
    \includegraphics[width=\linewidth]{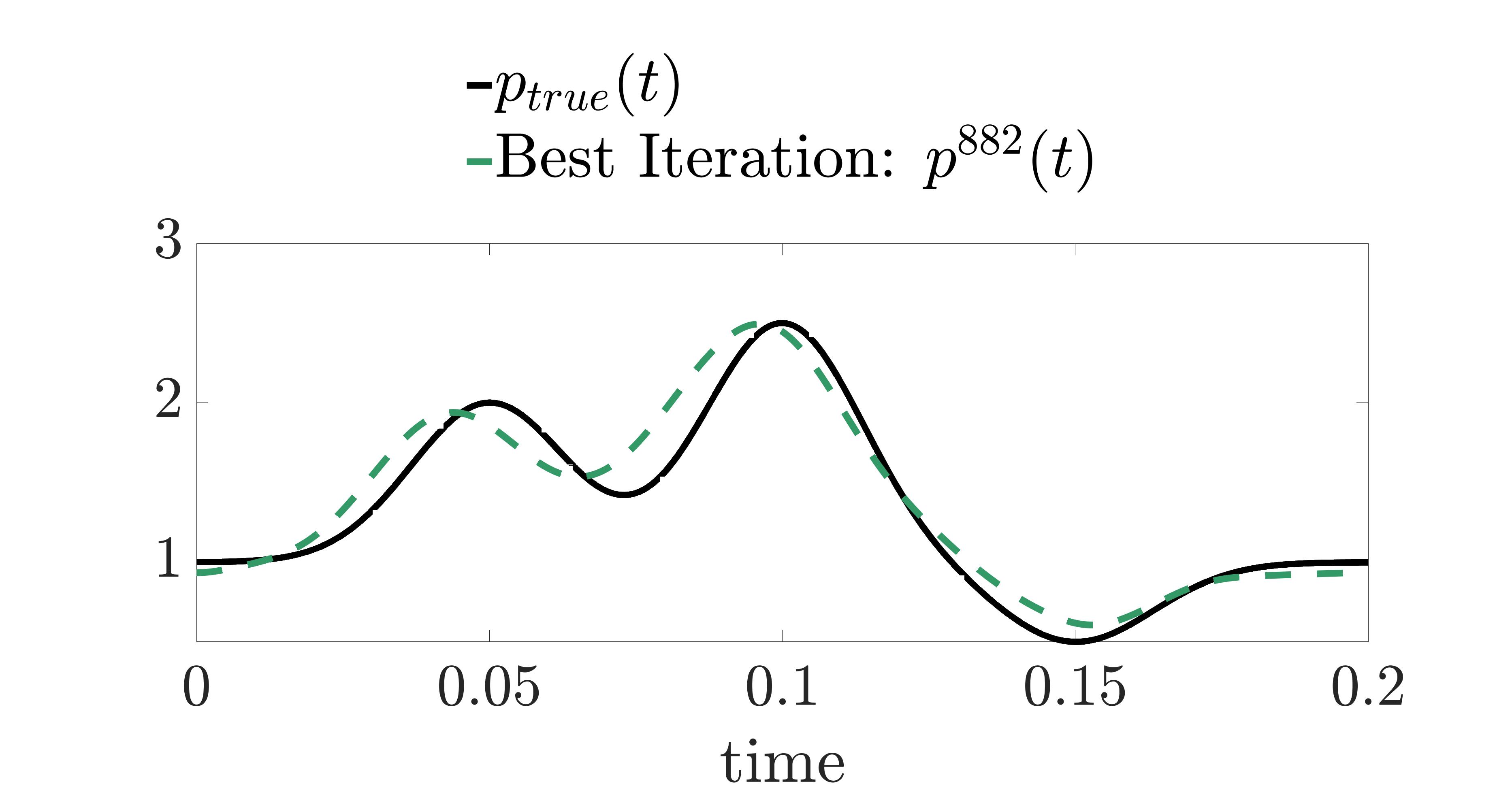}
\end{subfigure}
\begin{subfigure}[t]{0.32\textwidth}
    \centering
    \includegraphics[width=\linewidth]{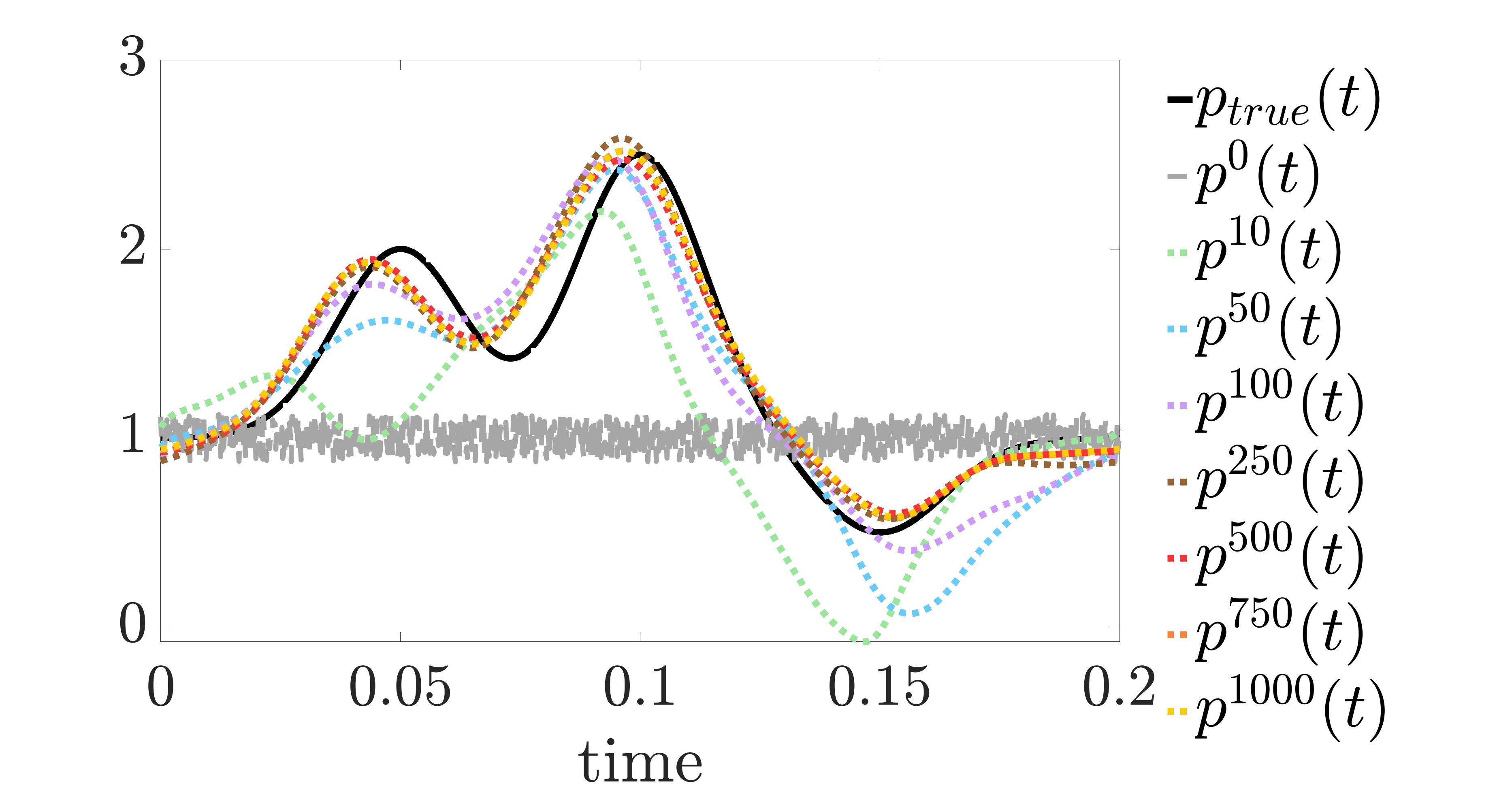}
\end{subfigure}
\begin{subfigure}[t]{0.32\textwidth}
    \centering
    \includegraphics[width=\linewidth]{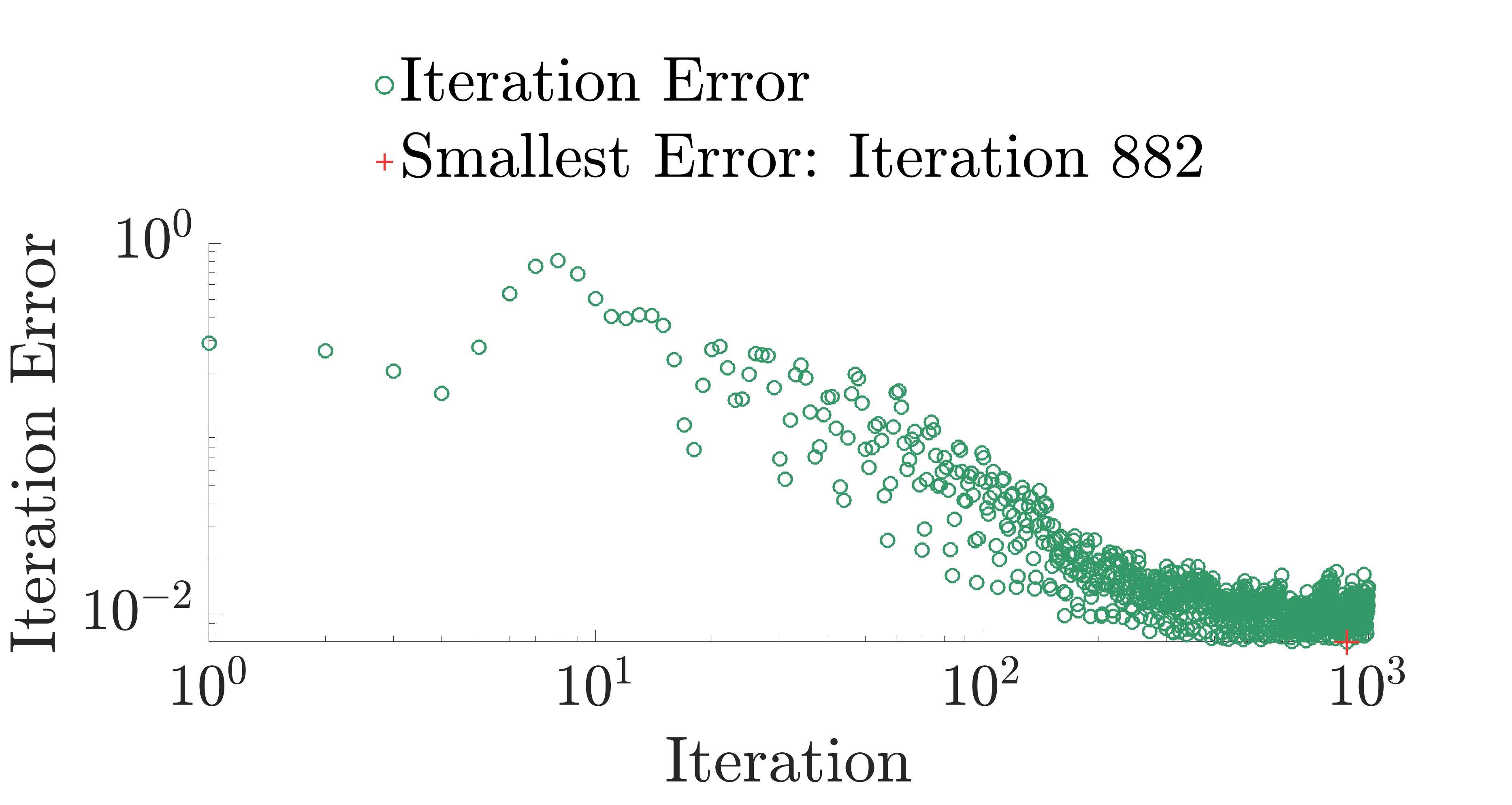}
\end{subfigure}
\caption{Results for Case \eqref{case_swe_4b}. Left: plots of the true $p$ and the $p$ corresponding with the smallest residue error, at iteration 882; Middle: plots of the true $p$, the noisy initial guess, and various iteration values for $p$; Right: iteration errors on a log-log scale.}
    \label{fig:swe_4B}
\end{figure}

\begin{figure}
    \centering
    \includegraphics[width=\linewidth]{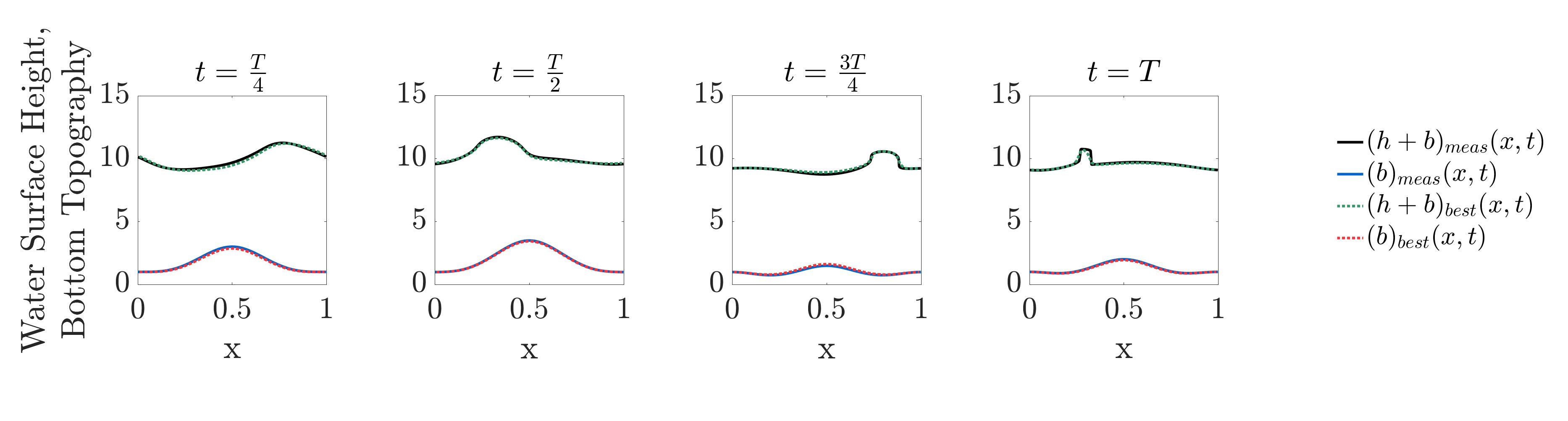}
    \includegraphics[width=\linewidth]{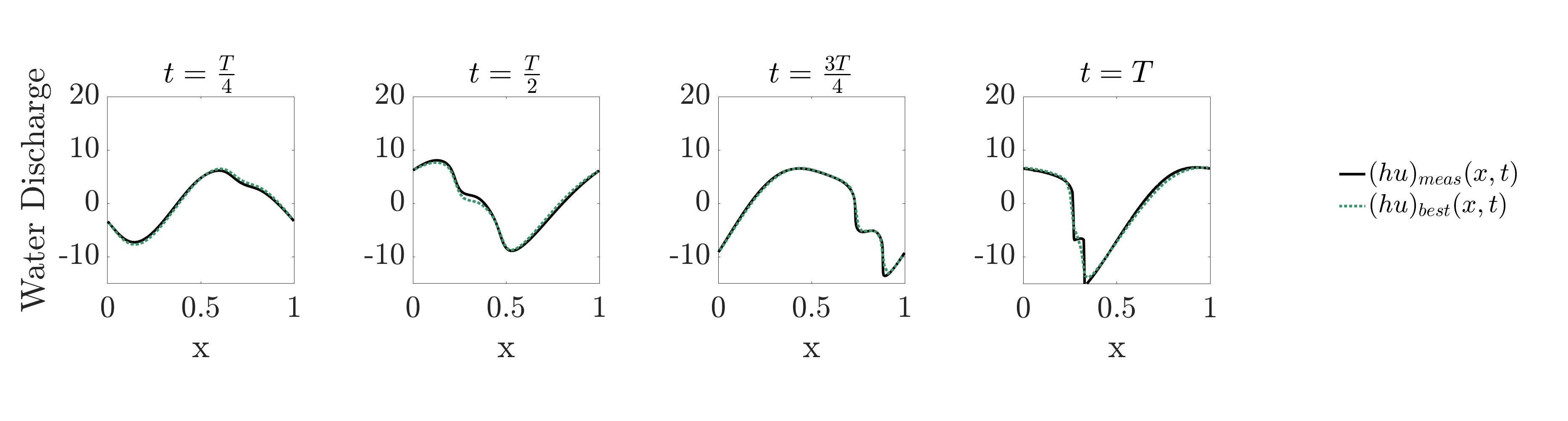}
    \caption{Comparison between the measured forward solutions and the numerical results from the iteration with the smallest residue error for Case \eqref{case_swe_4b}. The results shown are for 4 different time snapshots. The bottom topography function, $b$, water surface height, $h+b$ (top row), and the water discharge, $hu$ (bottom row), are compared.}
    \label{fig:swe_4B_2}
\end{figure}

A learning rate of $\ell=0.02$ was implemented with the regularization parameters $\gamma_L = 1 \times 10^{-4}$ and $\gamma_H = 1 \times 10^{-6}$. Results for case \eqref{case_swe_4a} are found in Figures \ref{fig:swe_4A} and \ref{fig:swe_4A_2} while Figures \ref{fig:swe_4B} and \ref{fig:swe_4B_2} contain the results for test \eqref{case_swe_4b}. Note that periodic boundary conditions are employed, therefore the discontinuities can pass the right boundary and re-enter the domain through the left boundary at some time between $3T/4$ and $T$, which means the measured data include the information of discontinuities. In both cases, we observe that our algorithm can recover the exact function $p_{true}(t)$ well, and the results are comparable with those containing smooth data only. We have included the comparison between the measured data of the water surface height, bottom topography function, water discharge, and the corresponding numerical solutions at the iteration with the smallest residue, at different times $t=\frac{T}{4}, \frac{T}{2}, \frac{3T}{4}$, and $T$, from which we can observe the numerical solutions match the measured data well. This elucidates the numerical scheme developed for the inverse problems can recover the true $p(t)$ well even when discontinuities develop in the solutions of the forward problem and hence in the measured data.    

\section{Conclusion}\label{sec_conclusion}

In this paper we constructed and validated an adjoint-based approach for recovering the bottom topography function in the source term of the one-dimensional SWEs, from the noisy measurement data at two boundaries of the domain. 
One novelty of this work is that the reconstruction of the bottom topography function is accomplished with only boundary data from a single measurement event. 
The adjoint scheme was determined by a linearization of the forward system, and has been derived for general hyperbolic balance laws. 
Another contribution of this work is the inclusion of two regularization terms. These extra regularization terms in the numerical approach aided in convexifying and handling the ill-posedness of the problem. The bottom topography function was recovered through an iterative process using a three-operator splitting descent method. Extensive numerical tests were carried out, which demonstrated that a variety of shapes for the true $p(t)$ function could be recovered regardless of the noisy initial guess. 
As a followup, we would like to develop an inverse algorithm to recover the general bottom topography $B(x,t)$ and its extension in higher dimensional problems. 
It is also worthwhile to extend this framework to study the inverse problems associated with other hyperbolic balance laws.